\documentclass[twoside]{amsart}

\usepackage[a4paper,top=2.54cm,bottom=2.54cm,left=3.17cm,right=3.17cm,%
            includehead,includefoot]{geometry}

\usepackage{amsmath,amssymb,amsfonts,amsthm}
\usepackage{graphicx}
\usepackage{subfigure}
\usepackage{float}
\usepackage{comment}
\usepackage{xcolor}
\usepackage[numbers,square,sort&compress]{natbib}
\usepackage{hyperref}
\usepackage{mathrsfs}  \hypersetup{colorlinks,citecolor=blue,linkcolor=blue,urlcolor=black,breaklinks=true}
\usepackage{booktabs}
\usepackage{colortbl}
\usepackage{caption}
\usepackage{enumitem}
\usepackage{epstopdf}
\usepackage{algorithm}
\usepackage{algpseudocode}
\usepackage{array}
\usepackage{bbding}
\usepackage{fancyhdr} 
\usepackage{fancyvrb} 
\usepackage{longtable}
\usepackage{listings}
\usepackage{rotating,rotfloat} 
\usepackage{yhmath} 
\usepackage{amsmath}
\numberwithin{equation}{section}
\usepackage{indentfirst}
\usepackage{fancyhdr}
\usepackage{titlesec}
\usepackage{graphicx}

\titleformat{\section}
{\normalfont\Large\bfseries\MakeUppercase}
{\thesection}{1em}{}

\titleformat{\subsection}
{\normalfont\large\bfseries}
{\thesubsection}{1em}{}

\theoremstyle{definition}
\newtheorem{theorem}{Theorem}[section]
\newtheorem{lemma}[theorem]{Lemma}
\newtheorem{corollary}[theorem]{Corollary}
\newtheorem{proposition}[theorem]{Proposition}

\newtheorem{remark}[theorem]{Remark}

\numberwithin{equation}{section}

\newcommand{\bbm}{\begin{bmatrix}}
\newcommand{\ebm}{\end{bmatrix}}

\title{A Double Centralizer Theorem in Quantized Weyl Superalgebras
and a First Fundamental Theorem for Quantum General Linear Supergroups}
\author{Qianqi Wang}
\address{
School of Mathematical Sciences,
East China Normal University,
Shanghai 200241, China
}
\email{51275500030@stu.ecnu.edu.cn}
\date{September 2026}

\begin{document}
\date{}

\maketitle

\begin{abstract}

We construct a quantized Weyl superalgebra $\mathscr{PD}_{r|s}^{k|l}$ as a deformed twisted tensor product of the quantized coordinate superalgebra $\mathscr{P}_{r|s}^{k|l}$ and the superalgebra $\mathscr{D}_{r|s}^{k|l}$ of constant-coefficient quantum differential operators. We establish a PBW basis for $\mathscr{PD}_{r|s}^{k|l}$ and show that it acts faithfully on $\mathscr{P}_{r|s}^{k|l}$. Let $\mathscr{L}$ and $\mathscr{R}$ be the images of $\mathbf{U}_q(\mathfrak{gl}_{k|l})$ and $\mathbf{U}_q(\mathfrak{gl}_{r|s})$ under their natural actions on $\mathscr{P}_{r|s}^{k|l}$. We prove that their intersections with $\mathscr{PD}_{r|s}^{k|l}$ are mutual supercentralizers in $\mathscr{PD}_{r|s}^{k|l}$, identify these intersections with the corresponding invariant subsuperalgebras, and show that they are generated by explicit quantum polarization operators. This establishes an operator commutant version of the first fundamental theorem of invariant theory for quantum general linear supergroups. More generally, we establish a first fundamental theorem for the action of $\mathbf{U}_q(\mathfrak{gl}_{r|s})$ on subsuperalgebras $\mathscr{A}_{p|t,p'|t',r|s}$ of $\mathscr{PD}_{r|s}^{k|l}$ by determining explicit generators for their invariant subsuperalgebras.
\end{abstract}

\section{Introduction}

The first fundamental theorem (FFT) of invariant theory
is one of the foundational results in classical invariant
theory, describing generators for algebras of polynomial
invariants under classical group actions. Its classical formulation for general linear groups was
systematically developed by Weyl
\cite{W39}.
With the development of quantum groups, considerable
attention has been devoted to quantum analogues of
classical invariant theory.
In particular, Lehrer, H.~Zhang and R.~B.~Zhang \cite{LZZ11}
constructed flat deformations of classical polynomial
algebras as module algebras for
$U_q(\mathfrak{gl}_n)$, $U_q(\mathfrak{o}_N)$
and $U_q(\mathfrak{sp}_{2n})$,
using twisted tensor products defined by the
universal $R$-matrix. They proved that the corresponding invariant
subalgebras are generated by explicit quadratic
invariants and determined their commutation relations,
thereby establishing noncommutative analogues
of the classical FFT.

The study of invariant theory has also been extended
to Lie superalgebras and quantum supergroups. Lehrer, H.~Zhang and R.~B.~Zhang \cite{LZZ20}
proved that the ribbon graph functors associated with
$U_q(\mathfrak{gl}_{m|n})$ and the modified
$U_q(\mathfrak{osp}_{2\ell+1|2n})$ are full,
thereby establishing a tensor categorical formulation
of the FFT. Subsequently, Y.~Zhang \cite{Zh20} established the FFT and SFT
for $U_q(\mathfrak{gl}_{m|n})$ in a noncommutative
polynomial setting.
Using quantum super Howe duality, he determined
explicit generators of the invariant subalgebras
and realised these subalgebras as quotients of
braided supersymmetric algebras.

Howe \cite{Ho89} gave a uniform formulation of the FFT
for classical groups using Weyl and Clifford algebras,
describing invariant operator algebras in terms of
classical reductive dual pairs. For the pair $(\mathfrak{gl}_m,\mathfrak{gl}_n)$, this yields a double centralizer property in the
algebra of polynomial-coefficient differential operators. More precisely, the images of the universal enveloping algebras $U(\mathfrak{gl}_m)$ and $U(\mathfrak{gl}_n)$ under their natural actions on polynomial functions
on $\operatorname{Mat}_{m\times n}$ are mutual centralizers in this differential operator algebra. The equivalence between this operator commutant
formulation and the polynomial formulation of the FFT
is explained in \cite{Ho95}.

In \cite{LSS24}, Letzter, Sahi and Salmasian provide a quantum analogue
of the operator commutant formulation of
the classical FFT. More precisely, they
constructed quantized Weyl algebras using
deformed twisted tensor products, established a double centralizer property
for the images of two quantum general linear
algebras inside the quantized Weyl algebra
and determined explicit generators for the
corresponding centralizer subalgebras.

These developments naturally lead to the question
of whether a similar operator commutant formulation
exists for quantum general linear supergroups.
The purpose of our paper is to address
this question by constructing quantized Weyl
superalgebras and establishing the corresponding
double centralizer property.

We now describe our construction and main results. Let $\mathbb{K}=\mathbb{C}(q)$, where $q$ is an
indeterminate, and let $k,l,r,s$ be nonnegative integers.
For nonnegative integers $a,b$, set
$
\mathbf{I}_{a|b}:=\{1,\ldots,a+b\}
$, with parity $[i]=0$ for $1\leq i\leq a$ and
$[i]=1$ for $a<i\leq a+b$.
Following the construction of Letzter, Sahi and
Salmasian \cite{LSS24}, we define the quantized Weyl
superalgebra $\mathscr{PD}_{r|s}^{k|l}$ as a deformed
twisted tensor product of the quantum coordinate
superalgebra $\mathscr{P}_{r|s}^{k|l}$ of
$(k|l)\times(r|s)$ matrices and the superalgebra
$\mathscr{D}_{r|s}^{k|l}$ of constant-coefficient
quantum differential operators.
The twisting uses the universal $R$-matrices
associated with both
$\mathbf{U}_q(\mathfrak{gl}(k|l))$ and
$\mathbf{U}_q(\mathfrak{gl}(r|s))$.
The superalgebra $\mathscr{PD}_{r|s}^{k|l}$ is generated by
$
t_{i,c},\ \partial_{i,c},
\ i\in\mathbf{I}_{k|l},\
c\in\mathbf{I}_{r|s},
$
where $t_{i,c}$ and $\partial_{i,c}$ are the
standard generators of
$\mathscr{P}_{r|s}^{k|l}$ and
$\mathscr{D}_{r|s}^{k|l}$, respectively, and both
have parity $[i]+[c]$.
The defining relations are given in Section \ref{subsection: PD}.

These superalgebras carry compatible
$\mathbf{U}_q(\mathfrak{gl}(k|l))
\otimes\mathbf{U}_q(\mathfrak{gl}(r|s))$-module
superalgebra structures.
Moreover, $\mathscr{PD}_{r|s}^{k|l}$ acts naturally
and faithfully on $\mathscr{P}_{r|s}^{k|l}$.
Thus we obtain homomorphisms of superalgebras
\[
\begin{aligned}
\phi_U:\,
\mathbf{U}_q(\mathfrak{gl}(k|l))
\otimes\mathbf{U}_q(\mathfrak{gl}(r|s))
&\longrightarrow
\operatorname{End}_{\mathbb{K}}
\bigl(\mathscr{P}_{r|s}^{k|l}\bigr),\\
\phi_{PD}:\,
\mathscr{PD}_{r|s}^{k|l}
&\hookrightarrow
\operatorname{End}_{\mathbb{K}}
\bigl(\mathscr{P}_{r|s}^{k|l}\bigr).
\end{aligned}
\]
We identify $\mathscr{PD}_{r|s}^{k|l}$ with its
image under $\phi_{PD}$ and set
\[
\begin{aligned}
\mathscr{L}
&:=
\phi_U\bigl(
\mathbf{U}_q(\mathfrak{gl}(k|l))\otimes 1
\bigr),
&
\mathscr{R}
&:=
\phi_U\bigl(
1\otimes\mathbf{U}_q(\mathfrak{gl}(r|s))
\bigr),\\
\mathscr{L}_{\bullet}
&:=
\mathscr{L}\cap\mathscr{PD}_{r|s}^{k|l},
&
\mathscr{R}_{\bullet}
&:=
\mathscr{R}\cap\mathscr{PD}_{r|s}^{k|l}.
\end{aligned}
\]
These intersections consist of the operators in
the two quantum supergroup images that belong
to the quantized Weyl superalgebra.
Define the quantum polarization operators
\[
\begin{aligned}
L_{i,j}
&:=
\sum_{c=1}^{r+s}
(-1)^{[c]([i]+1)}
t_{i,c}\partial_{j,c},
&& i,j\in\mathbf{I}_{k|l},\\
R_{i,j}
&:=
\sum_{c=1}^{k+l}
(-1)^{[c]([j]+1)}
t_{c,i}\partial_{c,j},
&& i,j\in\mathbf{I}_{r|s}.
\end{aligned}
\]
These can be viewed as super analogues of the quantum
polarization operators introduced in \cite{LSS24},
the parity signs disappear when $l=s=0$,
and the formulas reduce to those of \cite{LSS24}.
Each operator $L_{i,j}$ or $R_{i,j}$ has parity
$[i]+[j]$.
We denote by $\langle L_{i,j}\rangle$ and
$\langle R_{i,j}\rangle$ the unital subsuperalgebras
generated by the respective families.

Let $\epsilon_L$ and $\epsilon_R$ be the counits of
$\mathbf{U}_q(\mathfrak{gl}(k|l))$ and
$\mathbf{U}_q(\mathfrak{gl}(r|s))$, respectively.
Write
\[
\begin{aligned}
(\mathscr{PD}_{r|s}^{k|l})_{(\epsilon_L)}
:=
\bigl\{
f\in\mathscr{PD}_{r|s}^{k|l}\ \big|\
&(x\otimes 1)\cdot f=\epsilon_L(x)f
\ \text{for all }
x\in\mathbf{U}_q(\mathfrak{gl}(k|l))
\bigr\},
\\[1ex]
(\mathscr{PD}_{r|s}^{k|l})_{(\epsilon_R)}
:=
\bigl\{
f\in\mathscr{PD}_{r|s}^{k|l}\ \big|\
&(1\otimes y)\cdot f=\epsilon_R(y)f
\ \text{for all }
y\in\mathbf{U}_q(\mathfrak{gl}(r|s))
\bigr\},
\end{aligned}
\]
where $\cdot$ denotes the module superalgebra action
constructed above.

For subsuperalgebras $Y,Z$ of a superalgebra $X$,
let $Y^Z$ denote the supercentralizer of $Z$ in $Y$:
\[
Y^Z
:=
\bigoplus_{a\in\mathbb{Z}_2}
\left\{
y\in Y_a\ \middle|\
yz=(-1)^{ab}zy
\text{ for all }b\in\mathbb{Z}_2
\text{ and }z\in Z_b
\right\},
\]
where $Y_a$ and $Z_b$ are the corresponding
homogeneous components.
Our first main result is the following, which establishes an operator commutant
version of the FFT for quantum general linear
supergroups. 

\noindent\textbf{Main Theorem 1.} 
With the notation above, we have
\begin{equation*}
    \begin{aligned}
        (\mathscr{PD}_{r|s}^{k|l})^{\mathscr{R}_{\bullet}}
        &=
        (\mathscr{PD}_{r|s}^{k|l})_{(\epsilon_R)}
        =
        \mathscr{L}_{\bullet}
        =
        \langle L_{i,j}\rangle,\\
        (\mathscr{PD}_{r|s}^{k|l})^{\mathscr{L}_{\bullet}}
        &=
        (\mathscr{PD}_{r|s}^{k|l})_{(\epsilon_L)}
        =
        \mathscr{R}_{\bullet}
        =
        \langle R_{i,j}\rangle.
    \end{aligned}
\end{equation*}
In \cite{Zh20}, Zhang constructed the braided tensor product
$
\mathcal P_{r|s}^{k|l}
=
\mathcal M_{m|n}^{k|l}
\otimes_R
\overline{\mathcal M}_{m|n}^{r|s}
$
and determined generators of its
$U_q(\mathfrak{gl}(m|n))$-invariant subalgebra.
The defining relations of this algebra are homogeneous,
and its classical limit is the supersymmetric algebra
$
S(
\mathbb C^{k|l}\otimes\mathbb C^{m|n}
\oplus
\mathbb C^{r|s}\otimes(\mathbb C^{m|n})^*).
$
In contrast, the defining relations in our quantized Weyl
superalgebra $\mathscr{PD}_{r|s}^{k|l}$ contain constant terms.
A suitable integral form of $\mathscr{PD}_{r|s}^{k|l}$ specializes
at $q=1$ to an algebra isomorphic to the Weyl-Clifford superalgebra (see Section \ref{subsection: Integrals and specializations}). From this point of view, our construction gives a quantum deformation of the
Weyl-Clifford superalgebra, rather than
of a supersymmetric algebra. 

The second main result (Theorem \ref{theorem: AeR is generated by  Lij}) is a FFT for the quantum general linear supergroup
$U_q(\mathfrak{gl}(r|s))$ acting on the subalgebras $\mathscr A_{p|t,p'|t',r|s}$ of
$\mathscr{PD}_{r|s}^{k|l}$.

Let $p,t,p',t',r,s$ be nonnegative integers
and set
$
k=\max\{p,p'\}$, $l=\max\{t,t'\}.
$
For $0\leq a\leq k$ and $0\leq b\leq l$, we define an
index set
\begin{equation*}
    \begin{aligned}
    \tilde{\mathbf I}_{a|b,k|l}:=\{k+l-a-b+1,\ldots,k+l\},
    \end{aligned}
\end{equation*}
The parity of each index is inherited from
$\mathbf I_{k|l}$.
Let $\mathscr A_{p|t,p'|t',r|s}$ be the 
subsuperalgebra of $\mathscr{PD}_{r|s}^{k|l}$ generated by
\begin{equation*}
    \begin{aligned}
        t_{i,j},
\quad\partial_{i',j},
\quad
i\in\tilde{\mathbf I}_{p'|t',k|l},
i'\in\tilde{\mathbf I}_{p|t,k|l},
j\in\mathbf I_{r|s},
    \end{aligned}
\end{equation*}
This subalgebra is stable under the action of
$1\otimes\mathbf{U}_q(\mathfrak{gl}(r|s))$.
We denote its invariant subalgebra by
\[
\begin{aligned}
(\mathscr A_{p|t,p'|t',r|s})_{(\epsilon_R)}
:=
\{
f\in\mathscr A_{p|t,p'|t',r|s}\mid 
(1\otimes y)\cdot f=\epsilon_R(y)f\
\text{for all }y\in U_q(\mathfrak{gl}(r|s))
\}.
\end{aligned}
\]
For
$i\in\tilde{\mathbf I}_{p|t,k|l}$ and
$j\in\tilde{\mathbf I}_{p'|t',k|l}$,
the quantum polarization operator
$
L_{i,j}
$ belongs to $\mathscr A_{p|t,p'|t',r|s}$.
Our next result shows that these operators generate all
its invariants.

\noindent\textbf{Main Theorem 2.}
With the notation above, we have

\begin{equation*}
    \begin{aligned}
        (\mathscr A_{p|t,p'|t',r|s})_{(\epsilon_R)}
=
\left\langle
L_{i,j}\ \middle|\
i\in\tilde{\mathbf I}_{p|t,k|l},\
j\in\tilde{\mathbf I}_{p'|t',k|l}
\right\rangle.
    \end{aligned}
\end{equation*}

The main ideas of our proofs follow those of \cite{LSS24}.
However, the extension to the super setting involves more than
simply adding parity signs. To prove the FFT
(Theorem~\ref{theorem: AeR is generated by  Lij}),
we use a different argument based on specialization at $q=1$.
Using an integral form constructed from the PBW basis,
we combine the classical super FFT with dimension comparisons
under specialization to prove generation of the graded
invariant algebra. We then lift the generators to the
original invariant algebra through a filtration argument.

\subsection*{Acknowledgment}
We thank Li Luo for proposing the problem addressed in this paper and for offering valuable advice. Wek also thank Xirui Yu for many insightful discussions and helpful suggestions.

\section{Preliminaries}\label{Section:Preliminaries}

\subsection{The superalgebra $\mathbf{U}_q(\mathfrak{gl}(m|n))$}

Denote $\mathbf{I}_{m|n}=\{1,2,...,m+n\}$ and
 $\mathbf{I}_{m|n}'=\mathbf{I}_{m|n} \setminus \{m+n\}$ for $m,n\in \mathbb{Z}_{>0}$.
For any vector superspace $V=V_0\oplus V_1$, let $[~]:V_0\cup V_1\xrightarrow[]{}\mathbb{Z}_2$ be the parity map defined by $[x]=\bar{i}$ if $x\in V_i$. Let $\{ \epsilon_a \}$ be a basis of a vector space $A$ with a bilinear form $(\epsilon_a,\epsilon_b)=(-1)^{[a]}\delta_{ab}$  where $[~]:\mathbb{Z}\xrightarrow[]{}\mathbb{Z}_2={\{\bar{0},\bar{1}\}}$ is the parity map defined by
$[a]=
\begin{cases}
\bar{0} & \text{if }a\leq m,\\
\bar{1},  & \text{if } a>m.
\end{cases}$

The general linear quantum supergroup $\mathbf{U}_q(\mathfrak{gl}(m|n))$ is defined as the unital associative superalgebra over 
 $\mathbb{K}=\mathbb{C}(q)$, generated by 
 $$\{K_a,K_a^{-1},E_{b,b+1},E_{b+1,b} \mid a\in \mathbf{I}_{m|n},b\in \mathbf{I}'_{m|n}\},$$
 where  $E_{m,m+1},E_{m+1,m}$ are odd and the others are even,
 subject to the following relations:

$$
K_aK_a^{-1}=K_a^{-1}K_a=1,\quad K_aK_b=K_bK_a,
$$
$$
K_aE_{b,b\pm1}K_a^{-1}=q_a^{\delta_{ab}-\delta_{a,b\pm1}}E_{b,b\pm1},
$$
$$
[E_{a,a+1},E_{b+1,b}]=\delta_{ab}\frac{K_aK_{a+1}^{-1}-K_a^{-1}K_{a+1}}{q_a-q_a^{-1}},
$$
$$
E_{a,a+1}E_{b,b+1}=E_{b,b+1}E_{a,a+1},\quad
E_{a+1,a}E_{b+1,b}=E_{b+1,b}E_{a+1,a},\quad|a-b|\geq2,
$$
$$
(E_{a,a+1})^2 E_{b,b+1} - (q + q^{-1}) E_{a,a+1} E_{b,b+1} E_{a,a+1} + E_{b,b+1}(E_{a,a+1})^2 = 0, \quad |a-b| = 1, \ a \neq m,
$$
$$
(E_{a+1,a})^2 E_{b+1,b} - (q + q^{-1}) E_{a+1,a} E_{b+1,b} E_{a+1,a} + E_{b+1,b}(E_{a+1,a})^2 = 0, \quad |a-b| = 1, \ a \neq m,
$$
$$
(E_{m,m+1})^2 = (E_{m+1,m})^2 = [E_{m-1,m+2}, E_{m,m+1}] = [E_{m+2,m-1}, E_{m+1,m}] = 0,$$
where $q_a=q^{(\epsilon_a,\epsilon_a)}=q^{(-1)^{[a]}}$, $[x,y]=xy-(-1)^{[x][y]}yx$
and $E_{m-1,m+2} ,E_{m+2,m-1}$ 
are special cases of the elements:
$$E_{a,b}=E_{a,c}E_{c,b}-q_c^{-1}E_{c,b}E_{a,c}, \quad E_{b,a}=E_{b,c}E_{c,a}-q_cE_{c,a}E_{b,c},~a<c<b.$$

It is well known that the strcture of Hopf superalgebra of $\mathbf{U}_q(\mathfrak{gl}(m|n))$ is given by:

\noindent comultiplication
$$\Delta(E_{a,a+1})=E_{a,a+1}\otimes K_aK_{a+1}^{-1}+1\otimes E_{a,a+1},$$
$$\Delta(E_{a+1,a})=E_{a+1,a}\otimes 1+ K_a^{-1}K_{a+1}\otimes E_{a+1,a},$$
$$\Delta(K_a^{\pm1})=K_a^{\pm1}\otimes K_a^{\pm1},$$
counit
$$
\epsilon(E_{a,a+1}) = \epsilon(E_{a+1,a}) = 0,\quad \epsilon(K^{\pm 1}_{b}) = 1,
$$
and antipode
$$
S(E_{a,a+1}) = -E_{a,a+1}K^{-1}_{a}K_{a+1},\quad
S(E_{a+1,a}) = -K_{a}K^{-1}_{a+1}E_{a+1,a},
\quad
S(K^{\pm 1}_{a}) = K^{\mp 1}_{a}.
$$
We will use Sweedler's notation $\Delta(x)=\sum x_1\otimes x_2$ for all $x\in\mathbf{U}_q(\mathfrak{gl}(m|n))$.

Then, let $[\cdot,\cdot]_{q^{\pm}}$ be the deformed supercommutators defined in \cite[Formula 2.13]{KT91}, that is, 
\begin{equation*}
    \begin{aligned}
        [x,y]_{q^{\pm}}=xy-(-1)^{[x][y]}q^{\pm(\mathrm{wt}(x),\mathrm{wt}(y))}yx, \quad \forall x,y\in \mathbf{U}_q(\mathfrak{gl}(m|n)),
    \end{aligned}
\end{equation*}
where $\mathrm{wt}(x)$ is the weight of $x\in \mathbf{U}_q(\mathfrak{gl}(m|n))$ by the adjoint action of the Cartan elements. More precisely, $\mathrm{wt}(x)=\beta$ if $K_axK_a^{-1}=q^{(\epsilon_a,\beta)}x$ for $a\in \mathbf{I}_{m|n}$.
For any weight $\lambda=\sum_{i=1}^{m+n}a_i\epsilon_i$, $a_i\in \mathbb{Z}$, we denote $K_{\lambda}:=\prod_{i=1}^{m+n}K_{i}^{a_i}$. For brevity, we write $E_a:=E_{a,a+1},F_a:=E_{a+1,a}$. 
Let $\Delta$ be the set of roots of $\mathfrak{gl}(m|n)$ and $\Delta^+$ be the set of positive roots of $\mathfrak{gl}(m|n)$. Explicitly, $\Delta=\{\epsilon_i-\epsilon_j,1\leq i\neq j\leq m+n\} $ and $\Delta^+=\{\epsilon_i-\epsilon_j,1\leq i< j\leq m+n\}.$
The standard  root vectors of $\mathbf{U}_q(\mathfrak{gl}(m|n))$ are
\begin{equation}\label{equation: definition of root vectors}
    \begin{aligned}
        E_{\epsilon_i-\epsilon_j}:=[E_i,[\cdots[E_{j-2},E_{j-1}]_q]_q]_q,\quad F_{\epsilon_i-\epsilon_j}:=[F_{j-1},[\cdots[F_{i+1},F_i]_{q^{-1}}]_{q^{-1}}]_{q^{-1}},
    \end{aligned}
\end{equation}
where $1\leq i<j\leq m+n$. Then, for any $\gamma\in \Delta^+$, there exists a functor $a_\gamma=a_{\gamma}(q)$ satisfies the relation:
\begin{equation}\label{equation:relations for the deformed supercommutators}
    \begin{aligned}
        [E_\gamma,F_\gamma]=a_\gamma\frac{K_\gamma-K_{\gamma}^{-1}}{q-q^{-1}},
    \end{aligned}
\end{equation}
see \cite[Proposition 3.2]{KT91}.

\subsection{Universal R matrix}

Recall that the Drinfeld version of $\mathbf{U}_q(\mathfrak{gl}(m|n))$ admits a universal R-matrix $\mathfrak{R}$ in a completion of $\mathbf{U}_q(\mathfrak{gl}(m|n)) \otimes \mathbf{U}_q(\mathfrak{gl}(m|n)) $, but the Jimbo version of $\mathbf{U}_q(\mathfrak{gl}(m|n)) $ is not quasitriangular and
there exists no universal R-matrix belonging
 to $\mathbf{U}_q(\mathfrak{gl}(m|n)) \otimes \mathbf{U}_q(\mathfrak{gl}(m|n)) $ or any completion of it. Fortunately, the universal R-matrix $\mathfrak{R}$ of the Drinfeld version makes perfect sense within the Jimbo version of $\mathbf{U}_q(\mathfrak{gl}(m|n)) $ and gives rise to a functorial universal R-matrix in the category of locally finite $\mathbf{U}_q(\mathfrak{gl}(m|n))$-modules of type $1$, which gives a canonical braiding of this category, see \cite[Section 3.1]{Z98} and \cite[Section 3.5]{LZ17}.

 Let $\mathcal{C}^{(m|n)}$ be  the category of   locally finite $\mathbf{U}_q(\mathfrak{gl}(m|n))$-modules of type $1$. 
 Let $R^{(m|n)}$ be the functorial universal R-matrix in $\mathcal{C}^{(m|n)}$ .
 That means $R^{(m|n)}_{V,W}:V\otimes W\xrightarrow[]{}V\otimes W$ is an isomorphism of $\mathbf{U}_q(\mathfrak{gl}(m|n))$-modules  for any $V,W\in $ ob($\mathcal{C}^{(m|n)}$), and the family of isomorphisms
 \begin{equation}\label{definition of check R}
     \begin{aligned}
         \check{R}^{(m|n)}_{V,W}:=\tau_{V,W} \circ R^{(m|n)}_{V,W}:V\otimes W\xrightarrow[]{}W\otimes V
     \end{aligned}
 \end{equation}
 gives a braiding for $\mathcal{C}^{(m|n)}$, where $\tau_{V,W} (v\otimes w)=(-1)^{[v][w]}w\otimes v $ for  $v\in V,w\in W$.

Let $\underline{R}^{(m|n)}:=(R^{(m|n)}_{21})^{-1}=\tau\circ (R^{(m|n)})^{-1}\circ \tau$ and then $\check{\underline{R}}^{(m|n)}_{V,W}:=\tau_{V,W} \circ \underline{R}^{(m|n)}_{V,W}$ also gives a  braiding for $\mathcal{C}^{(m|n)}$.

$R^{(m|n)}$ satisfies the relations:
\begin{gather}
      R^{(m|n)}\Delta(x)=\tau\circ\Delta(x)R^{(m|n)}, \quad \forall x\in\mathbf{U}_q(\mathfrak{gl}(m|n)),\\
        (\Delta\otimes 1)R^{(m|n)}=R^{(m|n)}_{13}R^{(m|n)}_{23}, \quad (1\otimes\Delta)R^{(m|n)}=R^{(m|n)}_{13}R^{(m|n)}_{12},\\
       ( R^{(m|n)})^{-1}=(S\otimes1)R^{(m|n)},
\end{gather}
since $\mathfrak{R}$ satisfies the same relations by \cite[Section 2.3]{Zh20}.

Let the $q$-exponent be a formal series
$\mathrm{Exp}_{q}(x):=\sum\limits_{n\geq0} \frac{x^n}{(n)_{q}!}$,
where $(k)_q=\frac{1-q^k}{1-q}$ and $(n)_q!=(1)_q(2)_q\cdots(n)_q$.
Let $E_\gamma$ be the element defined in (\ref{equation: definition of root vectors}), $q_\gamma=(-1)^{[E_\gamma]}q^{-(\gamma,\gamma)}$. For any odd root $\gamma\in \Delta^+$, we have $q_{\gamma}=-1$ and $E_{\gamma}^2=F_{\gamma}^2=0$. To make the following formula well-defined, we adopt the convention
$
\exp_{-1}(x):=1+x
$
whenever $x^2=0$. By \cite[Theorem 9.1]{KT91}, for a given normal ordering in $\Delta^+$, we have

\begin{equation}
    \begin{aligned}
        R^{(m|n)}=K\cdot\prod_{\gamma\in\Delta^+}\mathrm{Exp}_{q_\gamma}((-1)^{[E_\gamma]}(a_\gamma)^{-1}(q-q^{-1})E_\gamma \otimes F_\gamma),
    \end{aligned}
\end{equation}
where $a_\gamma$ is the functor of $q$ satisfying (\ref{equation:relations for the deformed supercommutators}) and $K$ is the functor on the category $\mathcal{C}^{(m|n)}$ satisfying that for any $v\otimes w\in V\otimes W$, if $\mathrm{wt}(v)=\mu$ and $\mathrm{wt}(w)=\nu$, that means, $K_a v=q^{(\epsilon_a,\mu)} v$ and $K_a w=q^{(\epsilon_a,\nu)} w$ for any $a\in \mathbf{I}_{m|n}$, then
\begin{equation*}
    \begin{aligned}
        K\cdot(v\otimes w)=q^{(\mu,\nu)}v\otimes w.
    \end{aligned}
\end{equation*}
Note that for any $v\in V,w\in W,\gamma\in\Delta^+$, there exists $N\in\mathbb{N}$ such that $ (E_\gamma)^Nv=(F_\gamma)^Nv=(E_\gamma)^Nw=(F_\gamma)^Nw=0 $ for $k>N$, since $V,W$ are locally finite  $\mathbf{U}_q(\mathfrak{gl}(m|n))$-modules of type $1$.Therefore, the action of $R^{(m|n)}$ is well-defined. Using a construction completely analogous to that in the proof of \cite[Proposition 3.2.2]{LSS24}, we have the following result.

\begin{proposition} \label{proposition: the action of R is multipication by a element}
    If $V,W\in \mathrm{ob}(\mathcal{C}^{(m|n)})$ are finite dimensional, then there exists  elements $c_{V,W},c_{V,W}'$ of $ \mathbf{U}_q(\mathfrak{gl}(m|n))\otimes\mathbf{U}_q(\mathfrak{gl}(m|n)) $ such that the actions of $R^{(m|n)},\underline{R}^{(m|n)}$ on $V\otimes W$ is given by  multiplication with the elements $c_{V,W},c_{V,W}'$ respectively.
\end{proposition}

Recall that, for a Hopf superalgebra $H$, an associative superalgebra $M$ with unit is called a $H$-module superalgebra if  it is a $H$-module satisfying
\begin{equation}
    \begin{aligned}
        x.(ab)=\sum(-1)^{[x_2][a]}(x_1.a)(x_2.b) \quad \forall x\in H,\quad a,b\in M.
    \end{aligned}
\end{equation}

Let $A,B$ be $\mathbf{U}_q(\mathfrak{gl}(k|l))\otimes\mathbf{U}_q(\mathfrak{gl}(r|s))$-module superalgebras and they can be decomposed as
    \begin{equation}\label{equation: the dirct sum decomposition of A,B}
        \begin{aligned}
            A=\bigoplus_i(V_i\otimes V'_i) ,\quad B=\bigoplus_j(W_j\otimes W'_j),
        \end{aligned}
    \end{equation}
    where $V_i,W_j\in \mathrm{ob}(\mathcal{C}^{(k|l)}),V'_i,W'_j\in \mathrm{ob}(\mathcal{C}^{(r|s)})$  are finite dimensional. We define a family of $\mathbf{U}_q(\mathfrak{gl}(k|l))\otimes\mathbf{U}_q(\mathfrak{gl}(r|s))\otimes\mathbf{U}_q(\mathfrak{gl}(k|l))\otimes\mathbf{U}_q(\mathfrak{gl}(r|s)) $-modules isomorphisms, which mapping $W_j\otimes W'_j\otimes V_i\otimes V_i'$ onto $V_i\otimes V_i' \otimes W_j\otimes W'_j$:
    \begin{equation*}
        \begin{aligned}
            \check{R}'_{W_j\otimes W_j',V_i\otimes V_i'}:=(\check{\underline{R}}^{(k|l)}_{W_j,V_i})_{13}\circ (\check{\underline{R}}^{(r|s)}_{W'_j,V'_i})_{24}
            =
            \tau_{23}\circ
            (\check{\underline{R}}^{(k|l)}_{W_j,V_i})_{12}\circ (\check{\underline{R}}^{(r|s)}_{W'_j,V'_i})_{34}\circ
            \tau_{23}
        \end{aligned}
    \end{equation*}
where $(\check{\underline{R}}^{(k|l)}_{W_j,V_i})_{12}$ means that $\check{\underline{R}}^{(k|l)}_{W_j,V_i}$ acts on the first and the second components of $W_j\otimes W'_j\otimes V_i\otimes V_i'$ and the definition of  $(\check{\underline{R}}^{(r|s)}_{W'_j,V'_i})_{34}$ is similar.
Using these isomorphisms, we can construct a even linear isomorphism:
\begin{equation}\label{equation: the definition of R_{B,A}}
    \begin{aligned}
        \check{R}'_{B,A}:=\bigoplus_{i,j} \check{R}'_{W_j\otimes W_j',V_i\otimes V_i'}:
        B\otimes A\xrightarrow[]{}A\otimes B.
    \end{aligned}
\end{equation}

We now introduce the $\check{R}'$-twisted tensor product  $A\otimes_{\check{R}'}B$ which is  $A\otimes B$  as a vector space and is endowed with the multiplication
\begin{equation}\label{equation: the product of twisted product}
    \begin{aligned}
        m_{A,B}=(m_A\otimes m_B)\circ(id_A \otimes  \check{R}'_{B,A}\otimes id_B),
    \end{aligned}
\end{equation}
where $m_A,m_B$ are the multiplications of $A,B$ respectively.

\begin{remark}\label{remark: the module structure of twisted tensor product}
    We know that $A\otimes_{\check{R}'}B$ have the natural $\mathbf{U}\otimes \mathbf{U}$-module structure since $A,B$ are $\mathbf{U}$-module. Note that $H$ is a Hopf superalgebras since the tensor product of two Hopf superalgebras is also a Hopf superalgebras. Therefore, $A\otimes_{\check{R}'}B$ has a $\mathbf{U}$-module structure via the coproduct $\mathbf{U}\xrightarrow[]{}\mathbf{U}\otimes \mathbf{U}$.
\end{remark}
Next, we will use the following lemma to show that $A\otimes_{\check{R}'}B$ is also a $\mathbf{U}_q(\mathfrak{gl}(k|l))\otimes\mathbf{U}_q(\mathfrak{gl}(r|s))$-module
superalgebras.
\begin{lemma}\label{the property of R'}
     The isomorphism $\check{R}'_{B,A}$ satisfies the following relations:
    \begin{gather}
    \label{equation:first for R}
        \check{R}'_{B,A}(1\otimes a)=a\otimes1, \quad \check{R}'_{B,A}(b\otimes 1)=1\otimes b ,\quad\forall a\in A,b\in B,\\
    \label{equation:second for R}
        \check{R}'_{B,A}\circ(\mathrm{id}_B\otimes m_A)=
        (m_A\otimes \mathrm{id}_B)(\mathrm{id}_A\otimes\check{R}'_{B,A})(\check{R}'_{B,A}\otimes \mathrm{id}_A),\\
    \label{equation:third for R}
        \check{R}'_{B,A}\circ(m_B\otimes \mathrm{id}_A)=
        (\mathrm{id}_A\otimes m_B)(\check{R}'_{B,A}\otimes \mathrm{id}_B)(\mathrm{id}_B\otimes\check{R}'_{B,A}),
    \end{gather}           
\end{lemma}
\begin{proof}
   We first show that
$\check {R}'_{B,A}
=
\check{\underline R}^{(k|l)}_{B,A}
\circ\tau_{A,B}
\circ\check{\underline R}^{(r|s)}_{B,A}$. It suffices to verify this identity on each direct summand
$(W_j\otimes W_j')\otimes(V_i\otimes V_i')$ of $B\otimes A$. Since $\check{\underline R}^{(k|l)}$ is a braiding for $\mathcal{C}^{(k|l)}$,
its naturality implies that
\begin{equation*}
    \begin{aligned}
    \check{\underline R}^{(k|l)}_{B,A}|_{(W_j\otimes W_j')\otimes(V_i\otimes V_i')}
    =
    \check{\underline R}^{(k|l)}_{W_j\otimes W_j',V_i\otimes V_i'}=\tau_{23}\circ( \check{\underline R}^{(k|l)}_{W_j,V_i}\otimes \tau_{W_j',V_i'})\circ \tau_{23}.
    \end{aligned}
\end{equation*}
  Similarly, we have   $\check{\underline R}^{(r|s)}_{B,A}|_{(W_j\otimes W_j')\otimes(V_i\otimes V_i')}
    =
    \check{\underline R}^{(r|s)}_{W_j\otimes W_j',V_i\otimes V_i'}=\tau_{23}\circ( \tau_{W_j,V_i}\otimes\check{\underline R}^{(r|s)}_{W_j',V_i'} )\circ \tau_{23}.$ We note that $\tau_{23}\circ \tau_{A,B}\circ \tau_{23}=\tau_{V_i,W_j}\otimes\tau_{V_i',W_j'}$, so that
    \begin{equation*}
        \begin{aligned}
            &\check{\underline R}^{(k|l)}_{B,A}\circ\tau_{A,B}\circ\check{\underline R}^{(r|s)}_{B,A}|_{W_j\otimes W_j',V_i\otimes V_i'}
            \\
            =&\tau_{23}\circ( \check{\underline R}^{(k|l)}_{W_j,V_i}\otimes \tau_{W_j',V_i'})\circ(\tau_{V_i,W_j}\otimes\tau_{V_i',W_j'})\circ( \tau_{W_j,V_i}\otimes\check{\underline R}^{(r|s)}_{W_j',V_i'} )\circ \tau_{23}
            \\
            =&\tau_{23}\circ( \check{\underline R}^{(k|l)}_{W_j,V_i}\otimes\check{\underline R}^{(r|s)}_{W_j',V_i'} )\circ \tau_{23}=\check {R}'_{B,A}|_{W_j\otimes W_j',V_i\otimes V_i'}.
        \end{aligned}
    \end{equation*}
    By the compatibility of the braidings with the unit object
and their naturality with respect to the unit map
$\mathbb K\xrightarrow[]{} B$, we have
\begin{equation*}
    \begin{aligned}
        \check{\underline R}^{(k|l)}_{B,A}(1\otimes a)
        =a\otimes1,
        \qquad
        \check{\underline R}^{(r|s)}_{B,A}(1\otimes a)
        =a\otimes 1, \quad \forall a\in A.
    \end{aligned}
\end{equation*}
Hence, we have $\check {R}'_{B,A}(1\otimes a)
=
\check{\underline R}^{(k|l)}_{B,A}
\circ\tau_{A,B}
\circ\check{\underline R}^{(r|s)}_{B,A}(1\otimes a)=a \otimes 1$ for any  $a\in A.$. Similarly, we can obtain the other equality in (\ref{equation:first for R}). 
The proof of (\ref{equation:second for R}) and (\ref{equation:third for R}) are similar, so we only give the proof for  (\ref{equation:second for R}).
Since $A$ is a module superalgebra, its multiplication
$m_A:A\otimes A\to A$ is an even module homomorphism.
By the naturality of the two braidings, we have
\begin{equation*}
    \begin{aligned}
        \check{\underline R}^{(r|s)}_{B,A}
        \circ(\mathrm{id}_B\otimes m_A)
        &=(m_A\otimes\mathrm{id}_B)
        \circ\check{\underline R}^{(r|s)}_{B,A\otimes A},\\
        \check{\underline R}^{(k|l)}_{B,A}
        \circ(\mathrm{id}_B\otimes m_A)
        &=(m_A\otimes\mathrm{id}_B)
        \circ\check{\underline R}^{(k|l)}_{B,A\otimes A}.
    \end{aligned}
\end{equation*}
Moreover, the naturality of the super flip gives
\begin{equation*}
    \begin{aligned}
        \tau_{A,B}\circ(m_A\otimes\mathrm{id}_B)
        =
        (\mathrm{id}_B\otimes m_A)\circ\tau_{A\otimes A,B}.
    \end{aligned}
\end{equation*}

Thus, we have
\begin{equation}\label{equation: proof in relations for R}
    \begin{aligned}
        \check R'_{B,A}\circ(\mathrm{id}_B\otimes m_A)
        &=
        \check{\underline R}^{(k|l)}_{B,A}
        \circ\tau_{A,B}
        \circ\check{\underline R}^{(r|s)}_{B,A}
        \circ(\mathrm{id}_B\otimes m_A)
        \\
        &=
        \check{\underline R}^{(k|l)}_{B,A}
        \circ\tau_{A,B}
        \circ(m_A\otimes\mathrm{id}_B)
        \circ\check{\underline R}^{(r|s)}_{B,A\otimes A}
        \\
        &=
        \check{\underline R}^{(k|l)}_{B,A}
        \circ(\mathrm{id}_B\otimes m_A)
        \circ\tau_{A\otimes A,B}
        \circ\check{\underline R}^{(r|s)}_{B,A\otimes A}
        \\
        &=
        (m_A\otimes\mathrm{id}_B)
        \circ\check{\underline R}^{(k|l)}_{B,A\otimes A}
        \circ\tau_{A\otimes A,B}
        \circ\check{\underline R}^{(r|s)}_{B,A\otimes A}
        \\
        &=
        (m_A\otimes\mathrm{id}_B)\circ\check R'_{B,A\otimes A}.
    \end{aligned}
\end{equation}
By the hexagon identities for the two braidings, we have
\begin{equation*}
    \begin{aligned}
        \check{\underline R}^{(k|l)}_{B,A\otimes A}
        &=
        (\mathrm{id}_A\otimes\check{\underline R}^{(k|l)}_{B,A})
        \circ
        (\check{\underline R}^{(k|l)}_{B,A}\otimes\mathrm{id}_A),\\
        \check{\underline R}^{(r|s)}_{B,A\otimes A}
        &=
        (\mathrm{id}_A\otimes\check{\underline R}^{(r|s)}_{B,A})
        \circ
        (\check{\underline R}^{(r|s)}_{B,A}\otimes\mathrm{id}_A).
    \end{aligned}
\end{equation*}
Note that $\tau_{A\otimes A,B}
=
(\tau_{A,B}\otimes\mathrm{id}_A)
\circ(\mathrm{id}_A\otimes\tau_{A,B})$, we have
\begin{equation*}
    \begin{aligned}
        \check R'_{B,A\otimes A}
        =&
        \check{\underline R}^{(k|l)}_{B,A\otimes A}
        \circ\tau_{A\otimes A,B}
        \circ\check{\underline R}^{(r|s)}_{B,A\otimes A}
        \\
        =&
        (\mathrm{id}_A\otimes\check{\underline R}^{(k|l)}_{B,A})
        \circ(\check{\underline R}^{(k|l)}_{B,A}\otimes\mathrm{id}_A)
        \circ(\tau_{A,B}\otimes\mathrm{id}_A)
        \\
        &\quad\circ
        (\mathrm{id}_A\otimes\tau_{A,B})
        \circ(\mathrm{id}_A\otimes\check{\underline R}^{(r|s)}_{B,A})
        \circ(\check{\underline R}^{(r|s)}_{B,A}\otimes\mathrm{id}_A)
        \\=&
        \Bigl(
        \mathrm{id}_A\otimes
        \bigl(
        \check{\underline R}^{(k|l)}_{B,A}
        \circ\tau_{A,B}
        \circ\check{\underline R}^{(r|s)}_{B,A}
        \bigr)
        \Bigr)
        \circ
        \Bigl(
        \bigl(
        \check{\underline R}^{(k|l)}_{B,A}
        \circ\tau_{A,B}
        \circ\check{\underline R}^{(r|s)}_{B,A}
        \bigr)
        \otimes\mathrm{id}_A
        \Bigr)
        \\
        =&
        (\mathrm{id}_A\otimes\check R'_{B,A})
        \circ(\check R'_{B,A}\otimes\mathrm{id}_A).
    \end{aligned}
\end{equation*}
Combining the identity above with (\ref{equation: proof in relations for R}), we obtain (\ref{equation:second for R}).
\end{proof}

\begin{proposition}\label{proposition:twisted tensor is module superalgebra}
     $A\otimes_{\check{R}'}B$ is  a $\mathbf{U}_q(\mathfrak{gl}(k|l))\otimes\mathbf{U}_q(\mathfrak{gl}(r|s))$-module superalgebra.
\end{proposition}

\begin{proof}
    It follows from Lemma \ref{the property of R'} and \cite[Proposition 2.2]{VV94} that $A\otimes_{\check{R}'}B$ is an associative superalgebra. Given any $a_1,a_2\in A$ and $b_1,b_2\in B$, by the definition of $m_{AB}$ ,we have
    \begin{equation*}
        \begin{aligned}
            (a_1\otimes b_1)(a_2\otimes b_2)=(m_A\otimes m_B)(a_1\otimes \check{R}'_{B,A}(b_1\otimes a_2)\otimes b_2).
        \end{aligned}
    \end{equation*}
   We suppose that $\check{R}'_{B,A}(b_1\otimes a_2)=\sum a'\otimes b'$ where $a'\in A,b'\in B$, then we have
    \begin{equation*}
        \begin{aligned}
            (a_1\otimes b_1)(a_2\otimes b_2)=\sum a_1a'\otimes b' b_2.
        \end{aligned}
    \end{equation*}
It follows that for any $x\in\mathbf{U}_q(\mathfrak{gl}(k|l))\otimes\mathbf{U}_q(\mathfrak{gl}(r|s))$
    \begin{equation*}
        \begin{aligned}
           &x\cdot (a_1\otimes b_1)(a_2\otimes b_2)\\
           =&\sum(-1)^{[x_2]([a_1]+[a'])} x_1\cdot(a_1a')\otimes x_2\cdot(b' b_2)\\
           =&\sum(-1)^{[x_2]([a_1]+[a'])+[x_{12}][a_1]+[x_{22}][b']}((x_{11}\cdot a_1)(x_{12} \cdot a'))\otimes((x_{21}\cdot b')(x_{22}\cdot b_2)),
        \end{aligned}
    \end{equation*}
    where the second equality follows from the $\mathbf{U}_q(\mathfrak{gl}(k|l))\otimes\mathbf{U}_q(\mathfrak{gl}(r|s))$-module superalgebra structure of $A$ and $B$.
    By the (\ref{equation:first for R}), we have
    \begin{equation*}
       \begin{aligned}
           ((x_{11}\cdot a_1)(x_{12} \cdot a'))\otimes((x_{21}\cdot b')(x_{22}\cdot b_2))
           =
           ((x_{11}\cdot a_1)\otimes1)((x_{12}\cdot a')\otimes (x_{21}\cdot b'))(1\otimes (x_{22}\cdot b_2)).
       \end{aligned}
    \end{equation*}
    Since $\check{R}'_{B,A}$ is a $\mathbf{U}_q(\mathfrak{gl}(k|l))\otimes\mathbf{U}_q(\mathfrak{gl}(r|s))$-module isomorphism , we have
   \begin{equation*}
       \begin{aligned}
           \sum(x_{12}\cdot a')\otimes(x_{21}\cdot b')
           &=
           \sum(-1)^{[x_{21}][a']}(x_{12}\otimes x_{21})\cdot\check{R}'_{B,A}(b_1\otimes a_2)\\
           &=\sum(-1)^{[x_{21}][b_1]+[x_{21}][a']}\check{R}'_{B,A}( (x_{12}\cdot b_1)\otimes (x_{21}\cdot a_2)).
       \end{aligned}
   \end{equation*} 
We shall note that $[a']+[b']=[b_1]+[a_2]$, since  $\check{R}'_{B,A}$ is an isomorphism.
Then,
\begin{equation*}
    \begin{aligned}
         x\cdot (a_1\otimes b_1)(a_2\otimes b_2)
           &=\sum(-1)^r((x_{11}\cdot a_1)\otimes1)\check{R}'_{B,A}( (x_{12}\cdot b_1)\otimes (x_{21}\cdot a_2))(1\otimes (x_{22}\cdot b_2))\\
           &=\sum(-1)^r((x_{11}\cdot a_1)\otimes (x_{12}\cdot b_1))((x_{21}\cdot a_2)\otimes(x_{22}\cdot b_2))\\
           &=\sum(-1)^{r+[x_{12}][a_1]+[x_{22}][a_2]}(x_1\cdot(a_1\otimes b_1)(x_2\cdot(a_2\otimes b_2))),
    \end{aligned}
\end{equation*}
where $r={[x_2]([a_1]+[a'])+[x_{12}][a_1]+[x_{22}][b']+[x_{21}][b_1]+[x_{21}][a']} $. The second equality follows from the definition of $m_{AB}$ and the third equality follows from the module structure of $A\otimes_{\check{R}'}B$. By using the relation $[x_2]=[x_{21}]+[x_{22}]$, $[a']+[b']=[b_1]+[a_2]$, we have $r=[x_2]([a_1]+[b_1])$. Therefore,  $A\otimes_{\check{R}'}B$ is  a $\mathbf{U}_q(\mathfrak{gl}(k|l))\otimes\mathbf{U}_q(\mathfrak{gl}(r|s))$-module superalgebra.
\end{proof}

\section{The $q$-Weyl superalgebra $\mathscr{PD}_{r|s}^{k|l}$}\label{Section:PD}

\subsection{The superalgebra $\mathscr{P}_{m|n} $ and $\mathscr{D}_{m|n}$}

Let $V^{m|n}$ be the irreducible $\mathbf{U}_q(\mathfrak{gl}(m|n))$-module of highest weight $\epsilon_1$ with the standard  basis $\{ e_a \mid a\in \mathbf{I}_{m|n} \}$. The corresponding representation $\pi:\mathbf{U}_q(\mathfrak{gl}(m|n))\xrightarrow[]{} \mathrm{End}(V^{m|n})$ is determined  by:
\begin{equation}\label{equation: the module structure of V}
    \begin{aligned}      
        E_a\mapsto e_{a,a+1}, \quad
        F_a\mapsto e_{a+1,a},\quad
        K_a\mapsto I+(q_a-1)e_{a,a},
    \end{aligned}
\end{equation}
where $I$ is the identity matrix and  $e_{i,j}$ are the elementary matrix units associated to the standard  basis $\{ e_a \mid a\in \mathbf{I}_{m|n} \}$.

Let $\partial_{i,j}\in \mathbf{U}_q(\mathfrak{gl}(m|n))^{\circ}$ denote  the matrix coefficients of $V^{m|n}$, defined by 
\begin{align*}
    \langle \partial_{a,b},x\rangle = \langle e_a^*,x\cdot e_b\rangle = \pi(x)_{ab}, \quad\forall x\in \mathbf{U}_q(\mathfrak{gl}(m|n)),
\end{align*}
where $\{ e_a^* \mid a\in \mathbf{I}_{m|n} \}$ is the dual basis of $\{ e_a \mid a\in \mathbf{I}_{m|n} \}$. We set the $\mathbb{Z}_2$-grading of $\partial_{i,j}$ by $[\partial_{i,j}]=[i]+[j]$.

Denote $\mathscr{D}_{m|n}$ be the subsuperalgebra of $ \mathbf{U}_q(\mathfrak{gl}(m|n))^{\circ}$ generated by $\partial_{i,j}$, $i,j\in \mathbf{I}_{m|n}$. The multiplication of $\mathscr{D}_{m|n}$ is defined by
\begin{equation*}
    \begin{aligned}
        \langle\partial \partial',x\rangle=
        \sum_{(x)}(-1)^{[\partial'][x_1]} \langle\partial,x_1\rangle\langle\partial',x_2\rangle,
    \end{aligned}
\end{equation*}
where $\partial,\partial'\in D_{m|n}$, $x\in \mathbf{U}_q(\mathfrak{gl}(m|n))$.

\begin{lemma}[{\cite[Lemma 2.6]{Zh20}}]\label{lemma: generate relations of D}
The superalgebra $\mathscr{D}_{m|n}$ is generated by $\partial_{i,j}$, $i,j\in \mathbf{I}_{m|n}$, subject to the following relations:
\begin{equation}\label{equation:relations for D}
    \begin{aligned}
        (\partial_{ab})^{2} &= 0, &\qquad [a]+[b]=\bar{1},\\
        \partial_{ac}\partial_{bc} 
        &= (-1)^{([a]+[c])([b]+[c])}\, q_{c}\,\partial_{bc}\partial_{ac},
        &a>b,\\
        \partial_{ab}\partial_{ac} 
        &= (-1)^{([a]+[c])([a]+[b])}\, q_{a}\,\partial_{ac}\partial_{ab},
        & b>c,\\
        \partial_{ac}\partial_{bd} 
        &= (-1)^{([a]+[c])([b]+[d])}\,\partial_{bd}\partial_{ac},
        & a>b,\ c<d,\\
        \partial_{ac}\partial_{bd} 
        &= (-1)^{([a]+[c])([b]+[d])}\,\partial_{bd}\partial_{ac}\\
          &\quad + (-1)^{[a]([b]+[d])+[b][d]}\,(q-q^{-1})\,\partial_{bc}\partial_{ad},
        &a>b,\ c>d.
      \end{aligned}
  \end{equation}
\end{lemma}
Note that $\mathscr{D}_{m|n}$ has a bi-superalgebra structure with comultiplication and counit given by
\begin{align*}
    \Delta(\partial_{ab})=\sum_{c\in\mathbf{I}_{m|n} }(-1)^{([a]+[c])([c]+[b])}\partial_{a,c}\otimes\partial_{c,b},\ \ \epsilon(\partial_{a,b})=\delta_{a,b}.
\end{align*}

Similarly, $(V^{m|n})^*$ is the irreducible $\mathbf{U}_q(\mathfrak{gl}(m|n))$-module of highest weight $-\epsilon_{m+n}$ with the standard basis $\{ e_a' \mid a\in \mathbf{I}_{m|n} \}$. The corresponding representation $\tilde{\pi}:\mathbf{U}_q(\mathfrak{gl}(m|n))\xrightarrow[]{} \mathrm{End}((V^{m|n})^*)$ is determined  by:
\begin{equation}\label{equation: the module structure of V*}
    \begin{aligned}  
        E_a\mapsto -q_a^{-1}e_{a+1,a},\quad
        F_a\mapsto (-1)^{\delta_{a,m}+1}q_ae_{a,a+1},\quad
        K_a\mapsto I+(q_a^{-1}-1)e_{a,a}.
    \end{aligned}
\end{equation}
where  $I$ is the identity matrix and $e_{i,j}$ are the elementary matrix units associated to the standard  basis $\{ e'_a \mid a\in \mathbf{I}_{m|n} \}$.

Let $t_{i,j}\in \mathbf{U}_q(\mathfrak{gl}(m|n))^{\circ}$ be the matrix coefficients of $(V^{m|n})^*$, defined by 
\begin{align*}
    \langle t_{a,b},x\rangle = \langle(e_a')^*,x\cdot e_b'\rangle= \tilde{\pi}(x)_{ab},\quad\forall x\in \mathbf{U}_q(\mathfrak{gl}(m|n)).
\end{align*}

Denote $\mathscr{P}_{m|n}$ be the subsuperalgebra of $ \mathbf{U}_q(\mathfrak{gl}(m|n))^{\circ}$ generated by $t_{i,j}$, $i,j\in \mathbf{I}_{m|n}$. The multiplication of $\mathscr{P}_{m|n}$ is defined by
\begin{equation*}
    \begin{aligned}
        \langle tt',x\rangle=
        \sum_{(x)}(-1)^{[t'][x_1]} \langle t,x_1\rangle\langle t',x_2\rangle,
    \end{aligned}
\end{equation*}
where $t,t'\in P_{m|n}$, $x\in \mathbf{U}_q(\mathfrak{gl}(m|n))$.

\begin{lemma}[{\cite[Remark 2.12]{Zh20}}]\label{lemma: generate relation of P}
The superalgebra $\mathscr{P}_{m|n}$ is generated by $t_{i,j}$, $i,j\in \mathbf{I}_{m|n}$, subject to the following relations:
\begin{equation}\label{equation:relations for P}
    \begin{aligned}
        (t_{ab})^{2} &= 0,
        &\qquad [a]+[b]=\bar{1},\\[4pt]
        t_{ac}t_{bc}
        &= (-1)^{([a]+[c])([b]+[c])}\, q_c\, t_{bc}t_{ac},
        &\qquad a<b,\\[4pt]
        t_{ab}t_{ac}
        &= (-1)^{([a]+[c])([a]+[b])}\, q_a\, t_{ac}t_{ab},
        &\qquad b<c,\\[4pt]
        t_{ac}t_{bd}
        &= (-1)^{([a]+[c])([b]+[d])}\, t_{bd}t_{ac},
        &\qquad a<b,\ c>d,\\
        t_{ac}t_{bd}&= (-1)^{([a]+[c])([b]+[d])}\, t_{bd}t_{ac}  \\
          &\quad +(-1)^{[a]([b]+[d])+[b][d]}\,(q-q^{-1})\, t_{bc}t_{ad},
        &\qquad a<b,\ c<d.
    \end{aligned}
\end{equation}
    
\end{lemma}
Note that $\mathscr{P}_{m|n}$ has a bi-superalgebra structure with comultiplication and counit given by
\begin{align*}
    \Delta(t_{ab})=\sum_{c\in\mathbf{I}_{m|n} }(-1)^{([a]+[c])([c]+[b])}t_{a,c}\otimes t_{c,b},\ \ \epsilon(t_{a,b})=\delta_{a,b}.
\end{align*}


\begin{proposition}\label{proposition: the basis of D}(\cite[Proposition 3.8]{Zh20})
   The monomials
    \begin{equation}\label{equation: the basis of D}
        \begin{aligned}
            \partial_{1,1}^{a_{1,1}}\partial_{1,2}^{a_{1,2}}\cdots \partial_{1,m+n}^{a_{1,m+n}}\cdots \partial_{m+n,1}^{a_{m+n,1}} \cdots \partial_{m+n,m+n}^{a_{m+n,m+n}}
        \end{aligned}
    \end{equation}
    form a $\mathbb{K}$-basis of $\mathscr{D}_{m|n}$, where $a_{i,j}\in \mathbb{Z}_{\geq 0} $ whenever $ [i]+[j]=\bar{0}$  and $a_{i,j}\in \{0,1\} $ whenever $ [i]+[j]=\bar{1}$.
\end{proposition}

\begin{proposition}\label{proposition: the basis of P}
   The monomials
    \begin{equation}
        \begin{aligned}
             t_{m+n,m+n}^{a_{m+n,m+n}}\cdots t_{m+n,1}^{a_{m+n,1}}\cdots t_{1,m+n}^{a_{1,m+n}}\cdots t_{1,1}^{a_{1,1}}
        \end{aligned}
    \end{equation}
    form a $\mathbb{K}$-basis of $\mathscr{P}_{m|n}$, where $a_{i,j}\in \mathbb{Z}_{\geq 0} $ whenever $ [i]+[j]=\bar{0}$  and $a_{i,j}\in \{0,1\} $ whenever $ [i]+[j]=\bar{1}$.
\end{proposition}
\begin{proof}
    By Lemma \ref{lemma: generate relations of D} and \ref{lemma: generate relation of P}, there exists a superalgebra anti-isomorphism $\Upsilon:\mathscr{D}_{m|n}\xrightarrow[]{}\mathscr{P}_{m|n}$ such that
    $\Upsilon(\partial_{i,j})=(-1)^{[j]([i]+[j])}t_{j,i}$.
    Note that 
    \begin{equation*}
        \begin{aligned}
            \Upsilon( \partial_{1,1}^{a_{1,1}}\partial_{1,2}^{a_{2,1}}\cdots \partial_{1,m+n}^{a_{m+n,1}}\cdots \partial_{m+n,1}^{a_{1,m+n}} \cdots \partial_{m+n,m+n}^{a_{m+n,m+n}})=(-1)^{k}   t_{m+n,m+n}^{a_{m+n,m+n}}\cdots t_{m+n,1}^{a_{m+n,1}}\cdots t_{1,m+n}^{a_{1,m+n}}\cdots t_{1,1}^{a_{1,1}}
        \end{aligned}
    \end{equation*}
    for some $k\in \mathbb{Z}$, the assertion follows from Proposition \ref{proposition: the basis of D}
\end{proof}
Note that $\mathscr{P}_{m|n}$ (resp. $\mathscr{D}_{m|n}$) are $\mathbb{Z}$-graded with $\mathrm{deg}~t_{a,b}=1$ (resp. $\mathrm{deg}~\partial_{a,b}=1$). For any non-negative integer $d$, let $\mathscr{P}_{m|n}^{(d)}$ (resp.$\mathscr{D}_{m|n}^{(d)}$ ) be the homogeneous subspaces of degree $d$.


\subsection{The $\mathbf{U}_q(\mathfrak{gl}(m|n))\otimes \mathbf{U}_q(\mathfrak{gl}(m|n))$-module structure of $\mathscr{P}_{m|n}$ and $\mathscr{D}_{m|n}$}\label{subsection: module structure of P,D}

Firstly, we define two left actions $\mathcal{L}_{\mathscr{D}}$ and $\mathcal{R}_\mathscr{D}$ of $\mathbf{U}_q(\mathfrak{gl}(m|n))$ on $\mathscr{D}_{m|n}$. We define $\mathcal{R}_{\mathscr{D}}$ by the left translation:
\begin{align*}
    \langle \mathcal{R}_{\mathscr{D}}(x)f,y\rangle=\langle f,y\cdot x\rangle=\sum(-1)^{[y][f_2]}\langle f_1,y\rangle \langle f_2,x\rangle,
\end{align*}
where $\Delta(f)=\sum f_1\otimes f_2$. Then we have
\begin{equation*}
    \begin{aligned}
        \mathcal{R}_{\mathscr{D}}(x)f=\sum (-1)^{[x]([f]+[x])}\langle f_2,x\rangle f_1,
    \end{aligned}
\end{equation*}
 since $\langle f_1,y\rangle=0 $ unless $[f_1]=[y]$ and $[f]=[f_1]+[f_2]$. In addition, $\mathscr{D}_{m|n}$ is a $\mathbf{U}_q(\mathfrak{gl}(m|n))$-module superalgebra under $R_{\mathscr{D}}$.

There is an anti-automorphism $\natural: \mathbf{U}_q(\mathfrak{gl}(m|n))\xrightarrow[]{}\mathbf{U}_q(\mathfrak{gl}(m|n))^{\mathrm{op}}$ given by
\begin{equation*}
    \begin{aligned}
        E_a \mapsto q_{a+1}K_aK_{a+1}^{-1}F_a, \quad
        \quad F_a\mapsto(-1)^{\delta_{a,m}}q_{a+1}^{-1}E_aK_a^{-1}K_{a+1},\quad
        K_a\mapsto K_a.
    \end{aligned}
\end{equation*}

We define $\mathcal{L}_{\mathscr{D}}$ by the right translation:
\begin{align*}
   \langle \mathcal{L}_{\mathscr{D}}(x)f,y\rangle=(-1)^{[x][y]}\langle f, x^\natural y\rangle=
   \sum(-1)^{[x][y]+[f_2][x]}\langle f_1,x^{\natural}\rangle \langle f_2,y\rangle,
\end{align*}
where $\Delta(f)=\sum f_1\otimes f_2$.
Then we have 
\begin{equation*}
    \begin{aligned}
        \mathcal{L}_{\mathscr{D}}(x)f=\sum\langle f_1,x^{\natural}\rangle f_2,
    \end{aligned}
\end{equation*}
 since $\langle f_2,y\rangle=0 $ unless $[f_2]=[y]$.
 Also, $\mathscr{D}_{m|n}$ is a $\mathbf{U}_q(\mathfrak{gl}(m|n))$-module superalgebra under $\mathcal{L}_{\mathscr{D}}$.
Note that, for any $x,y\in \mathbf{U}_q(\mathfrak{gl}(m|n))$, $f\in \mathscr{D}
 $, we have 
 \begin{equation*}
     \begin{aligned}
         \mathcal{L}_{\mathscr{D}}(x)\mathcal{R}_{\mathscr{D}}(y)f=(-1)^{[x][y]}\mathcal{R}_{\mathscr{D}}(y)\mathcal{L}_{\mathscr{D}}(x)f=\sum(-1)^{[y]([f]+[y])}\langle f_1,x^{\natural}\rangle \langle f_3,y \rangle f_2,
     \end{aligned}
 \end{equation*}
 where $(\Delta\otimes1)\circ\Delta(f)=\sum f_1\otimes f_2\otimes f_3$. Therefore, we can make $\mathscr{D}_{m|n} $ a $\mathbf{U}_q(\mathfrak{gl}(m|n))\otimes \mathbf{U}_q(\mathfrak{gl}(m|n))$ module superalgebra by defining $(x\otimes y).f =\mathcal{L}_{\mathscr{D}}(x)\mathcal{R}_{\mathscr{D}}(y)f$.

 Explicitly, we have
 \begin{equation}\label{equation: formulas of the action on D}
     \begin{aligned}
        (E_a\otimes1)\partial_{i,j}&=(-1)^{\delta_{a,m}[j]}\delta_{i,a+1}\partial_{a,j},\quad (1\otimes E_a)\partial_{i,j}=\delta_{j,a+1}\partial_{i,a},\\
        (F_a\otimes1)\partial_{i,j}&=(-1)^{\delta_{a,m}[j]}\delta_{i,a}\partial_{a+1,j}, \quad(1\otimes F_a)\partial_{i,j}=\delta_{j,a}\partial_{i,a+1},\\
        (K_b\otimes1)\partial_{i,j}&=q_b^{\delta_{i,b}}\partial_{i,j},\quad\quad\quad\quad\quad\quad(1\otimes K_b)\partial_{i,j}=q_b^{\delta_{j,b}}\partial_{i,j},\\
     \end{aligned}
 \end{equation}
where $a\in \mathbf{I}_{m|n}'$, $b\in \mathbf{I}_{m|n}$.

Next, we define two left actions $\mathcal{L}_{\mathscr{P}}$ and $\mathcal{R}_\mathscr{P}$ of $\mathbf{U}_q(\mathfrak{gl}(m|n))$ on $\mathscr{P}_{m|n}$.

 We define $ \mathcal{R}_{\mathscr{P}} $ and $ \mathcal{L}_{\mathscr{P}} $  as
\begin{equation}\label{equation: formulas of Rp,Lp}
    \begin{aligned}
        \mathcal{R}_{\mathscr{P}}(x)f&=\sum(-1)^{[x]([f]+[x])}\langle  f_2,x\rangle f_1,\\
        \mathcal{L}_{\mathscr{P}}(x)f&=\sum \langle  f_1\circ S,(S^{-1}(x))^{\natural}\rangle f_2.
    \end{aligned}
\end{equation}
By directly calculating, we have
\begin{equation}\label{lem:S}
    \natural\circ S^{-1}=S\circ\natural\ \ \mathrm{i.e.}\ \ (S^{-1}(x))^\natural=S(x^\natural),~\forall x\in \mathbf{U}_q(\mathfrak{gl}(m|n)).
\end{equation}

The antipode $S$ of $\mathbf{U}_q(\mathfrak{gl}(m|n))$ induces an isomorphism of Hopf superalgebras  
\begin{equation}\label{equation: definition  of Phi}
    \begin{aligned}
      \Phi:\mathbf{U}_q(\mathfrak{gl}(m|n))^{\circ}& \xrightarrow{} (\mathbf{U}_q(\mathfrak{gl}(m|n))^{\circ})^{\mathrm{op,cop}},\\
        u&\mapsto u \circ S.
    \end{aligned}
\end{equation}
For any sequence $\{c_a \mid a\in \mathbf{I}_{m|n}\}$ of elements of $\mathbb{K}$, let $\eta_{c_a}$ be the automorphism of $\mathbf{U}_q(\mathfrak{gl}(m|n))$ defined by
\begin{equation*}
    \begin{aligned}
        E_a\mapsto c_aE_a, \quad F_a\mapsto c_a^{-1}F_a ,\quad K_a\mapsto K_a.
    \end{aligned}
\end{equation*}
By (\ref{lem:S}) and $S^2=\eta_{q_aq_{a+1}}$, we obtain more concise formulas:
\begin{equation*}
    \begin{aligned}
 \mathcal{L}_{\mathscr{P}}(x)f&=\sum \langle f_1, S^2(x^{\natural})\rangle f_2=\sum\langle f_1\circ\eta_{q_a q_{a+1}}, x^\natural\rangle f_2.
    \end{aligned}
\end{equation*}
Similar to $\mathcal{L}_{\mathscr{D}}$ and $\mathcal{R}_\mathscr{D}$, we have
\begin{equation*}
    \begin{aligned}
        \mathcal{L}_{\mathscr{P}}(x)\mathcal{R}_{\mathscr{P}}(y)f=(-1)^{[x][y]}\mathcal{R}_{\mathscr{P}}(y)\mathcal{L}_{\mathscr{P}}(x)f=\sum(-1)^{[y]([f]+[y])}\langle f_1\circ
 \eta_{q_aq_{a+1}},x^{\natural}\rangle \langle f_3,y \rangle f_2,
    \end{aligned}
\end{equation*}
 where $(\Delta\otimes1)\circ\Delta(f)=\sum f_1\otimes f_2\otimes f_3$. Therefore, we can also make $\mathscr{P}_{m|n} $ a $\mathbf{U}_q(\mathfrak{gl}(m|n))\otimes \mathbf{U}_q(\mathfrak{gl}(m|n))$ module superalgebra by defining $(x\otimes y).f =\mathcal{L}_{\mathscr{P}}(x)\mathcal{R}_{\mathscr{P}}(y)f$. 
Explicitly, we have
 \begin{equation}\label{equation: formulas of the action on P}
     \begin{aligned}
        (E_a\otimes1)t_{i,j}&=(-1)^{\delta_{a,m}[j]+1}q_a^{-1}\delta_{i,a}t_{a+1,j},
        \quad \quad(1\otimes E_a)t_{i,j}=-q_a^{-1}\delta_{j,a}t_{i,a+1},\\
        (F_a\otimes1)t_{i,j}&=(-1)^{\delta_{a,m}(1+[j])+1}q_a\delta_{i,a+1}t_{a,j},
        \quad
        (1\otimes F_a)t_{i,j}=(-1)^{\delta_{a,m}+1}q_a\delta_{j,a+1}t_{i,a}, \\
        (K_b\otimes1)t_{i,j}&=q_b^{-\delta_{i,b}}t_{i,j},
        \quad\quad\quad\quad\quad\quad\quad\ \ \quad(1\otimes K_b)t_{i,j}=q_b^{-\delta_{j,b}}t_{i,j},
     \end{aligned}
 \end{equation}
where $a\in \mathbf{I}_{m|n}'$, $b\in \mathbf{I}_{m|n}$.

\subsection{The superalgebra $\mathscr{P}^{k|l}_{r|s}$ and $\mathscr{D}^{k|l}_{r|s}$ }\label{subsection: P and D}
In this part, let  $k,l,r,s$ be non-negative integers satisfying $k,r\leq m$ and $l,s\leq n$, we define the superalgebra $\mathscr{P}^{k|l}_{r|s}$   and $\mathscr{D}^{k|l}_{r|s}$ as subsuperalgebras of $\mathscr{P}_{m|n}$  and $\mathscr{D}_{m|n}$ respectively, and give the $ \mathbf{U}_q(\mathfrak{gl}(k|l))\otimes \mathbf{U}_q(\mathfrak{gl}(r|s))$ module structure of them. Let $\hat{\mathbf{I}}_{k|l}=\{a\in \mathbb{Z}_{\geq 0}\mid m-k+1\leq a \leq m+l\}$ and $\hat{\mathbf{I}}_{r|s}=\{a\in \mathbb{Z}_{\geq 0}\mid m-r+1\leq a \leq m+s\}$.

Let $\mathscr{P}^{k|l}_{r|s}$ be the subsuperalgebra of $\mathscr{P}_{m|n}$  generated by $t_{i,j}$, $i \in  \hat{\mathbf{I}}_{k|l} $, $j\in  \hat{\mathbf{I}}_{r|s}$ and $\mathscr{D}^{k|l}_{r|s}$ be the subsuperalgebra of $\mathscr{D}_{m|n}$  generated by $\partial_{i,j}$, $i \in  \hat{\mathbf{I}}_{k|l} $, $j\in  \hat{\mathbf{I}}_{r|s}$.

Note that we can view $ \mathbf{U}_q(\mathfrak{gl}(k|l))$ as a Hopf subsuperalgebra of $\mathbf{U}_q(\mathfrak{gl}(m|n))$ generated by $\{E_a,F_a,K_b\mid a\in \hat{\mathbf{I}}_{k|l},a\neq m+l,b\in\hat{\mathbf{I}}_{k|l} \}$ via the monomorphism of superalgebras
\begin{equation}\label{equation: embedding of quantum groups}
    \begin{aligned}
        \kappa^{k|l}_{m|n}:\mathbf{U}_q(\mathfrak{gl}(k|l))& \xrightarrow{} \mathbf{U}_q(\mathfrak{gl}(m|n)),\\
        E_i\mapsto E_{m-k+i}, \quad F_i&\mapsto F_{m-k+i} ,\quad K_i\mapsto K_{m-k+i}.
    \end{aligned}
\end{equation}
Then the action of $\mathbf{U}_q(\mathfrak{gl}(m|n))\otimes \mathbf{U}_q(\mathfrak{gl}(m|n))$ on $\mathscr{P}_{m|n}$(respectively, $\mathscr{D}_{m|n}$) can induce the action of $ \mathbf{U}_q(\mathfrak{gl}(k|l))\otimes \mathbf{U}_q(\mathfrak{gl}(r|s))$ on  $\mathscr{P}^{k|l}_{r|s}$ (respectively, $\mathscr{D}^{k|l}_{r|s}$ ) via $\kappa^{k|l}_{m|n}\otimes\kappa^{r|s}_{m|n}$.
\begin{proposition}
     The superalgebras $\mathscr{P}^{k|l}_{r|s}$ and $\mathscr{D}^{k|l}_{r|s}$ are $ \mathbf{U}_q(\mathfrak{gl}(k|l))\otimes \mathbf{U}_q(\mathfrak{gl}(r|s))$-module superalgebras.
\end{proposition}
\begin{proof}
   We only prove this for  $\mathscr{P}^{k|l}_{r|s}$ .
   It suffices to show that  $\mathscr{P}^{k|l}_{r|s}$ is stable under the action of $ \mathbf{U}_q(\mathfrak{gl}(k|l))\otimes \mathbf{U}_q(\mathfrak{gl}(r|s))$, since  $\mathscr{P}_{m|n}$ is a $\mathbf{U}_q(\mathfrak{gl}(m|n)\otimes \mathbf{U}_q(\mathfrak{gl}(m|n))$-module superalgebra.
   Let $x$ be a generator of $ \mathbf{U}_q(\mathfrak{gl}(k|l))$, then $\kappa^{k|l}_{m|n}(x)\in \{E_a,F_a,K_b\mid a\in \hat{\mathbf{I}}_{k|l},a\neq k+l,b\in\hat{\mathbf{I}}_{k|l} \} $ and  for $i\in \hat{\mathbf{I}}_{k|l}$, $j\in \hat{\mathbf{I}}_{r|s}$, we have $(x\otimes1)\cdot t_{ij}=(\kappa^{k|l}_{m|n}(x)\otimes 1)\cdot t_{i,j}\in \mathscr{P}^{k|l}_{r|s}$ by (\ref{equation: formulas of the action on P}). 
   Since  $\mathscr{P}_{m|n}$ is a $\mathbf{U}_q(\mathfrak{gl}(m|n)\otimes \mathbf{U}_q(\mathfrak{gl}(m|n))$-module superalgebra,  $\mathscr{P}^{k|l}_{r|s}$ is stable under the action of $ \mathbf{U}_q(\mathfrak{gl}(k|l))\otimes 1$. 
   Similarly, we can prove that $\mathscr{P}^{k|l}_{r|s}$ is stable under the action of $ 1\otimes \mathbf{U}_q(\mathfrak{gl}(r|s))$ by (\ref{equation: formulas of the action on D}).
\end{proof}

\begin{proposition}\label{proposition: equivalent definition of P,D}
   The superalgebras $\mathscr{P}^{k|l}_{r|s}$ and  $\mathscr{D}^{k|l}_{r|s}$ can be viewed as quotient superalgebras as follows.
    
   (1) $\mathscr{P}^{k|l}_{r|s}$ is isomorphic to the quotient of the  free superalgebra generated by $t_{i,j}$ , $i\in \mathbf{I}_{k|l},j\in \mathbf{I}_{r|s}$ with $[t_{i,j}]=[i]+[j]$ modulo (\ref{equation:relations for P}).
    
    (2) $\mathscr{D}^{k|l}_{r|s}$ is isomorphic to the quotient of the  free superalgebra generated by $\partial_{i,j}$ , $i\in \mathbf{I}_{k|l},j\in \mathbf{I}_{r|s} $ with $[\partial_{i,j}]=[i]+[j]$  modulo (\ref{equation:relations for D}).
\end{proposition}
\begin{proof}
    We just give the proof of (1), since the proof of (2) is similar. Let $\mathscr{A}$ be the quotient of the  free superalgebra generated by $t_{i,j}$ with $i\in \mathbf{I}_{k|l},j\in \mathbf{I}_{r|s}$ modulo (\ref{equation:relations for P}). It follows from (\ref{equation:relations for P}) that there exists a superalgebra homomorphism 
    \begin{equation*}
        \begin{aligned}
            \Theta:\mathscr{A}\xrightarrow[]{} \mathscr{P}^{k|l}_{r|s},\quad
            t_{ij}\mapsto t_{m-k+i,m-r+j}.
        \end{aligned}
    \end{equation*}
 Moreover, the homomorphism $\Theta$ is surjective since $\mathscr{P}^{k|l}_{r|s}$ is generated by  $t_{i,j}$, $i \in  \hat{\mathbf{I}}_{k|l} $, $j\in  \hat{\mathbf{I}}_{r|s}$. By using the relations (\ref{equation:relations for P}),  we can express every element of $\mathscr{A}$ as a linear combination of 
 the monomials
    \begin{equation}
        \begin{aligned}
            t_{k+l,r+s}^{a_{k+l,r+s}}\cdots t_{k+l,1}^{a_{k+l,1}}\cdots t_{1,r+s}^{a_{1,r+s}}\cdots t_{1,1}^{a_{1,1}},
        \end{aligned}
    \end{equation}
     where $a_{i,j}\in \mathbb{Z}_{\geq 0} $ whenever $ [i]+[j]=\bar{0}$  and $a_{i,j}\in \{0,1\} $ whenever $ [i]+[j]=\bar{1}$. It follows from  Proposition \ref{proposition: the basis of P} that $\Theta$ maps the monomials to  a linearly independent set of elements of $\mathscr{P}^{k|l}_{r|s}$, so the monomials are linearly independent and form a $\mathbb{K}$-basis of  $\mathscr{A}$. Hence, $\Theta$ is injective. Therefore, $\Theta:\mathscr{A}\xrightarrow[]{} \mathscr{P}^{k|l}_{r|s}$ is an superalgebra isomorphism.

\end{proof}

In what follows, we always view  $\mathscr{P}^{k|l}_{r|s}$ and  $\mathscr{D}^{k|l}_{r|s}$ as the quotient superalgebras in proposition \ref{proposition: equivalent definition of P,D} and the $ \mathbf{U}_q(\mathfrak{gl}(k|l))\otimes \mathbf{U}_q(\mathfrak{gl}(r|s))$-module superalgebras structures of  the quotient superalgebras are induced by those of  $\mathscr{P}^{k|l}_{r|s}$ and  $\mathscr{D}^{k|l}_{r|s}$. Then, the explicit formulas of the actions on the generators $t_{i,j},\partial_{i,j}$, $i\in \mathbf{I}_{k|l},j\in \mathbf{I}_{r|s}$ are given by (\ref{equation: formulas of the action on P}) and (\ref{equation: formulas of the action on D}). 

By the proof of Proposition \ref{proposition: equivalent definition of P,D}, we obtain the bases of $\mathscr{P}^{k|l}_{r|s}$ and $\mathscr{D}^{k|l}_{r|s}$.

\begin{proposition}\label{proposition: the basis of P and D}
     The monomials
    \begin{equation}
        \begin{aligned}
            \partial_{1,1}^{a_{1,1}}\partial_{1,2}^{a_{1,2}}\cdots \partial_{1,r+s}^{a_{1,r+s}}\cdots \partial_{k+l,1}^{a_{k+l,1}} \cdots \partial_{k+l,r+s}^{a_{k+l,r+s}}
        \end{aligned}
    \end{equation}
    form a $\mathbb{K}$-basis of $\mathscr{D}^{k|l}_{r|s}$ and
    the monomials
    \begin{equation}
        \begin{aligned}
            t_{k+l,r+s}^{a_{k+l,r+s}}\cdots t_{k+l,1}^{a_{k+l,1}}\cdots t_{1,r+s}^{a_{1,r+s}}\cdots t_{11}^{a_{11}}
        \end{aligned}
    \end{equation}
    form a $\mathbb{K}$-basis of $\mathscr{P}^{k|l}_{r|s}$ , where $a_{i,j}\in \mathbb{Z}_{\geq 0} $ whenever $ [i]+[j]=\bar{0}$  and $a_{i,j}\in \{0,1\} $ whenever $ [i]+[j]=\bar{1}$. 
\end{proposition}

\subsection{multiplicity-free decomposition of $\mathscr{P}^{k|l}_{r|s}$ and $\mathscr{D}^{k|l}_{r|s}$}

Let $\lambda=(\lambda_1,\lambda_2,...,\lambda_{m+n})$ be a $(m,n)$-hook partition, which means $\lambda_i\geq \lambda_{i+1}$ for $1\leq\ i\leq m+n-1$ and $\lambda_{m+1}\leq n$.
We denote by $\Lambda_{m|n}$  the set of $(m,n)$-hook partitions and denote by
 $L_{\lambda}^{m|n}$  the irreducible $\mathbf{U}_q(\mathfrak{gl}(m|n))$-module with the highest weight $\lambda^\natural=\sum\limits_{a=1}^{m}\lambda_a\epsilon_a+\sum\limits_{b=1}^{n}\langle\lambda_{b}'-m\rangle\epsilon_{m+b}$, where $\lambda'=(\lambda'_1,\lambda_2',...,\lambda_{m+n}')$ is the transpose partition of $\lambda$ and $\langle\lambda_{b}'-m\rangle=\max\{ \lambda_{b}'-m,0\}$.

As $\mathbf{U}_q(\mathfrak{gl}(k|l))\otimes\mathbf{U}_q(\mathfrak{gl}(r|s))$-modules, $\mathscr{P}^{k|l}_{r|s}$ is isomorphic to the superalgebra $\overline{M}^{k|l}_{r|s}$ under the action $\mathcal{L}\otimes\mathcal{R}$ and  $\mathscr{D}^{k|l}_{r|s}$ is isomorphic to the superalgebra ${M}^{k|l}_{r|s}$ under the action $\mathcal{L}\otimes\mathcal{R}$, where  $\overline{M}^{k|l}_{r|s}$, ${M}^{k|l}_{r|s}$, $\mathcal{L}$ and $\mathcal{R}$ are   defined in \cite[\S 2.4 \& \S 3.1]{Zh20}.
This can be checked by  calculating the formulas of the action of $\mathbf{U}_q(\mathfrak{gl}(k|l))\otimes\mathbf{U}_q(\mathfrak{gl}(r|s))$ on the generators of $\overline{M}^{k|l}_{r|s}$ and ${M}^{k|l}_{r|s}$ like (\ref{equation: formulas of the action on D}) and (\ref{equation: formulas of the action on P}). Therefore, by \cite[Theorem 2.8 \& Theorem 2.10]{Zh20}, we have the following multiplicity-free decompositions.

\begin{proposition}\label{proposition: multiplicity-free decompositions of P and D}
    As  $\mathbf{U}_q(\mathfrak{gl}(k|l))\otimes\mathbf{U}_q(\mathfrak{gl}(r|s))$-modules,
    \begin{equation*}
        \begin{aligned}
            \mathscr{P}^{k|l}_{r|s}\cong \bigoplus_{\lambda \in \Lambda_{k|l}\cap \Lambda_{r|s}}(L_{\lambda}^{k|l})^*\otimes (L_{\lambda}^{r|s})^*,\quad
              \mathscr{D}^{k|l}_{r|s}\cong \bigoplus_{\lambda \in \Lambda_{k|l}\cap \Lambda_{r|s}}L_{\lambda}^{k|l}\otimes L_{\lambda}^{r|s}.
        \end{aligned}
    \end{equation*}
\end{proposition}
Note that $\mathscr{P}^{k|l}_{r|s}$ (respectively, $\mathscr{D}^{k|l}_{r|s}$ ) is $\mathbb{Z}$-graded as the subspace of  $\mathscr{P}_{m|n}$ (respectively, $\mathscr{D}_{m|n}$). Let $\mathscr{P}^{k|l,(d)}_{r|s}$ (respectively, $\mathscr{D}^{k|l,(d)}_{r|s}$ ) be the homogeneous subspaces of degree $d$. Then, we have
\begin{equation*}
    \begin{aligned}
         \mathscr{P}^{k|l,(d)}_{r|s}\cong \bigoplus_{\lambda \in\Lambda_{k|l}\cap \Lambda_{r|s},\lvert\lambda\rvert=d}(L_{\lambda}^{k|l})^*\otimes (L_{\lambda}^{r|s})^*,\quad
              \mathscr{D}^{k|l,(d)}_{r|s}\cong \bigoplus_{\lambda \in\Lambda_{k|l}\cap \Lambda_{r|s},\lvert\lambda\rvert=d}L_{\lambda}^{k|l}\otimes L_{\lambda}^{r|s}.
    \end{aligned}
\end{equation*}

\subsection{The superalgebra $\mathscr{PD}^{k|l,gr}_{r|s}$ and $\mathscr{PD}^{k|l}_{r|s}$}\label{subsection: PD}

It follows from Proposition \ref{proposition:twisted tensor is module superalgebra} and Proposition \ref{proposition: multiplicity-free decompositions of P and D} that we can obtain a $\mathbf{U}_q(\mathfrak{gl}(k|l))\otimes\mathbf{U}_q(\mathfrak{gl}(r|s))$-module superalgebra 
$\mathscr{PD}^{gr}=\mathscr{PD}^{k|l,gr}_{r|s}:=\mathscr{P}^{k|l}_{r|s}\otimes_{\check{R}'}\mathscr{D}^{k|l}_{r|s}$ with $\check{R}'_{\mathscr{D}^{k|l}_{r|s},\mathscr{P}^{k|l}_{r|s}}$.
We consider the linear map
\begin{equation*}
    \begin{aligned}
        m:\mathscr{P}^{k|l,(1)}_{r|s}\otimes \mathscr{D}^{k|l,(1)}_{r|s}&\xrightarrow[]{} T(\mathscr{P}^{k|l,(1)}_{r|s}\oplus \mathscr{D}^{k|l,(1)}_{r|s}),\\
        f\otimes g&\mapsto fg,
    \end{aligned}
\end{equation*}
where $T(X)$ denotes the tensor algebra on $X$. Then, it follows from \cite[Proposition 2.4.8]{LSS24} that the homomorphism of superalgebras 
\begin{equation}\label{equation: an homomorphism from T(P+D) to PD}
    \begin{aligned}
        \zeta : T(\mathscr{P}^{k|l,(1)}_{r|s}\oplus \mathscr{D}^{k|l,(1)}_{r|s})&\xrightarrow[]{}\mathscr{PD}^{k|l,gr}_{r|s},\\
        f\mapsto f\otimes1,\quad& g\mapsto1\otimes g \quad \mathrm{for}~f\in\mathscr{P}^{k|l,(1)} ,~g\in\mathscr{D}^{k|l,(1)}
    \end{aligned}
\end{equation}
induces an isomorphism of superalgebras
$\mathscr{PD}^{k|l,gr}_{r|s}\cong T(\mathscr{P}^{k|l,(1)}_{r|s}\oplus \mathscr{D}^{k|l,(1)}_{r|s})/I$, where
$I$ is the two-sided ideal of $T(\mathscr{P}^{k|l,(1)}_{r|s}\oplus \mathscr{D}^{k|l,(1)}_{r|s})$ generated by (\ref{equation:relations for D}), (\ref{equation:relations for P}) and 
\begin{equation}\label{equation: the definition of PD^gr}
    \begin{aligned}
        gf-m\circ\check{R}'_{ \mathscr{D}^{k|l}_{r|s},\mathscr{P}^{k|l}_{r|s}}(g\otimes f), \quad \forall f\in \mathscr{P}^{k|l,(1)}_{r|s},g\in \mathscr{D}^{k|l,(1)}_{r|s}.
    \end{aligned}
\end{equation}

We consider a bilinear map:
\begin{equation*}
    \begin{aligned}
        \psi :\mathscr{D}^{k|l,(1)}_{r|s}\times \mathscr{P}^{k|l,(1)}_{r|s}&\xrightarrow[]{}\mathbb{K},\\
       (\partial_{a,b},t_{c,d})&\mapsto (-1)^{[a]+[b]+[a][b]} \delta_{a,c} \cdot\delta_{b,d}
    \end{aligned}
\end{equation*}
and define $\mathscr{PD}=\mathscr{PD}^{k|l}_{r|s}:=T(\mathscr{P}^{k|l,(1)}_{r|s}\oplus \mathscr{D}^{k|l,(1)}_{r|s})/I_{\psi}$, where $I_\psi$ is the two-sided ideal of $T(\mathscr{P}^{k|l,(1)}_{r|s}\oplus \mathscr{D}^{k|l,(1)}_{r|s})$ generated by (\ref{equation:relations for D}), (\ref{equation:relations for P}) and 
\begin{equation}\label{equation: the definition of PD}
    \begin{aligned}
        gf-m\circ\check{R}'_{ \mathscr{D}^{k|l}_{r|s},\mathscr{P}^{k|l}_{r|s}}(g\otimes f)-\psi(g,f), \quad \forall f\in \mathscr{P}^{k|l,(1)}_{r|s}, g\in \mathscr{D}^{k|l,(1)}_{r|s}.
    \end{aligned}
\end{equation}

\begin{proposition}\label{proposition: the relations for PD}

    (1) $\mathscr{PD}^{k|l}_{r|s}$ is generated by $t_{a,b}$, $\partial_{a,b}$, $a\in \mathbf{I}_{k|l},b\in \mathbf{I}_{r|s}$ modulo (\ref{equation:relations for D}), (\ref{equation:relations for P}) and the following relations:
    \begin{align}
    \label{equation: the first relation for PD}
        \partial_{a,c}t_{b,d}=&t_{b,d}\partial_{a,c}\cdot(-1)^{([a]+[c])([b]+[d])}, (a\neq b,c\neq d),\\
    \label{equation: the second relation for PD}
        \partial_{a,b}t_{a,c}=&q_at_{a,c}\partial_{a,b}\cdot(-1)^{([a]+[c])([a]+[b])}
        +
        (q-q^{-1})\sum_{a'>a}t_{a',c}\partial_{a',b}\cdot(-1)^{[a']+([a]+[b])([c]+[a'])}, (b\neq c)\\
    \label{equation: the third relation for PD}
        \partial_{a,c}t_{b,c}=&q_c t_{b,c}\partial_{a,c}\cdot(-1)^{([b]+[c])([a]+[c])}
        +
        (q-q^{-1})\sum_{c'>c}t_{b,c'}\partial_{a,c'}\cdot(-1)^{[c']+([b]+[c])([a]+[c'])},(a\neq b) \\
    \label{equation: the fourth relation for PD}
        \partial_{a,b}t_{a,b}=&(-1)^{[a]+[b]+[a][b]}+q_aq_b t_{a,b}\partial_{a,b}\cdot(-1)^{[a]+[b]}\\
        &+
        (q-q^{-1})\sum_{b'>b}q_at_{a,b'}\partial_{a,b'}\cdot(-1)^{[b']+([a]+[b])([a]+[b'])}\notag\\
        &+(q-q^{-1})\sum_{a'>a}q_bt_{a',b}\partial_{a',b}\cdot(-1)^{[a']+([a]+[b])([b]+[a'])} \notag\\
        &+(q-q^{-1})^2\sum_{a'>a}\sum_{b'>b}t_{a',b'}\partial_{a',b'}\cdot(-1)^{[a']+[b']+([a]+[b'])([b]+[a'])}.\notag
    \end{align}

    (2) $\mathscr{PD}^{k|l,gr}_{r|s}$ is generated by $t_{i,j}$, $\partial_{i,j}$, $i\in \mathbf{I}_{k|l},j\in \mathbf{I}_{r|s}$ modulo the same relations  of  $\mathscr{PD}^{k|l}_{r|s}$ except  that the constant term $(-1)^{[a]+[b]+[a][b]}$ on the right side of (\ref{equation: the fourth relation for PD}) is removed.
\end{proposition}

\begin{proof}
    By the definition of $\mathscr{PD}$, it suffices to compute (\ref{equation: the definition of PD}) 
    for  any $g=\partial_{a,c}$ and $f=t_{b,d}$.    
    Note that $V^{k|l}\otimes V^{r|s}$ and $ (V^{k|l})^*\otimes(V^{r|s})^*$ have natural $\mathbf{U}_q(\mathfrak{gl}(k|l))\otimes\mathbf{U}_q(\mathfrak{gl}(r|s))$-module structures via (\ref{equation: the module structure of V}) and (\ref{equation: the module structure of V*}). It follows from (\ref{equation: formulas of the action on P}) and (\ref{equation: formulas of the action on D}) that
    $\mathscr{P}^{k|l,(1)}_{r|s}\cong (V^{k|l})^*\otimes(V^{r|s})^*$ and $ \mathscr{D}^{k|l,(1)}_{r|s}\cong V^{k|l}\otimes V^{r|s}$ as $\mathbf{U}_q(\mathfrak{gl}(k|l))\otimes\mathbf{U}_q(\mathfrak{gl}(r|s))$-modules  via the isomorphisms $\Psi_1: t_{i,j}\mapsto (-1)^{[i][j]}e_i'\otimes e_j'$, $\Psi_2:\partial_{i,j}\mapsto (-1)^{[i][j]}e_i\otimes e_j$. By using the isomorphisms, we can obtain $\check{R}'_{\mathscr{D}^{k|l}_{r|s},\mathscr{P}^{k|l}_{r|s}}(\partial_{a,c}\otimes t_{b,d})$ by computing $\check{R}'_{V^{k|l}\otimes V^{r|s},(V^{k|l})^*\otimes(V^{r|s})^*}(e_{a}\otimes e_c\otimes e_b'\otimes e_d') $ for $a,b\in\mathbf{I}_{k|l}$ ,$c,d\in \mathbf{I}_{r|s}$. By the (\ref{equation: the definition of R_{B,A}}), we have
    \begin{equation*}
    \begin{aligned}
        \check{R}'_{V^{k|l}\otimes V^{r|s},(V^{k|l})^*\otimes(V^{r|s})^*}=(\check{\underline{R}}^{(k|l)}_{V^{k|l},(V^{k|l})^*})_{13}\circ (\check{\underline{R}}^{(r|s)}_{V^{r|s},(V^{r|s})^*})_{24}.
    \end{aligned}
    \end{equation*}
    
    For simplicity, we will abuse notations in the following computations. We write $\check{R}'$ for any $\check{R}'_{B,A}$  and write $\check{\underline{R}}$ for both $\check{\underline{R}}^{(k|l)}_{V^{k|l},(V^{k|l})^*}$ and $\check{\underline{R}}^{(r|s)}_{V^{r|s},(V^{r|s})^*}$ whenever this causes no ambiguity. By using \cite[Theorem 1]{Z92}, we have
    \begin{equation*}
    \begin{aligned}
        \underline{R}^{(k|l)}_{V^{k|l},(V^{k|l})^*}=\sum_{j\in\mathbf{I}_{k|l}}q_j e_{j,j}\otimes e_{j,j}+\sum_{\substack{i\neq j\\i,j\in\mathbf{I}_{k|l}}}e_{i,i}\otimes e_{j,j}+(q-q^{-1})\sum_{\substack{i>j\\ i,j\in\mathbf{I}_{k|l}}}e_{i,j}\otimes e_{i,j}(-1)^{[j]},\\
         \underline{R}^{(r|s)}_{V^{r|s},(V^{r|s})^*}=\sum_{j\in\mathbf{I}_{r|s}}q_j e_{j,j}\otimes e_{j,j}+\sum_{\substack{i\neq j\\i,j\in\mathbf{I}_{r|s}}}e_{i,i}\otimes e_{j,j}+(q-q^{-1})\sum_{\substack{i>j\\ i,j\in\mathbf{I}_{r|s}}}e_{i,j}\otimes e_{i,j}(-1)^{[j]}.
    \end{aligned}
    \end{equation*}
    Then for any $e_a\in V^{k|l},e_b'\in (V^{k|l})^*,e_c\in V^{r|s}, e_d\in  (V^{r|s})^* $, we have
    \[
    \check{\underline{R}}(e_a\otimes e_b') =
    \left\{
    \begin{aligned}
        &(-1)^{[a][b]}e_b'\otimes e_a\quad &a\neq b,\\
        &(-1)^{[a]}q_ae_a'\otimes e_a+(q-q^{-1})\sum_{i>a}e_i'\otimes e_i\cdot(-1)^{[i]([a]+1)}\quad&a=b.
    \end{aligned}
    \right.
    \]
      \[
    \check{\underline{R}}(e_c\otimes e_d') =
    \left\{
    \begin{aligned}
        &(-1)^{[c][d]}e_d'\otimes e_c\quad &c\neq d,\\
        &(-1)^{[c]}q_ce_c'\otimes e_c+(q-q^{-1})\sum_{j>c}e_j'\otimes e_j\cdot(-1)^{[j]([c]+1)}\quad&c=d.
        \end{aligned}
    \right.
    \]
    Using the two formulas above, we can write 
    $\check{R}'(e_a\otimes e_c\otimes e_b'\otimes e_d')$ as a sum of the elements of the form  $(e_i'\otimes e_j'\otimes e_p\otimes e_q)$ by considering four separate cases. We only provide the detailed calculation for one case. Assume that $b\neq c$, we have
    \begin{equation*}
        \begin{aligned}
            \check{R'}(e_a\otimes e_b \otimes e_a' \otimes e_c')
            =&\tau_{23}\circ  \check{\underline{R}}_{12}\check{\underline{R}}_{34}(e_a\otimes e_a'\otimes e_b \otimes e_c')\cdot(-1)^{[a][b]}\\
            =&\tau_{23}((-1)^{[a]}q_a e_a'\otimes e_a
            \\&+(q-q^{-1})\sum_{i>a}e_i'\otimes e_i \cdot(-1)^{[i]([a]+1)})\otimes(e_c'\otimes e_b)\cdot(-1)^{([a]+[c])[b]})\\
            =&q_a e_a'\otimes e_c'\otimes e_a\otimes e_b \cdot (-1)^{([a]+[c])([a]+[b])}\\
            &+(q-q^{-1})\sum_{i>a} e_i'\otimes e_c'\otimes e_i \otimes e_b \cdot (-1)^{[i]+([a]+[c])([b]+[i])}.
        \end{aligned}
    \end{equation*}
    By using $\Psi_1,\Psi_2$, we have 
    \begin{equation*}
        \begin{aligned}
            \check{R'}(\partial_{a,b}\otimes t_{a,c})
            &=q_at_{a,c}\otimes\partial_{a,b}\cdot(-1)^{([a]+[c])([a]+[b])}
            \\
            &+
            (q-q^{-1})\sum_{i>a}t_{i,c}\otimes\partial_{i,b}\cdot(-1)^{[i]+([a]+[b])([c]+[i])}.
        \end{aligned}
    \end{equation*}
    Then we obtain the (\ref{equation: the second relation for PD}) by substituting the above formula into (\ref{equation: the definition of PD}).

    Comparing (\ref{equation: the definition of PD^gr}) to (\ref{equation: the definition of PD}), it follows that $\mathscr{PD}^{k|l,gr}_{r|s}$ is generated by $t_{i,j}$, $\partial_{i,j}$, $i\in \mathbf{I}_{k|l},j\in \mathbf{I}_{r|s}$ modulo the same relations  except  that the constant term $1$ on the right side of (\ref{equation: the fourth relation for PD}) is removed.
\end{proof}
In what follows, we always view $\mathscr{PD}^{k|l,gr}_{r|s}$ and $\mathscr{PD}^{k|l}_{r|s}$ as the quotient superalgebra of $T(\mathscr{P}^{k|l,(1)}_{r|s}\oplus \mathscr{D}^{k|l,(1)}_{r|s})$. In particular, the monomials in the following two propositions  are viewed as elements  of the quotient superalgebra of  $T(\mathscr{P}^{k|l,(1)}_{r|s}\oplus \mathscr{D}^{k|l,(1)}_{r|s})$.

We set
$
    \Xi_{k|l,r|s}:=\{A=(a_{i,j})\in \mathrm{Mat}_{(k+l)\times(r+s)}(\mathbb{N})\mid  a_{i,j}=0 \ \mathrm{or} \ 1 
    \ \mathrm{for}\ [i]+[j]=1\}
$. For $A,B\in \Xi_{k|l,r|s} $, we define
 \begin{equation*}
     \begin{aligned}
         t^{(A)}&:=t_{k+l,r+s}^{a_{k+l,r+s}}\cdots t_{k+l,1}^{a_{k+l,1}}\cdots t_{1,r+s}^{a_{1,r+s}}\cdots t_{11}^{a_{11}},\\
         \partial^{(B)}&:=\partial_{11}^{b_{11}}\partial_{12}^{b_{12}}\cdots \partial_{1,r+s}^{b_{1,r+s}}\cdots \partial_{k+l,1}^{b_{k+l,1}} \cdots \partial_{k+l,r+s}^{b_{k+l,r+s}}.
     \end{aligned}
 \end{equation*}

\begin{proposition}\label{proposition: the basis of PD^gr}
The monomials    
\begin{equation}\label{equation: the basis of PD^gr}
    \begin{aligned}
         t^{(A)}\otimes\partial^{(B)},\quad A,B\in \Xi_{k|l,r|s}
    \end{aligned}
\end{equation}
form a $\mathbb{K}$-basis of $\mathscr{PD}^{k|l,gr}_{r|s}$.
\end{proposition}
\begin{proof}
  By the definition  of  $\mathscr{PD}^{k|l,gr}_{r|s}$, it is isomorphic to $\mathscr{P}^{k|l}_{r|s}\otimes\mathscr{D}^{k|l}_{r|s}$  as a vector space. Then we can get a $\mathbb{K}$-basis of $\mathscr{P}^{k|l}_{r|s}\otimes\mathscr{D}^{k|l}_{r|s}$ by  Proposition~\ref{proposition: the basis of P and D}, which is  corresponding to (\ref{equation: the basis of PD^gr}) via the isomorphism of superalgebras
$\mathscr{PD}^{k|l,gr}_{r|s}\cong T(\mathscr{P}^{k|l,(1)}_{r|s}\oplus \mathscr{D}^{k|l,(1)}_{r|s})/I$ induced by the homomorphism $\zeta$ given in (\ref{equation: an homomorphism from T(P+D) to PD}).
\end{proof}

We next find the PBW basis of $\mathscr{PD}^{k|l}_{r|s}$ by using the  Bergman's Diamond Lemma \cite[Theorem 1.2]{B78}.
\begin{proposition}\label{proposition: the basis of PD}
The monomials    
\begin{equation}\label{equation: the basis of PD}
    \begin{aligned}
         t^{(A)}\partial^{(B)},\quad A,B\in \Xi_{k|l,r|s}
    \end{aligned}
\end{equation}
form a $\mathbb{K}$-basis of $\mathscr{PD}^{k|l}_{r|s}$.
\end{proposition}

\begin{proof}
   We can express any elements of  $\mathscr{PD}^{k|l}_{r|s}$ be a linear combination of (\ref{equation: the basis of PD}) by using (\ref{equation: the first relation for PD})-(\ref{equation: the fourth relation for PD}). 
Let
\[
X=\{t_{a,b},\partial_{a,b}
\mid a\in \mathbf{I}_{k|l},\ b\in \mathbf{I}_{r|s}\},
\]
,$\langle X\rangle$ the semigroup with $1$ on $X$ and $\mathbb K\langle X\rangle$  the free
associative $\mathbb K$-algebra on $X$.
Let $W=\{W_\sigma\}$(respectively, $F=\{f_{\sigma}\}$) denote the set of elements on the left-hand(respectively, right-hand) sides of the relations  (\ref{equation:relations for D}),(\ref{equation:relations for P}),(\ref{equation: the first relation for PD}-\ref{equation: the fourth relation for PD}).                 .
We define a reduction system
\begin{equation*}
    \begin{aligned}
        S=\{\sigma=(W_\sigma,f_\sigma)|W_\sigma\in W,f_{\sigma}\in F\}
    \end{aligned}
\end{equation*}
on $\mathbb K\langle X\rangle$.
Thus,$\mathscr{PD}_{r|s}^{k|l}
\cong
\mathbb K\langle X\rangle/I$,
where $I$ is the two-sided ideal generated by
$W_\sigma-f_\sigma$ for $\sigma\in S$.

 In order to use the  Bergman's Diamond Lemma \cite[Theorem 1.2]{B78}, we first define a semigroup partial ordering
on $\langle X\rangle$ which is compatible with
$S$ and has descending chain condition.

Equip $X$ with the total ordering $\prec$
determined by
\begin{equation*}
    \begin{aligned}
    &t_{a,b}\prec\partial_{c,d},\quad\forall a,c\in \mathbf{I}_{k|l},b,d\in \mathbf{I}_{r|s},\\
    &t_{a,b}\prec t_{c,d}
    \iff (a,b)>_{\mathrm{lex}}(c,d),\\
    &\partial_{a,b}\prec\partial_{c,d}
    \iff (a,b)<_{\mathrm{lex}}(c,d),
    \end{aligned}
\end{equation*}
where $<_{\mathrm{lex}}$ denotes the lexicographic order. 
Next, we extend $\prec$ to a
semigroup partial ordering $\prec$ on $\langle X\rangle$. For any $u=x_1\cdots x_n$, $v=y_1\cdots y_m$, we set $u\prec v$ if $n<m$ or if $m=n$ but for some $i\in \mathbb{N}$, $x_j=y_j$ ($j<i$) and $x_i \prec y_i$.

Since $X$ is finite, this is a semigroup
ordering with descending chain condition.

It follows directly from (\ref{equation:relations for D}),(\ref{equation:relations for P}),(\ref{equation: the first relation for PD}-\ref{equation: the fourth relation for PD}) that every monomial occurring
with nonzero coefficient in $f_\sigma$ is
strictly smaller than $W_\sigma$.
Hence $\prec$ is compatible with $S$.

We next show that all ambiguities of $S$ are resolvable.
Since every $W_\sigma$ has length two and
each reducible pair has a unique reduction
formula, there are no inclusion ambiguities.
The overlap ambiguities are precisely those
arising from monomials of the following forms:
\begin{equation*}
    \begin{aligned}
        t_{a,u}t_{b,v}t_{c,w},\qquad
\partial_{a,u}\partial_{b,v}\partial_{c,w},\qquad
\partial_{a,u}t_{b,v}t_{c,w},\qquad
\partial_{a,u}\partial_{b,v}t_{c,w}.
    \end{aligned}
\end{equation*}
The ambiguities involving only the $t_{i,j}$
or only the $\partial_{i,j}$ are resolvable
by Proposition \ref{proposition: the basis of P and D} and  \cite[Theorem 1.2]{B78}.
It remains to show that the overlap
ambiguities arising from the monomials
\[
\partial_{a,u}t_{b,v}t_{c,w}
\quad\text{and}\quad
\partial_{a,u}\partial_{b,v}t_{c,w},
\]
are resolvable.
We present a detailed calculation verifying the resolvability of an overlap ambiguity in one particular case. 

Let $(\sigma,\tau,\partial_{a,u},t_{a,v},t_{a,v})$ be the overlap ambiguity such that $W_{\sigma}=\partial_{a,u}t_{a,v}$,  $W_{\tau}=t_{a,v}t_{a,v}$, where $u<v$, $[a]+[v]=\bar{1}$. It follows from (\ref{equation:relations for P}) that
\begin{equation*}
    \begin{aligned}
        \partial_{a,u}(t_{a,v}t_{a,v})=0.
    \end{aligned}
\end{equation*}
On the other hand, by the relations  $[a]+[v]=\bar{1}$, (\ref{equation:relations for P}), (\ref{equation: the first relation for PD}) and (\ref{equation: the second relation for PD}), we have

\begin{equation*}
    \begin{aligned}
        (\partial_{a,u}t_{a,v})t_{a,v}=&q_a(-1)^{[a]+[u]}t_{a,v}(\partial_{a,u}t_{a,v})
 +(q-q^{-1})\sum_{h>a}(-1)^{[h]+([a]+[u])([v]+[h])}t_{h,v}(\partial_{h,u}t_{a,v})\\
=&q_a^{2}(t_{a,v}^{2})\partial_{a,u}
 +q_a(q-q^{-1})\sum_{h>a}(-1)^{[a]+[u]+[h]+([a]+[u])([v]+[h])}
 (t_{a,v}t_{h,v})\partial_{h,u}\\
&\quad +(q-q^{-1})\sum_{h>a}(-1)^{[h]+([a]+[u])([v]+[h])+[h]+[u]}
 t_{h,v}t_{a,v}\partial_{h,u}\\
=&q_aq_v(q-q^{-1})\sum_{h>a}
 (-1)^{[a]+[u]+[h]+([a]+[u])([v]+[h])+[h]+[v]}
 t_{h,v}t_{a,v}\partial_{h,u}\\
&\quad +(q-q^{-1})\sum_{h>a}(-1)^{[h]+([a]+[u])([v]+[h])+[h]+[u]}
 t_{h,v}t_{a,v}\partial_{h,u}\\
={}&(q-q^{-1})\sum_{h>a}
 (-1)^{[h]+([a]+[u])([v]+[h])+[h]+[u]}
 \Bigl(q_aq_v(-1)^{[a]+[v]}+1\Bigr)t_{h,v}t_{a,v}\partial_{h,u}\\
=&0=\partial_{a,u}(t_{a,v}t_{a,v}).
    \end{aligned}
\end{equation*}

Thus, the overlap ambiguity $(\sigma,\tau,\partial_{a,u},t_{a,v},t_{a,v})$ is resolvable.

To verify the resolvability of all overlap ambiguities, we classify these ambiguities  arising from the monomials $
\partial_{a,u}\partial_{b,v}t_{c,w}$ according to
the relative orders of their row and column
indices. If $a>b$, there are five possible relative orders of $a,b,c$, namely,
\begin{equation*}
    \begin{aligned}
        c>a>b,\quad a=c>b, \quad a>c>b,\quad a>b=c,\quad a>b>c.
    \end{aligned}
\end{equation*}
We need to consider these five cases separately in our computations.
In total, we need to consider $89$ possible relative orderings of the indices $(a,b,c,u,v,w)$. In total, we need to consider $89$ possible relative orderings of the indices $(a,b,c,u,v,w)$.
Similarly, we classify these ambiguities arising from the monomials $
\partial_{a,u}t_{b,v}t_{c,w}$ according to
the relative orders of their row and column
indices. There are also $89$ cases. Thus the computations are rather extensive, we develop a Python program to verify the resolvability of these overlap ambiguities. The Python programs used for the symbolic verification of all overlap ambiguities are available at the following GitHub repository: \url{https://github.com/Wangqq33/PBW_Basis_verification}.
 We emphasize that the verification is entirely symbolic, with all indices and summations treated symbolically, rather than being restricted to specific cases with small values of (k,l,r,s). 

\end{proof}

\begin{remark}
    We also tested replacing the constant term with \(1\) in (\ref{equation: the fourth relation for PD}), which resulted in unresolved overlap ambiguities. Thus, the chosen constant term is essential for resolving the overlap ambiguities and establishing the PBW basis.
\end{remark}

Let $p,t,p',t',r,s$ be nonnegative integers and set $k=\mathrm{max}\{p,p'\}, l=\mathrm{max}\{t,t'\}$. For $a\leq k,b\leq l$, we set $\tilde{\mathbf{I}}_{a|b,k|l}=\{k+l-a-b+1,\cdots,k+l\}$. We define $\mathscr{A}_{p|t,p'|t',r|s}$(respectively, $\mathscr{A}_{p|t,p'|t',r|s}^{gr}$) be the subalgebra of  $\mathscr{PD}^{k|l}_{r|s}$(respectively,  $\mathscr{PD}^{k|l,gr}_{r|s}$) generated by $t_{i,j}$ and $\partial_{i',j}$ (respectively, $t_{i,j}\otimes 1$ and $1\otimes\partial_{i',j}$), where $i\in \tilde{\mathbf{I}}_{p|t,k|l}$, $i'\in\tilde{\mathbf{I}}_{p'|t',k|l}$, $j\in \mathbf{\mathbf{I}}_{r|s}$.
\begin{proposition}
    We set $\tilde{\Xi}_{p|t,r|s}^{k|l}:=\{A=(a_{i,j})\in\Xi_{k|l,r|s}\mid a_{i,j}=0 \ \mathrm{if} \ i\notin  \tilde{\mathbf{I}}_{p|t,k|l}\}$, then the monomials
    \begin{equation}\label{equation: the basis of PD}
    \begin{aligned}
         t^{(A)}\partial^{(B)}(\mathrm{respectively, t^{(A)}\otimes\partial^{(B)}}),\quad A\in \tilde{\Xi}_{p|t,r|s}^{k|l},\ B\in \tilde{\Xi}_{p'|t',r|s}^{k|l}
    \end{aligned}
\end{equation}
form a $\mathbb{K}$-basis of  $\mathscr{A}_{p|t,p'|t',r|s}$(respectively,  $\mathscr{A}_{p|t,p'|t',r|s}^{gr}$).
\end{proposition}
\begin{proof}
    We prove the assertion only for $\mathscr{A}_{p|t,p'|t',r|s}$, as the proof for $\mathscr{A}_{p|t,p'|t',r|s}^{gr}$ is analogous. By Proposition \ref{proposition: the basis of PD}, it suffices to show that $\mathscr{A}_{p|t,p'|t',r|s}= \mathrm{span}\{ t^{(A)}\partial^{(B)},A\in \tilde{\Xi}_{p|t,r|s}^{k|l},\ B\in \tilde{\Xi}_{p'|t',r|s}^{k|l}\}$. Take an arbitrary product of the generators of $\mathscr{A}_{p|t,p'|t',r|s}$, by using the relations (\ref{equation:relations for D}),(\ref{equation:relations for P}),(\ref{equation: the first relation for PD})-(\ref{equation: the fourth relation for PD}), we can express it as a linear combination of $t^{(A)}\partial^{(B)},\ A,B\in \Xi_{k|l,r|s}$. By observe these relations, we see that only terms with weakly larger row or column indices can arise. Thus, $\mathscr{A}_{p|t,p'|t',r|s}\subseteq\mathrm{span}\{ t^{(A)}\partial^{(B)},A\in \tilde{\Xi}_{p|t,r|s}^{k|l},\ B\in \tilde{\Xi}_{p'|t',r|s}^{k|l}\}$. The reverse inclusion is immediate.
\end{proof}

\subsection{The module superalgebra structure of  $\mathscr{PD}^{k|l,gr}_{r|s}$ and $\mathscr{PD}^{k|l}_{r|s}$}

Next, we introduce the $\mathbf{U}_q(\mathfrak{gl}(k|l))\otimes\mathbf{U}_q(\mathfrak{gl}(r|s))$-module superalgebra structure of $\mathscr{PD}^{k|l,gr}_{r|s}$ and $\mathscr{PD}^{k|l}_{r|s}$. We note that $\mathscr{P}^{k|l,(1)}_{r|s}\oplus \mathscr{D}^{k|l,(1)}_{r|s}$ is $\mathbf{U}_q(\mathfrak{gl}(k|l))\otimes\mathbf{U}_q(\mathfrak{gl}(r|s))$-stable by (\ref{equation: formulas of the action on D}) and (\ref{equation: formulas of the action on P}). The $\mathbf{U}_q(\mathfrak{gl}(k|l))\otimes\mathbf{U}_q(\mathfrak{gl}(r|s))$-module structure on $\mathscr{P}^{k|l,(1)}_{r|s}\oplus \mathscr{D}^{k|l,(1)}_{r|s}$ can canonically induce a $\mathbf{U}_q(\mathfrak{gl}(k|l))\otimes\mathbf{U}_q(\mathfrak{gl}(r|s))$-module superalgebra structure on $T(\mathscr{P}^{k|l,(1)}_{r|s}\oplus \mathscr{D}^{k|l,(1)}_{r|s})$. Moreover, it can induce $\mathbf{U}_q(\mathfrak{gl}(k|l))\otimes\mathbf{U}_q(\mathfrak{gl}(r|s))$-module superalgebra structures on $\mathscr{PD}^{k|l,gr}_{r|s}$ and $\mathscr{PD}^{k|l}_{r|s}$ via the isomorphisms by the following lemma. It is easy to check that the structures of the $\mathbf{U}_q(\mathfrak{gl}(k|l))\otimes\mathbf{U}_q(\mathfrak{gl}(r|s))$-module of $\mathscr{PD}^{k|l,gr}_{r|s}$ given in the Remark \ref{remark: the module structure of twisted tensor product} coincide with those given above. 
\begin{lemma}
   The ideals $I$ and $I_{\psi}$ of $T(\mathscr{P}^{k|l,(1)}_{r|s}\oplus \mathscr{D}^{k|l,(1)}_{r|s})$ are $\mathbf{U}_q(\mathfrak{gl}(k|l))\otimes\mathbf{U}_q(\mathfrak{gl}(r|s))$-stable.
\end{lemma}
\begin{proof}
    We just prove it for $I_{\psi}$. Let $I_\mathscr{P}$ and $I_\mathscr{D}$ be the ideals of $T(\mathscr{P}^{k|l,(1)}_{r|s}\oplus \mathscr{D}^{k|l,(1)}_{r|s})$ generated by (\ref{equation:relations for D}) and (\ref{equation:relations for P}) respectively. Then $\mathscr{P}^{k|l}_{r|s}\cong T(\mathscr{P}^{k|l,(1)}_{r|s})/I_{\mathscr{P}}$ and $\mathscr{D}^{k|l}_{r|s}\cong T(\mathscr{D}^{k|l,(1)}_{r|s})/I_{\mathscr{D}}$. The $\mathbf{U}_q(\mathfrak{gl}(k|l))\otimes\mathbf{U}_q(\mathfrak{gl}(r|s))$-module structures on $\mathscr{P}^{k|l,(1)}_{r|s}$ and $\mathscr{D}^{k|l,(1)}_{r|s}$ can canonically induce $\mathbf{U}_q(\mathfrak{gl}(k|l))\otimes\mathbf{U}_q(\mathfrak{gl}(r|s))$-module superalgebra structures on $ T(\mathscr{P}^{k|l,(1)}_{r|s})$ and $ T(\mathscr{D}^{k|l,(1)}_{r|s})$ which are compatible with the  $\mathbf{U}_q(\mathfrak{gl}(k|l))\otimes\mathbf{U}_q(\mathfrak{gl}(r|s))$-module structure on $T(\mathscr{P}^{k|l,(1)}_{r|s}\oplus \mathscr{D}^{k|l,(1)}_{r|s})$. Since the canonical projections $T(\mathscr{P}^{k|l,(1)}_{r|s})\xrightarrow[]{} \mathscr{P}^{k|l}_{r|s}$ and $T(\mathscr{D}^{k|l,(1)}_{r|s})\xrightarrow[]{} \mathscr{D}^{k|l}_{r|s}$ are $\mathbf{U}_q(\mathfrak{gl}(k|l))\otimes\mathbf{U}_q(\mathfrak{gl}(r|s))$-module homomorphism,
    $I_{\mathscr{P}}$ and $I_{\mathscr{D}}$ are  $\mathbf{U}_q(\mathfrak{gl}(k|l))\otimes\mathbf{U}_q(\mathfrak{gl}(r|s))$-stable.

    Then, it suffices to show that the ideal generated by (\ref{equation: the definition of PD}) is $\mathbf{U}_q(\mathfrak{gl}(k|l))\otimes\mathbf{U}_q(\mathfrak{gl}(r|s))$-stable. We shall note that  for any $x\in\mathbf{U}_q(\mathfrak{gl}(k|l))\otimes\mathbf{U}_q(\mathfrak{gl}(r|s))$, $f\in \mathscr{P}^{k|l,(1)}_{r|s}$ and $g\in \mathscr{D}^{k|l,(1)}_{r|s}$ we have 
    \begin{equation*}
        \begin{aligned}
            \epsilon(x)\psi (g,f)=\sum (-1)^{[g][x_2]}\psi(x_1\cdot g,x_2\cdot f).
        \end{aligned}
    \end{equation*}
    Then,
    \begin{equation*}
        \begin{aligned}
            &x\cdot(  gf-m\circ\check{R}'_{ \mathscr{D}^{k|l}_{r|s},\mathscr{P}^{k|l}_{r|s}}(g\otimes f)-\psi (g,f))\\
            &=\sum(-1)^{[x_2][g]}((x_1\cdot g)(x_2\cdot f)-m\circ\check{R}'_{ \mathscr{D}^{k|l}_{r|s},\mathscr{P}^{k|l}_{r|s}}(x_1\cdot g\otimes x_2\cdot f)-\psi(x_1\cdot g,x_2\cdot f)),
        \end{aligned}
    \end{equation*}
    since  $x$ act on $\mathbb{K}$ by  $\epsilon(x)$. The right hand side of the equation is a sum of (\ref{equation: the definition of PD}). Therefore, $I_{\psi}$ is $\mathbf{U}_q(\mathfrak{gl}(k|l))\otimes\mathbf{U}_q(\mathfrak{gl}(r|s))$-stable.    
\end{proof}

By $\mathscr{PD}^{k|l,gr}_{r|s}\cong T(\mathscr{P}^{k|l,(1)}_{r|s}\oplus \mathscr{D}^{k|l,(1)}_{r|s})/I$ and $\mathscr{PD}^{k|l}_{r|s}\cong T(\mathscr{P}^{k|l,(1)}_{r|s}\oplus \mathscr{D}^{k|l,(1)}_{r|s})/I_{\psi}$, we have the following proposition.
\begin{proposition}\label{proposition: PD^gr and PD are module superalgs}
The superalgebras $\mathscr{PD}^{k|l,gr}_{r|s}$ and $\mathscr{PD}^{k|l}_{r|s}$ are $\mathbf{U}_q(\mathfrak{gl}(k|l))\otimes\mathbf{U}_q(\mathfrak{gl}(r|s))$-module superalgebras.
\end{proposition} 
Next we show that $\mathscr{PD}^{k|l,gr}_{r|s}\cong \mathscr{PD}^{k|l}_{r|s}$ as   $\mathbf{U}_q(\mathfrak{gl}(k|l))\otimes\mathbf{U}_q(\mathfrak{gl}(r|s))$-modules. We consider the linear map $P$ defined by
\begin{equation}\label{equation: the module isomorphism between PDgr and PD}
    \begin{aligned}
        P:\mathscr{PD}^{k|l,gr}_{r|s}&\xrightarrow[]{}\mathscr{PD}^{k|l}_{r|s},\\
        f\otimes g&\mapsto fg \quad\mathrm{for}~ f\in \mathscr{P}^{k|l}_{r|s},~ g\in \mathscr{D}^{k|l}_{r|s}.
    \end{aligned}
\end{equation}
\begin{proposition}
    The map $P$ is an isomorphism of  $\mathbf{U}_q(\mathfrak{gl}(k|l))\otimes\mathbf{U}_q(\mathfrak{gl}(r|s))$-modules.
\end{proposition}
\begin{proof}
    Since $\mathscr{PD}^{k|l,gr}_{r|s}$ and $\mathscr{PD}^{k|l}_{r|s}$ are     $\mathbf{U}_q(\mathfrak{gl}(k|l))\otimes\mathbf{U}_q(\mathfrak{gl}(r|s))$-module superalgebras and their module structure are both induced by the  $\mathbf{U}_q(\mathfrak{gl}(k|l))\otimes\mathbf{U}_q(\mathfrak{gl}(r|s))$-module structure on  $\mathscr{P}^{k|l,(1)}_{r|s}\oplus \mathscr{D}^{k|l,(1)}_{r|s}$, $P$ is a homomorphism of $\mathbf{U}_q(\mathfrak{gl}(k|l))\otimes\mathbf{U}_q(\mathfrak{gl}(r|s))$-module. It follows from Proposition \ref{proposition: the basis of PD^gr} and \ref{proposition: the basis of PD} that $P$ maps a basis of $\mathscr{PD}^{k|l,gr}_{r|s}$ to a basis of $\mathscr{PD}^{k|l}_{r|s}$. Therefore, $P$ is an isomorphism.
\end{proof}

Let 
\begin{equation*}
    \begin{aligned}
        (\mathscr{PD}^{k|l,gr}_{r|s})_{(\epsilon_R)}:=\{ f\in \mathscr{PD}^{k|l,gr}_{r|s}\mid x\cdot f= \epsilon_R(x) f, \ \forall x\in \mathbf{U}_q(\mathfrak{gl}(r|s)) \},
    \end{aligned}
\end{equation*}
 where $x\cdot f$ denotes $(1\otimes x)\cdot f$ and $\epsilon_R(x)$ denotes the counit of $\mathbf{U}_q(\mathfrak{gl}(r|s))$. We  define $(\mathscr{PD}^{k|l}_{r|s})_{(\epsilon_R)}$ similarly. Next we find elements  in $ (\mathscr{PD}^{k|l,gr}_{r|s})_{(\epsilon_R)}$ and $(\mathscr{PD}^{k|l}_{r|s})_{(\epsilon_R)}$. Let
\begin{equation}\label{equation: definition of L_i,j}
    \begin{aligned}
        L_{i,j}^{gr}&:=\sum_{c=1}^{r+s}(-1)^{[c]([i]+1)} t_{i,c} \partial_{j,c}, \quad R_{i,j}^{gr}:=\sum_{c=1}^{k+l}(-1)^{[c]([j]+1)}t_{c,i}\partial_{c,j} \in  \mathscr{PD}^{k|l,gr}_{r|s},\\
        L_{i,j}&:=\sum_{c=1}^{r+s}(-1)^{[c]([i]+1)} t_{i,c} \partial_{j,c},
        \quad  R_{i,j}:=\sum_{c=1}^{k+l}(-1)^{[c]([j]+1)}t_{c,i}\partial_{c,j}
        \in  \mathscr{PD}^{k|l}_{r|s},
    \end{aligned}
\end{equation}
 where $i,j\in \mathbf{I}_{k|l}$.
\begin{proposition}\label{proposition: Lijgr Lij are invariants}
    The elements $L_{i,j}^{gr} \in (\mathscr{PD}^{k|l,gr}_{r|s})_{(\epsilon_R)}$, $L_{i,j} \in (\mathscr{PD}^{k|l}_{r|s})_{(\epsilon_R)}$.
\end{proposition}
\begin{proof}
    We just provide the proof for $L_{ij}$ since the proof for $L_{ij}$ follows by exactly the same argument. It suffices to check $x\cdot L_{ij}=\epsilon_R(x)L_{ij}$ for any $x=E_a,F_a ,K_b\in \mathbf{U}_q(\mathfrak{gl}(r|s))$, $a\in \mathbf{I}'_{r|s}$, $b\in \mathbf{I}_{r|s}$. Since  $\mathscr{PD}^{k|l}_{r|s}$ is a  $\mathbf{U}_q(\mathfrak{gl}(k|l))\otimes\mathbf{U}_q(\mathfrak{gl}(r|s))$-module superalgebra, we have 
    \begin{equation*}
        \begin{aligned}
            &E_a\cdot L_{i,j}\\
            =&\sum_{c=1}^{r+s}((-1)^{[c]([i]+1)}(E_a\cdot t_{i,c})(K_a K_{a+1}^{-1}\cdot \partial_{j,c})+(-1)^{[c]([i]+1)+([i]+[c])([a]+[a+1])}t_{i,c}(E_a\cdot\partial_{j,c}))\\
            =&(-1)^{[a]([i]+1)}(E_a \cdot t_{i,a})(K_aK_{a+1}^{-1}\cdot\partial_{j,a})+(-1)^{[a+1]([i]+1)+([i]+[a+1])([a]+[a+1])}t_{i,a+1}(E_a\cdot\partial_{j,a+1})\\
            =&(-1)^{[a]([i]+1)} ((-q_a^{-1} t_{i,a+1})(q_a\partial_{j,a})+t_{i,a+1}\partial_{j,a})=0=\epsilon(E_a)L_{ij}.
        \end{aligned}
    \end{equation*}
    

Similarly, we have  $F_a\cdot L_{i,j}=0$ and $K_a\cdot L_{i,j}=L_{i,j}$.
\end{proof}

\begin{lemma}\label{lemma: use Lk+l,k+l to obtain Lij}
    The elements $L_{i,j}$ belongs to the $(\mathbf{U}_q(\mathfrak{gl}(k|l))\otimes1)$-submodule of   $\mathscr{PD}^{k|l}_{r|s}$  that is generated  by $L_{k+l,k+l}$. 
\end{lemma}
\begin{proof}
    Let $M$ be the $(\mathbf{U}_q(\mathfrak{gl}(k|l))\otimes1)$-submodule of   $\mathscr{PD}^{k|l}_{r|s}$  that is generated  by $L_{k+l,k+l}$. For simplicity of notation, we write $x\cdot f$ for $(x\otimes1)\cdot f$ in the following proof. Since  $\mathscr{PD}^{k|l}_{r|s}$ is a  $\mathbf{U}_q(\mathfrak{gl}(k|l))\otimes\mathbf{U}_q(\mathfrak{gl}(r|s))$-module superalgebra, for any $a\in \mathbf{I}_{k|l}$, we have 
    \begin{equation*}
        \begin{aligned}
             E_a\cdot L_{i,j}=\sum_{c=1}^{r+s}(-1)^{r_{i,c}}(E_a\cdot t_{i,c})(K_a K_{a+1}^{-1}\cdot \partial_{j,c})+(-1)^{r_{i,c}+([i]+[c])([a]+[a+1])}t_{i,c}(E_a\cdot\partial_{j,c})
        \end{aligned}
    \end{equation*}
    By (\ref{equation: formulas of the action on D}) and (\ref{equation: formulas of the action on P}), we have 
        \begin{align}
        \label{equation: Ea act on Li,a+1}
            E_a\cdot L_{i,a+1}&=(-1)^{([a]+[a+1])[i]}L_{i,a},&\mathrm{for}~i\neq a
        \end{align}
    Similarly, we have
        \begin{align}
        \label{equation: Fa act on La+1,j}
            F_a\cdot L_{a+1,j}&=-(-1)^{([a]+[a+1])[a+1]}q_aL_{a,j},&\mathrm{for}~j\neq a,\\
        \label{equation: Ea act on La,a+1}
            E_a\cdot L_{a,a+1}&=-q_a^{-1}q_{a+1}^{-1}L_{a+1,a+1}+L_{a,a},&\mathrm{otherwise}.
        \end{align}
    From (\ref{equation: Ea act on Li,a+1}) and (\ref{equation: Fa act on La+1,j}), 
    it follows that all $L_{i,j}$ with $i\neq j$ can be obtained from 
    $L_{k+l,k+l}$ by repeatedly applying $E_a$ and $F_a$. Thus, $L_{i,j}\in M$ for $i\neq j$. By (\ref{equation: Ea act on La,a+1}) and by induction on $a$, we obtain $L_{i,i}\in M$.
\end{proof}

\subsection{The $\mathscr{PD}^{k|l}_{r|s}$-module structure of $\mathscr{P}^{k|l}_{r|s}$}\label{the module structure of PD}

In this subsection, we will give the  $\mathscr{PD}^{k|l}_{r|s}$-module structure of $\mathscr{P}^{k|l}_{r|s}$ and show that $\mathscr{P}^{k|l}_{r|s}$ a faithful  $\mathscr{PD}^{k|l}_{r|s}$-module. Then we can view  $\mathscr{PD}^{k|l}_{r|s}$ as a subsuperalgebra of $\mathrm{End}_{\mathbb{K}}(\mathscr{P}^{k|l}_{r|s})$.
Let $\mathscr{I}$ be the left ideal of $\mathscr{PD}^{k|l}_{r|s}$ generated by $\partial_{a,b}$, $a\in \mathbf{I}_{k|l}$, $b\in \mathbf{I}_{r|s}$.
Then $\mathscr{P}^{k|l}_{r|s}\cong \mathscr{PD}^{k|l}_{r|s}/\mathscr{I}$ by Proposition \ref{proposition: the basis of P and D} and \ref{proposition: the basis of PD}. Then we define the action of $\mathscr{PD}^{k|l}_{r|s}$ on $\mathscr{P}^{k|l}_{r|s}$ as 
\begin{equation*}
    \begin{aligned}
        D\cdot f=Df\in \mathscr{P}^{k|l}_{r|s}\cong \mathscr{PD}^{k|l}_{r|s}/\mathscr{I},
    \end{aligned}
\end{equation*}
where $D\in \mathscr{PD}^{k|l}_{r|s}$ and $f\in\mathscr{P}^{k|l}_{r|s}$. It makes $\mathscr{P}^{k|l}_{r|s}$ a $\mathscr{PD}^{k|l}_{r|s}$-module.
Before proving that  $\mathscr{P}^{k|l}_{r|s}$ a faithful  $\mathscr{PD}^{k|l}_{r|s}$-module, we need some lemmas.
Recall that $\mathscr{P}^{k|l,(i)}_{r|s}$ is the homogeneous subspaces of degree $i$ of $\mathscr{P}^{k|l}_{r|s}$. Let  $\mathscr{P}^{k|l,(\leq i)}_{r|s}:=\oplus_{j=0}^i\mathscr{P}^{k|l,(j)}_{r|s}$ for $i\geq 0$.

\begin{lemma}\label{lemma: the property 1 of the action of partial_ij}
    We have $\partial_{a,b}\cdot \mathscr{P}^{k|l,(\leq i)}_{r|s} \subseteq  \mathscr{P}^{k|l,(\leq i-1)}_{r|s}$ for any $a\in \mathbf{I}_{k|l}$, $b\in \mathbf{I}_{r|s} $ and $i\in \mathbb{Z}_{>0}$.
\end{lemma}
\begin{proof}
   The proof proceeds by induction on $i$. Fix arbitrary $a\in \mathbf{I}_{k|l}$, $b\in \mathbf{I}_{r|s} $. For $i=1$, by the (\ref{equation: the first relation for PD})-(\ref{equation: the fourth relation for PD}), we have $\partial_{a,b}\cdot t_{a',b'}=(-1)^{[a]+[b]+[a][b]}\delta_{a,a'}\delta_{b,b'}\in \mathbb{K}$ for any $a'\in \mathbf{I}_{k|l}$, $b'\in \mathbf{I}_{r|s} $. Thus, it holds for $i=1$. 
    
    We assume that it holds for $i$. 
    It suffices to check that $\partial_{a,b}\cdot (t_{a_1,b_1}\cdots t_{a_{i+1},b_{i+1}})\in \mathscr{P}^{k|l,(\leq i)}_{r|s}$ for any $a_p\in \mathbf{I}_{k|l}$, $b_p\in \mathbf{I}_{r|s}$. 
    By using  the (\ref{equation: the first relation for PD})-(\ref{equation: the fourth relation for PD}), $\partial_{a,b}t_{a_1,b_1}$ can be expressed as a linear combination of the elements of $\mathscr{PD}^{k|l}_{r|s}$ as the form of $t_{p,q}\partial_{p',q'}$ or $1$. 
    Thus, we have $\partial_{a,b}\cdot (t_{a_1,b_1}\cdots t_{a_{i+1},b_{i+1}}) \in \mathscr{P}^{k|l,(\leq 1)}_{r|s} \cdot \mathscr{P}^{k|l,(\leq i-1)}_{r|s}\subseteq  \mathscr{P}^{k|l,(\leq i)}_{r|s} $ by the induction hypothesis.
\end{proof}

We define a total order $\prec$ for the pairs $(a,b)$ where $a\in \mathbf{I}_{k|l}$, $b\in \mathbf{I}_{r|s} $ such that $(a,b)\prec(c,d)$ if either $a+b<c+d$ or $a+b=c+d$ and $a<c$. For  $a\in \mathbf{I}_{k|l}$, $b\in \mathbf{I}_{r|s} $ and $i\in \mathbb{Z}_{>0}$, we set 
\begin{equation*}
    \begin{aligned}
        k_{a,b,i}=(-1)^{[a]+[b]+[a][b]}\sum_{j=0}^{i-1}(-1)^{j([a]+[b])}q_a^jq_b^j.
    \end{aligned}
\end{equation*}
\begin{lemma} \label{lemma: the property 2 of the action of partial_ij}
    If $i,p\in \mathbb{Z}_{>0}$ and $(a_j,b_j)\prec(a,b)$ for $1\leq j\leq p$, then 
    \begin{equation}
        \begin{aligned}
            \partial_{a,b}\cdot(t_{a,b}^it_{a_1,b_1}\cdots t_{a_p,b_p})=k_{a,b,i}t_{a,b}^{i-1}t_{a_1,b_1}\cdots t_{a_p,b_p}.
        \end{aligned}
    \end{equation}
\end{lemma}
\begin{proof}
    Firstly, we show that  if  $(a_j,b_j)\prec(\tilde{a},\tilde{b})$ for $1\leq j\leq p$ then 
    \begin{equation}\label{equation: the action of partial on what is 0}
        \begin{aligned}
            \partial_{\tilde{a},\tilde{b}}\cdot(t_{a_1,b_1}\cdots t_{a_p,b_p})=0
        \end{aligned}
    \end{equation}
     by induction on $p$. It suffices to show that  $\partial_{\tilde{a},\tilde{b}}t_{a_1,b_1}\cdots t_{a_p,b_p}\in \mathscr{I}$.  
    We note that $(a_1,b_1)\prec(\tilde{a},\tilde{b})$ implies that $a_1\neq \tilde{a}$ or $b_1\neq \tilde{b}$, so  $\partial_{\tilde{a},\tilde{b}}t_{a_1,b_1}\in \mathscr{I}$ by the (\ref{equation: the first relation for PD})-(\ref{equation: the third relation for PD}). Moreover, we have  
    \begin{equation*}
        \begin{aligned}
            \partial_{\tilde{a},\tilde{b}}t_{a_1,b_1}\cdots t_{a_p,b_p}=kt_{a_1,b_1}\partial_{\tilde{a},\tilde{b}} t_{a_2,b_2}\cdots t_{a_p,b_p}+\delta_{\tilde{a},a_1}\sum_{a'>\tilde{a}} k_{a'}t_{a',b_1}\partial_{a',\tilde{b}} t_{a_2,b_2}\cdots t_{a_p,b_p}\\
            +\delta_{\tilde{b},b_1}\sum_{b'>\tilde{b}}k_{b'}  t_{a_1,b'}\partial_{\tilde{a},b'}t_{a_2,b_2}\cdots t_{a_p,b_p},
        \end{aligned}
    \end{equation*}
    for some $k,k_{a'},k_{b'}\in \mathbb{K}$.  Since $a'>a$ implies $(a_j,b_j)\prec(a,b)\prec(a',b)$ and $b'>b$ implies $(a_j,b_j)\prec(a,b)\prec (a,b')$, all terms on the right-hand side of the equality lie in the left ideal $\mathscr{I}$  by inductive hypotheses.

    Next, we prove the lemma by induction on $i$. It follows from the relation \ref{equation: the fourth relation for PD} that 
    \begin{equation}
        \begin{aligned}
        &\partial_{a,b}t_{a,b}^it_{a_1,b_1}\cdots t_{a_p,b_p}\\
        =&
            (-1)^{[a]+[b]+[a][b]}t_{a,b}^{i-1}t_{a_1,b_1}\cdots t_{a_p,b_p}+q_aq_b t_{a,b}\partial_{a,b}t_{a,b}^{i-1}t_{a_1,b_1}\cdots t_{a_p,b_p}\cdot(-1)^{[a]+[b]}\\
        +&
        (q-q^{-1})\sum_{b'>b}q_at_{a,b'}\partial_{a,b'}t_{a,b}^{i-1}t_{a_1,b_1}\cdots t_{a_p,b_p}\cdot(-1)^{[b']+([a]+[b])([a]+[b'])}\\
        +&(q-q^{-1})\sum_{a'>a}q_bt_{a',b}\partial_{a',b}t_{a,b}^{i-1}t_{a_1,b_1}\cdots t_{a_p,b_p}\cdot(-1)^{[a']+([a]+[b])([b]+[a'])} \notag\\
        +&(q-q^{-1})^2\sum_{a'>a}\sum_{b'>b}t_{a',b'}\partial_{a',b'}t_{a,b}^{i-1}t_{a_1,b_1}\cdots t_{a_p,b_p}\cdot(-1)^{[a']+[b']+([a]+[b'])([b]+[a'])}.
        \end{aligned}
    \end{equation}
    For $i=1$, by (\ref{equation:  the action of partial on what is 0}), we have  $\partial_{a,b}t_{a,b}t_{a_1,b_1}\cdots t_{a_p,b_p}=(-1)^{[a]+[b]+[a][b]}
    t_{a_1,b_1}\cdots t_{a_p,b_p}+f$ for some $f\in \mathscr{I}$, so the lemma holds for $i=1$. For $i\geq 2$, by inductive hypotheses, we have  
    \begin{equation}\label{equation : somthing not important}
        \begin{aligned}
            \partial_{a,b}t_{a,b}^{i-1}t_{a_1,b_1}\cdots t_{a_p,b_p}=(-1)^{[a]+[b]+[a][b]}\sum_{j=0}^{i-2}(-1)^{j([a]+[b])}q_a^jq_b^jt_{a,b}^{i-2}t_{a_1,b_1}\cdots t_{a_p,b_p}+ f'
        \end{aligned}
    \end{equation}
    for some $f'\in \mathscr{I}$. By (\ref{equation:  the action of partial on what is 0}) and (\ref{equation : somthing not important}), we have 
    \begin{equation*}
        \begin{aligned}
            \partial_{a,b}t_{a,b}^{i}t_{a_1,b_1}\cdots t_{a_p,b_p}
            =&(-1)^{[a]+[b]+[a][b]}t_{a,b}^{i-1}t_{a_1,b_1}\cdots t_{a_p,b_p}\\
            &+(-1)^{[a]+[b]+[a][b]}\sum_{j=0}^{i-2}(-1)^{(j+1)([a]+[b])}q_a^{j+1}q_b^{j+1}t_{a,b}^{i-1}t_{a_1,b_1}\cdots t_{a_p,b_p}+ f''\\
            =&(-1)^{[a]+[b]+[a][b]}\sum_{j=0}^{i-1}(-1)^{j([a]+[b])}q_a^jq_b^jt_{a,b}^{i-1}t_{a_1,b_1}\cdots t_{a_p,b_p} + f''
        \end{aligned}
    \end{equation*}
     for some $f''\in \mathscr{I}$. Thus, $\partial_{a,b}\cdot(t_{a,b}^it_{a_1,b_1}\cdots t_{a_p,b_p})=k_{a,b,i}t_{a,b}^{i-1}t_{a_1,b_1}\cdots t_{a_p,b_p}$.
\end{proof}

\begin{proposition}\label{proposition: P is faithful PD module}
     The superalgebra $\mathscr{P}^{k|l}_{r|s}$ is a faithful $\mathscr{PD}^{k|l}_{r|s}$-module. 
\end{proposition}
\begin{proof}
    The proof is completely analogous to that of \cite[Proposition 3.11.4]{LSS24} by using Lemma \ref{lemma: the property 1 of the action of partial_ij} and \ref{lemma: the property 2 of the action of partial_ij}.
\end{proof}
Let $\phi_{PD}:\mathscr{PD}^{k|l}_{r|s}\xrightarrow[]{} \mathrm{End}_{\mathbb{K}}(\mathscr{P}^{k|l}_{r|s})$ be the homomorphism of superalgebras induced by the $\mathscr{PD}^{k|l}_{r|s}$-module structure of $\mathscr{P}^{k|l}_{r|s}$. 
Proposition \ref{proposition: P is faithful PD module} follows that $\phi_{PD}$
is injective, so we can identify $\mathscr{PD}^{k|l}_{r|s}$ with the subsuperalgebra $\phi_{PD}(\mathscr{PD}^{k|l}_{r|s})$ of $\mathrm{End}_{\mathbb{K}}(\mathscr{P}^{k|l}_{r|s})$.

\section{the subsuperalgebras $\mathscr{L}_{\bullet},\mathscr{R}_{\bullet}$ of $\mathrm{End}_{\mathbb{K}}(\mathscr{P}^{k|l}_{r|s})$} \label{Section:L,R}

\subsection{the definitions of $\mathscr{L}_{\bullet},\mathscr{R}_{\bullet}$ and the module structures of  $\mathrm{End}_{\mathbb{K}}(\mathscr{P}^{k|l}_{r|s})$}
Let $H$ be a Hopf superalgebra and vector superspace $V$ be a $H$-module. Let $\phi: H\xrightarrow[ ]{}\mathrm{End}(V) $ be  the homomorphism of superalgebras induced by the $H$-module structure of $V$. We define an action of $H$ on  $\mathrm{End}(V)$ by
\begin{equation}\label{equation: the module structure of End(V)}
    \begin{aligned}
        x\cdot F=\sum(-1)^{[x_2][F]}\phi(x_1)F\phi(S(x_2)),
    \end{aligned}
\end{equation}
for $x\in H$ and $F\in \mathrm{End}(V)$. Let $\mathrm{ad}_y(x):=\sum(-1)^{[x][y_2]}y_1xS(y_2)$ be the adjoint  action of $H$. We can check that 
\begin{equation}\label{equation: the relation of adjoint action and phi}
    \begin{aligned}
        \phi(\mathrm{ad}_y(x))=y\cdot \phi(x),\quad \forall x,y\in H.
    \end{aligned}
\end{equation}
We set
\begin{equation*}
    \begin{aligned}
        &\mathrm{End}(V)_{(\epsilon)}
        =
        \{F\in \mathrm{End}(V) \mid x\cdot F=\epsilon(x) F, \quad\forall x\in H\},\\
        &\mathrm{End}(V)^{\phi(H)}
        =
        \{F\in \mathrm{End}(V)\mid FG=(-1)^{[F][G]}GF,\quad \forall G\in\phi(H)\}
    \end{aligned}
\end{equation*}
the  invariant subspaces of $\mathrm{End}_{\mathbb{K}}(V)$. In fact, they are equal. 
\begin{lemma}\label{lemma: invariant subspaces of  End(V)}
    We have $\mathrm{End}(V)_{(\epsilon)}
    =
    \mathrm{End}(V)^{\phi(H)}$.
\end{lemma}
\begin{proof}
    For any  $F\in  \mathrm{End}_{\mathbb{K}}(V)^{\phi(H)}$, $x\in H$, we have
    \begin{equation*}
        \begin{aligned}
            x\cdot F
            &=\sum(-1)^{[x_2][F]}\phi_{U}(x_1)F\phi(S(x_2))\\
            &=\sum(-1)^{[x_2][F]}\phi_{U}(x_1)\phi(S(x_2))F\cdot(-1)^{[x_2][F]}\\
            &=\phi(\sum x_1S(x_2))F
            =\epsilon(x)F,
        \end{aligned}
    \end{equation*}
    where the last equality follows from the antipode axiom. Thus, $
    \mathrm{End}(V)^{\phi(H)}\subseteq\mathrm{End}(V)_{(\epsilon)}
    $. Conversely, for any $F\in \mathrm{End}(V)_{(\epsilon)}$, $x\in H$, the counit axiom, the antipode axiom, and the fact that $\epsilon(x_2)=0$ unless $[x_2]=0$ follows that
    \begin{equation*}
        \begin{aligned}
            \phi(x)F&=\sum \phi(x_1\epsilon(x_2))F\\
            &=\sum(-1)^{[x_2][F]}\phi(x_1)F\phi(\epsilon(x_2))\\
            &=\sum(-1)^{([x_2]+[x_3])[F]}\phi(x_1)F\phi(S(x_2))\phi(x_3)\\
            &=\sum(-1)^{[x_2][F]}(x_1\cdot F)\phi(x_2)\\
            &=\sum(-1)^{([x_1]+[x_2])[F]}\epsilon(x_1)F\phi(x_2)\\
            &=(-1)^{[x][F]}F\phi(x),
        \end{aligned}
    \end{equation*}
where the   second  equality follows from the fact that $\epsilon(x_2)=0$ unless $[x_2]=0$.
\end{proof}

Let $\phi_{U}:\mathbf{U}_q(\mathfrak{gl}(k|l))\otimes \mathbf{U}_q(\mathfrak{gl}(r|s))\xrightarrow[]{} \mathrm{End}_{\mathbb{K}}(\mathscr{P}^{k|l}_{r|s})$ be the homomorphism of superalgebras induced by the $\mathbf{U}_q(\mathfrak{gl}(k|l))\otimes \mathbf{U}_q(\mathfrak{gl}(r|s))$-module structure of $\mathscr{P}^{k|l}_{r|s}$. 
We set
\begin{equation*}
    \begin{aligned}
        \mathscr{L}:=\phi_{U}(\mathbf{U}_q(\mathfrak{gl}(k|l))\otimes 1),\quad \mathscr{R}:=\phi_{U}(1\otimes \mathbf{U}_q(\mathfrak{gl}(r|s))).\\
    \end{aligned}
\end{equation*}
By identifying  $\mathscr{PD}^{k|l}_{r|s}$ with the  subsuperalgebra $\phi_{PD}(\mathscr{PD}^{k|l}_{r|s})$ of $\mathrm{End}_{\mathbb{K}}(\mathscr{P}^{k|l}_{r|s})$ ,we set 
\begin{equation}\label{definition: L. and R.}
    \begin{aligned}
        \mathscr{L}_{\bullet}:=\mathscr{L}\cap \mathscr{PD}^{k|l}_{r|s},\quad \mathscr{R}_{\bullet}:=\mathscr{R}\cap \mathscr{PD}^{k|l}_{r|s}.
    \end{aligned}
\end{equation}
Note that $\mathrm{End}_{\mathbb{K}}(\mathscr{P}^{k|l}_{r|s})$ is equipped with $\mathbf{U}_q(\mathfrak{gl}(k|l))\otimes \mathbf{U}_q(\mathfrak{gl}(r|s))$-module structure by (\ref{equation: the module structure of End(V)}).

\begin{lemma}\label{lemma: relations between invariant subpaces}
    We have    $(\mathscr{PD}^{k|l}_{r|s})_{(\epsilon_L)}=(\mathscr{PD}^{k|l}_{r|s})^{\mathscr{L}}$  and $(\mathscr{PD}^{k|l}_{r|s})_{(\epsilon_R)}=(\mathscr{PD}^{k|l}_{r|s})^{\mathscr{R}}$.
\end{lemma}
\begin{proof}
 By applying Lemma \ref{lemma: invariant subspaces of  End(V)} to $\mathscr{P}^{k|l}_{r|s}$, we have $\mathrm{End}_{\mathbb{K}}(\mathscr{P}^{k|l}_{r|s})_{(\epsilon_L)}=\mathrm{End}_{\mathbb{K}}(\mathscr{P}^{k|l}_{r|s})^{\mathscr{L}}$ and  $\mathrm{End}_{\mathbb{K}}(\mathscr{P}^{k|l}_{r|s})_{(\epsilon_R)}=\mathrm{End}_{\mathbb{K}}(\mathscr{P}^{k|l}_{r|s})^{\mathscr{R}}$ , where $\epsilon_L$(respectively, $\epsilon_R$) is the counit of $\mathbf{U}_q(\mathfrak{gl}(k|l))$ (respectively, $\mathbf{U}_q(\mathfrak{gl}(r|s))$).
Since $\mathscr{PD}^{k|l}_{r|s}$ is a subsuperalgebra of  $\mathrm{End}_{\mathbb{K}}(\mathscr{P}^{k|l}_{r|s})$, the action of $\mathbf{U}_q(\mathfrak{gl}(k|l))\otimes \mathbf{U}_q(\mathfrak{gl}(r|s))$ on $\mathrm{End}_{\mathbb{K}}(\mathscr{P}^{k|l}_{r|s})$  induces an action on $\mathscr{PD}^{k|l}_{r|s}$, which coincides with that defined in \S \ref{subsection: module structure of P,D}.    
\end{proof}

We set 
\begin{equation*}
    \begin{aligned}
        \mathring{U}_L:=\{x\in \mathbf{U}_q(\mathfrak{gl}(k|l))\mid \phi_U(x\otimes1)\in \mathscr{PD}^{k|l}_{r|s}\},
        \
        \mathring{U}_R:=\{x\in \mathbf{U}_q(\mathfrak{gl}(r|s))\mid \phi_U(1\otimes x)\in \mathscr{PD}^{k|l}_{r|s}\}.
    \end{aligned}
\end{equation*}
\begin{lemma}\label{lemma: invariant under the adjoint action}
    For any $y\in\mathbf{U}_q(\mathfrak{gl}(k|l))$, we have  $\mathrm{ad}_{y}(\mathring{U}_L)\subseteq \mathring{U}_L $ .
\end{lemma}
\begin{proof}
    It follows from (\ref{equation: the relation of adjoint action and phi}) that 
    \begin{equation*}
        \begin{aligned}
            \phi_U(\mathrm{ad}_y(x)\otimes 1)=\phi_U(\mathrm{ad}_{y\otimes 1}(x\otimes 1))=(y\otimes 1)\cdot\phi_U(x\otimes 1)\in\mathscr{PD}^{k|l}_{r|s},
        \end{aligned}
    \end{equation*}
    where $x\in \mathring{U}_L$, $y\in \mathbf{U}_q(\mathfrak{gl}(k|l))$. 
\end{proof}

\subsection{Some special elements in $\mathscr{L}_{\bullet}$}

\begin{lemma}\label{lemma: partial_ab act on products of some t_ij is zero}
    Assume that either $a\notin \{a_1,\cdots a_p\}$ or $b\notin \{ b_1,\cdots , b_p\}$, then 
    \begin{equation*}
        \begin{aligned}
            \partial_{a,b}\cdot(t_{a_1,b_1}\cdots t_{a_p,b_p})=0.
        \end{aligned}
    \end{equation*}
\end{lemma}
\begin{proof}
    The proof proceeds by the induction on $p$. It follows from (\ref{equation: the first relation for PD})-(\ref{equation: the third relation for PD}) that it holds for $p=1$.
    
    We now consider the general case, if $a\neq a_1$ and $b\neq b_1$ then
    \begin{equation*}
        \begin{aligned}
            \partial_{a,b}t_{a_1,b_1}\cdots t_{a_p,b_p}=(-1)^{([a]+[b])([a_1]+[b_1])}t_{a_1,b_1}\partial_{a,b}t_{a_2,b_{2}}\cdots t_{a_p,b_p}.
        \end{aligned}
    \end{equation*}
    It follows from  the induction hypothesis that the right-hand side of the equality lies in  $\mathscr{I}$, which means $\partial_{a,b}\cdot(t_{a_1,b_1}\cdots t_{a_p,b_p})=0$ in this case.
    The proof for the case $a=a_1$ and the case $b=b_1$ is similar, so we just give the proof for  one of them.
    If $a= a_1$ ,then $b\notin \{ b_1,\cdots , b_p\}$. By (\ref{equation: the second relation for PD}), we have
     \begin{equation*}
        \begin{aligned}
            \partial_{a,b}t_{a_1,b_1}\cdots t_{a_p,b_p}=(-1)^{([a]+[b])([a_1]+[b_1])}q_at_{a_1,b_1}\partial_{a,b}t_{a_2,b_{2}}\cdots t_{a_p,b_p}\\+\sum_{a'>a}c_{a'}t_{a',b_1}\partial_{a',b}t_{a_2,b_{2}}\cdots t_{a_p,b_p},
        \end{aligned}
    \end{equation*}
    where $c_{a'}\in \mathbb{K}$. Thus, the statement holds for this case by the induction hypothesis. 
\end{proof}

\begin{lemma}\label{lemma: the property 3 of the action of PD on P}
    Assume that $\partial_{a,b} \cdot f=f'$ for some $a\in \mathbf{I}_{k|l}$, $b\in \mathbf{I}_{r|s} $, $f,f'\in \mathscr{P}^{k|l}_{r|s}$, then we have 
    \begin{equation*}
        \begin{aligned}
            \partial_{a,b}\cdot(f t_{a_1,b_1}\cdots t_{a_p,b_{p}})=f' t_{a_1,b_1}\cdots t_{a_p,b_{p}}
        \end{aligned}
    \end{equation*}
    for any $a_j\in \mathbf{I}_{k|l}$, $b_j\in \mathbf{I}_{r|s}$ that satisfy either $a>\mathrm{max}\{a_j\}_{j=1}^p$ or $b>\mathrm{max}\{b_j\}_{j=1}^p$.
\end{lemma}
\begin{proof}
    By the definition of  $\mathscr{I}$ and (\ref{equation: the first relation for PD})-(\ref{equation: the fourth relation for PD}), we have 
    \begin{equation*}
        \begin{aligned}
            \partial_{a,b} f=f'+\sum_{a'\geq a,b'\geq b}g_{a',b'}\partial_{a',b'}\in \mathscr{PD}^{k|l}_{r|s}
        \end{aligned}
    \end{equation*}
    for some $g_{a',b'}\in \mathscr{PD}^{k|l}_{r|s}$. Then we have 
    \begin{equation*}
        \begin{aligned}
            \partial_{a,b} ft_{a_1,b_1}\cdots t_{a_p,b_{p}}=f't_{a_1,b_1}\cdots t_{a_p,b_{p}}+\sum_{a'\geq a,b'\geq b}g_{a',b'}\partial_{a',b'}t_{a_1,b_1}\cdots t_{a_p,b_{p}}.
        \end{aligned}
    \end{equation*}
    We notice that  either $a'\notin \{a_1,\cdots a_p\}$ or $b'\notin \{ b_1,\cdots , b_p\}$ , since either $a'\geq a>\mathrm{max}\{a_j\}_{j=1}^p$ or $b' \geq b>\mathrm{max}\{b_j\}_{j=1}^p$. Then it follows from Lemma   \ref{lemma: partial_ab act on products of some t_ij is zero} that the second term on the right-hand side of the equality lies in $\mathscr{I}$.
\end{proof}

\begin{lemma}\label{lemma: the action of the element D on products of some tij}
    Let $D_0:=1+(q_{k+l}^2-1)\sum_{c=1}^{r+s}(-1)^{[k+l]+[c]+[k+l][c]}t_{k+l,c}\partial_{k+l,c}$, then for any $a_j\in \mathbf{I}_{k|l}$,$b_j\in \mathbf{I}_{r|s}$, we have
    \begin{equation*}
        \begin{aligned}
            D_0\cdot(t_{a_1,b_1}\cdots t_{a_p,b_p})=q_{k+l}^{2\sum_{j=1}^p \delta_{k+l,a_j}} t_{a_1,b_1}\cdots t_{a_p,b_p}.
        \end{aligned}
    \end{equation*}
\end{lemma}

\begin{proof}
    By the (\ref{equation:relations for P}), it suffices to prove for the monomials 
    \begin{equation*}
        \begin{aligned}
            f=t_{k+l,r+s}^{a_{k+l,r+s}}\cdots t_{k+l,1}^{a_{k+l,1}}\cdots t_{1,r+s}^{a_{1,r+s}}\cdots t_{1,1}^{a_{1,1}}.
        \end{aligned}
    \end{equation*}
    where $a_{ij}\in \mathbb{Z}_{\geq 0 } $ whenever $[i]+[j]=\bar{0}$ and $a_{i,j}\in \{0,1\} $ whenever $ [i]+[j]=\bar{1}$.
    
    Firstly, we consider the case where  $a_{ij}=0$ for $i<k+l$. For $j>c$ , we have 
    \begin{equation*}
        \begin{aligned}
            t_{k+l,c}\partial_{k+l,c}t_{k+l,j}&=t_{k+l,c}((-1)^{([k+l]+[c])([k+l]+[j])}q_{k+l} t_{k+l,j}\partial_{k+l,c})=q_{k+l}^2 t_{k+l,j}t_{k+l,c}\partial_{k+l,c}.
        \end{aligned}
    \end{equation*}
    Then,
    \begin{equation*}
        \begin{aligned}
            &t_{k+l,c}\partial_{k+l,c}t_{k+l,r+s}^{a_{k+l,r+s}}\cdots t_{k+l,1}^{a_{k+l,1}}
            =q_{k+l}^{2\sum_{j=c+1}^{r+s}a_{k+l,j}} t_{k+l,r+s}^{a_{k+l,r+s}}\cdots t_{k+l,c+1}^{a_{k+l,c+1}} t_{k+l,c}\partial_{k+l,c}t_{k+l,c}^{a_{k+l,c}}\cdots t_{k+l,1}^{a_{k+l,1}}.
        \end{aligned}
    \end{equation*}
    It follows from Lemma \ref{lemma: the property 2 of the action of partial_ij} that
    \begin{equation*}
        \begin{aligned}
            t_{k+l,c}\partial_{k+l,c}\cdot(t_{k+l,c}^{a_{k+l,c}}\cdots t_{k+l,1}^{a_{k+l,1}})=(-1)^{[k+l]+[c]+[k+l][c]}\sum_{j=0}^{a_{k+l,c}-1}(q_{k+l}q_c(-1)^{[c]+[k+l]})^j
            t_{k+l,c}^{a_{k+l,c}}\cdots t_{k+l,1}^{a_{k+l,1}}.
        \end{aligned}
    \end{equation*}
    Since $a_{i,j}\in \{0,1\} $ whenever $ [i]+[j]=\bar{1}$, we may simplify the constant  
    \begin{equation*}
        \begin{aligned}
            \sum_{j=0}^{a_{k+l,c}-1}(q_{k+l}q_c(-1)^{[c]+[k+l]})^j &=\frac{1-(q_{k+l}q_c(-1)^{[c]+[k+l]})^{a_{k+l,c}}}{1-q_{k+l}q_c(-1)^{[c]+[k+l]}}
            =\begin{cases}
            \frac{1-q_{k+l}^{2a_{k+l,c}}}{1-q_{k+l}^2}, & [c]+[k+l]=0 ,\\
            a_{k+l,c} , & [c]+[k+l]=1.
            \end{cases}
        \end{aligned}
    \end{equation*}
    We note that $a_{k+l,c}(q_{k+l}^2-1)=q_{k+l}^{2a_{k+l,c}}-1$ when $[c]+[k+l]=1$. Then,
    \begin{equation*}
        \begin{aligned}
            (q_{k+l}^2-1)\sum_{c=1}^{r+s}(-1)^{[k+l]+[c]+[k+l][c]}t_{k+l,c}\partial_{k+l,c} \cdot f&=\sum_{c=1}^{r+s}q_{k+l}^{2\sum_{j=c+1}^{r+s}a_{k+l,j}} (q_{k+l}^{2a_{k+l,c}}-1)f\\&=(q_{k+l}^{2\sum_{j=1}^{r+s}a_{k+l,j}}-1)f.
        \end{aligned}
    \end{equation*}
    Thus, the statement holds for this special case. For  the general case, it follows from \ref{lemma: the property 3 of the action of PD on P} that 
    \begin{equation*}
        \begin{aligned}
            D_0\cdot f=(D_0\cdot (t_{k+l,r+s}^{a_{k+l,r+s}}\cdots t_{k+l,1}^{a_{k+l,1}}))h=q_{k+l}^{2\sum_{j=1}^{r+s}a_{k+l,j}}f,
        \end{aligned}
    \end{equation*}
    where $h=
t_{k+l-1,r+s}^{a_{k+l-1,r+s}}\cdots
      t_{k+l-1,1}^{a_{k+l-1,1}}\cdots
t_{1,r+s}^{a_{1,r+s}}\cdots t_{1,1}^{a_{1,1}}$. Thus, the statement holds.
\end{proof}
For any $a\in \mathbf{I}_{k|l}$, we denote $\lambda_{a}:=-2\epsilon_a-\sum_{i=a+1}^{k+l}2^{i-a}\epsilon_i$ and let $K_{\lambda_{a}}$ be the  Cartan elements corresponding $\lambda_{a}$ in $\mathbf{U}_q(\mathfrak{gl}(k|l))$. 
\begin{proposition}\label{proposition: some elements in L.}
    For any $a\in \mathbf{I}_{k|l}$, we have $\phi_U(K_{\lambda_a}\otimes1)\in \mathscr{L}_{\bullet}$ 
\end{proposition}

\begin{proof}
     By the definition of $\mathscr{L}_{\bullet}$, it suffices to show that $\phi_U(K_{\lambda_a}\otimes1)\in \mathscr{PD}^{k|l}_{r|s}$, i.e. $K_{\lambda_a}\in \mathring{U}_L$. We will prove it by the induction on $a$. For $a=k+l$, $K_{\lambda_{k+l}}=K_{k+l}^{-2}$, so we have 
     \begin{equation*}
         \begin{aligned}
             (K_{\lambda_{k+l}}\otimes1)\cdot (t_{a_1,b_1}\cdots t_{a_p,b_p})=q_{k+l}^{2\sum_{j=1}^p \delta_{k+l,a_j}}t_{a_1,b_1}\cdots t_{a_p,b_p}.
         \end{aligned}
     \end{equation*}
     Lemma~\ref{lemma: the action of the element D on products of some tij} follows that the actions of $D_0$ and $K_{\lambda_{k+l}}\otimes 1$ on $\mathscr{P}^{k|l}_{r|s}$ coincide. Then, by Proposition \ref{proposition: P is faithful PD module}, we have
     \begin{equation*}
         \begin{aligned}
             \phi_U( K_{\lambda_{k+l}}\otimes1)=D_0\in \mathscr{PD}^{k|l}_{r|s}.
         \end{aligned}
     \end{equation*}
     Assume that $a<k+l$ and $K_{\lambda_b}\in \mathring{U}_L$ for $b>a$, we want to show that $K_{\lambda_a}\in\mathring{U}_L$.
    By Lemma \ref{lemma: invariant under the adjoint action}, we can use adjoint actions to construct  elements $u_a$, $v_a$ of  $\mathring{U}_L$:
    \begin{equation*}
        \begin{aligned}
            u_a:=(1-q_{a+1}^2)^{-1}\mathrm{ad}_{E_a}(K_{\lambda_{a+1}}),\quad 
            v_a:=(1-q_{a+1}^{-2})^{-1}\mathrm{ad}_{F_{a}}(K_{\lambda_{a+1}}).
        \end{aligned}
    \end{equation*}
    Explicitly, we have $u_a=E_aK_{\lambda_{a+1}}K_a^{-1}K_{a+1}$ and $v_a=F_aK_{\lambda_{a+1}}$. Then
    \begin{equation*}
        \begin{aligned}
            [E_a,F_a]K_{2\lambda_{a+1}}K_a^{-1}K_{a+1}=q_a^{-1}q_{a+1}u_av_a-(-1)^{[a]+[a+1]}q_{a+1}^{-2}v_au_a\in \mathring{U}_L.
        \end{aligned}
    \end{equation*}
    On the other side, we have 
    \begin{equation*}
        \begin{aligned}
            [E_a,F_a]=\frac{K_aK_{a+1}^{-1}-K_a^{-1}K_{a+1}}{q_a-q_a^{-1}}.
        \end{aligned}
    \end{equation*}
    Note that $\lambda_a=2\lambda_{a+1}-2\epsilon_a+2\epsilon_{a+1}$ and $K_{2\lambda_{a+1}}=K_{\lambda_{a+1}}^2\in \mathring{U}_L$ , we have
    \begin{equation*}
        \begin{aligned}
            K_{\lambda_{a}}=K_{2\lambda_{a+1}}-(q_a-q_a^{-1})[E_a,F_a]K_{2\lambda_{a+1}}K_a^{-1}K_{a+1}\in \mathring{U}_L.
        \end{aligned}
    \end{equation*}
\end{proof}

\begin{lemma}\label{lemma: Lij belongs to L.}
    The elements $L_{i,j}\in \mathscr{L}_{\bullet}$ for $i,j\in \mathbf{I}_{k|l}$.
\end{lemma}
\begin{proof}
    By the proof of Proposition \ref{proposition: some elements in L.}, we have $D_0=\phi_U(K_{\lambda_{k+l}}\otimes1)\in \mathscr{L}_{\bullet}$. Then, $L_{k+l,k+l}=\sum\limits_{c=1}^{r+s}(-1)^{[c]([k+l]+1)}t_{k+l,c}\partial_{k+l,c}=\frac{(-1)^{[k+l]}}{q_{k+l}^2-1}(D_0-1)\in \mathscr{L}_{\bullet}$. By (\ref{equation: the relation of adjoint action and phi}), $(y\otimes1)\cdot \phi_U(x\otimes1)=\phi_{U}(\mathrm{ad}_y(x)\otimes 1)$ for any $x,y\in \mathbf{U}_q(\mathfrak{gl}(k|l))$, so $\mathscr{L}_{\bullet}$ is ($\mathbf{U}_q(\mathfrak{gl}(k|l))\otimes 1$)-stable. It follows from Lemma \ref{lemma: use Lk+l,k+l to obtain Lij} that $L_{i,j}\in \mathscr{L}_{\bullet}$.
\end{proof}

\begin{lemma}\label{lemma: ULring and Ki generate UL}
    The superalgebra $ \mathbf{U}_q(\mathfrak{gl}(k|l))$ is generated by $\mathring{U}_L$ and $\{K_i\}_{i=1}^{k+l}$.
\end{lemma}
\begin{proof}
    Let $A$ be the subalgebra of  $ \mathbf{U}_q(\mathfrak{gl}(k|l))$ generated by $\mathring{U}_L$ and $\{K_i\}_{i=1}^{k+l}$. It suffices to show that $\mathbf{U}_q(\mathfrak{gl}(k|l))\subseteq A$. Proposition \ref{proposition: some elements in L.} implies that $K_{\lambda_{a}}\in \mathring{U}_L$ for $a\in \mathbf{I}_{k|l}$, so we have $K_{\rho}=K_{\lambda_{1}}K_{\lambda_{2}}\cdots K_{\lambda_{k+l}}\in \mathring{U}_L$ where $\rho=-\sum_{i=1}^{k+l}2^i\epsilon_i$. By Lemma \ref{lemma: invariant under the adjoint action}, we have $\mathrm{ad}_{E_i}(K_{\rho})\in \mathring{U}_L$ and $\mathrm{ad}_{F_i}(K_{\rho})\in \mathring{U}_L$. Note that
    \begin{equation*}
        \begin{aligned}
            \mathrm{ad}_{E_i}(K_{\rho})=(1-q_i^{-2^i}q_{i+1}^{2^{i+1}})E_iK_{\rho}K_i^{-1}K_{i+1},\quad \mathrm{ad}_{F_i}(K_{\rho})=(1-q_i^{2^i}q_{i+1}^{-2^{i+1}})F_iK_{\rho},
        \end{aligned}
    \end{equation*}
    so the subalgebra $A$ contains $E_i,F_i$ for $1\leq i\leq k+l-1$. We have $K_i\in A$ for $i\in \mathbf{I}_{k|l}$ by definition. Then $K_i^{-1}=K_{\rho}K_{-\rho-\epsilon_i}\in A$ for $i\in \mathbf{I}_{k|l}$. Thus, $\mathbf{U}_q(\mathfrak{gl}(k|l))\subseteq A$.
\end{proof}

\begin{lemma}\label{lemma: relations between L,L. invariant subspaces}
    We have $\mathrm{End}_{\mathbb{K}}(\mathscr{P}^{k|l}_{r|s})^{\mathscr{L}_\bullet}=\mathrm{End}_{\mathbb{K}}(\mathscr{P}^{k|l}_{r|s})^{\mathscr{L}}$.
\end{lemma}
\begin{proof}
    Since $\mathscr{L}_\bullet\subseteq \mathscr{L}$, we have $\mathrm{End}_{\mathbb{K}}(\mathscr{P}^{k|l}_{r|s})^{\mathscr{L}}\subseteq\mathrm{End}_{\mathbb{K}}(\mathscr{P}^{k|l}_{r|s})^{\mathscr{L}_\bullet}$. It suffices to show that $\mathrm{End}_{\mathbb{K}}(\mathscr{P}^{k|l}_{r|s})^{\mathscr{L}_\bullet}\subseteq\mathrm{End}_{\mathbb{K}}(\mathscr{P}^{k|l}_{r|s})^{\mathscr{L}}$. For any $F\in \mathrm{End}_{\mathbb{K}}(\mathscr{P}^{k|l}_{r|s})^{\mathscr{L}_\bullet}$, we will show that $F$ supercommutes with $\phi_U(K_a\otimes 1)$ for any $a\in \mathbf{I}_{k|l}$ by induction on $a$. For $a=k+l$, $F$ supercommutes with $\phi_U(K_{-2\epsilon_{k+l}}\otimes 1)$ by Proposition \ref{proposition: some elements in L.}, so $F$ supercommutes with $\phi_U(K_{2\epsilon_{k+l}}\otimes 1)=\phi_U(K_{-2\epsilon_{k+l}}\otimes 1)^{-1}$. (\ref{equation: formulas of the action on P}) implies that the eigenspaces of $\phi_U(K_{2\epsilon_{k+l}}\otimes 1)$ and $\phi_U(K_{\epsilon_{k+l}}\otimes 1)$ are the same. Thus $F$ also supercommutes with $\phi_U(K_{\epsilon_{k+l}}\otimes 1)$. We assume that $F$ supercommutes with $\phi_U(K_b\otimes 1)$ for $b>a$. Again by  Proposition \ref{proposition: some elements in L.}, $F$ supercommutes with
    \begin{equation*}
        \begin{aligned}
            \phi_U(K_{2\epsilon_a}\otimes 1)=\phi_U(K_{\lambda_a}\otimes 1)^{-1}\prod_{i=a+1}^{k+l}\phi_U(K_i\otimes 1)^{-2^{i-a}}.
        \end{aligned}
    \end{equation*}
    Since the eigenspaces of $\phi_U(K_{2\epsilon_{a}}\otimes 1)$ and $\phi_U(K_{\epsilon_{a}}\otimes 1)$ are the same, $F$ supercommutes with $\phi_U(K_{\epsilon_{a}}\otimes 1)$. Thus, Lemma \ref{lemma: ULring and Ki generate UL} implies that $F\in \mathrm{End}_{\mathbb{K}}(\mathscr{P}^{k|l}_{r|s})^{\mathscr{L}}$. 
\end{proof}

 Thanks to Lemma \ref{lemma: relations between L,L. invariant subspaces} and Lemma \ref{lemma: relations between invariant subpaces}, we have the following result.

\begin{proposition}\label{proposition: invariant subspaces are equal}
     We have $(\mathscr{PD}^{k|l}_{r|s})^{\mathscr{L}_\bullet}=(\mathscr{PD}^{k|l}_{r|s})^{\mathscr{L}}=(\mathscr{PD}^{k|l}_{r|s})_{(\epsilon_L)}$.
\end{proposition}

\section{The Double Centralizer Property and a First Fundamental Theorem
for $\mathbf{U}_q(\mathfrak{gl}_{r|s})$}\label{Section: main results}
\subsection{Integral forms and specialization}\label{subsection: Integrals and specializations}
 Let $\mathfrak{A}:=\mathbb{C}[q,q^{-1}]_{(q-1)}$. Then $\mathfrak{A}$ is a local ring with maximal idea $\mathfrak{m}=(q-1)\mathfrak{A}$ and residue field $\mathfrak{A}/\mathfrak{m}\cong \mathbb{C}$.

We set
\begin{equation*}
\begin{aligned}
    \mathscr{A}_{p|t,p'|t',r|s}^{(d_1,d_2),gr}:=\mathscr{A}_{p|t,p'|t',r|s}^{gr}\cap \mathscr{PD}^{k|l,gr,(d_1,d_2)}_{r|s}=\bigoplus\limits_{\substack{ A\in \tilde{\Xi}_{p|t,r|s}^{k|l},\ B\in \tilde{\Xi}_{p'|t',r|s}^{k|l},\\|A|=d_1,|B|=d_2}}\mathbb{K}(t^{(A)}\otimes\partial^{(B)}).
\end{aligned}
\end{equation*}
The definition of $\mathscr{A}_{p|t,p'|t',r|s}^{(d_1,d_2)}$ is similar.
We define the following free $\mathfrak{A}$-submodules
by taking the $\mathfrak{A}$-spans of the PBW basis:
\begin{equation}\label{equation: direct sum decomposition of the integral form of A}
    \begin{aligned}
        &\mathscr{A}_{p|t,p'|t',r|s,\mathfrak{A}}^{(d_1,d_2),gr}:=\bigoplus\limits_{\substack{ A\in \tilde{\Xi}_{p|t,r|s}^{k|l},\ B\in \tilde{\Xi}_{p'|t',r|s}^{k|l},\\|A|=d_1,|B|=d_2}}\mathfrak{A}(t^{(A)}\otimes\partial^{(B)})\subseteq \mathscr{A}_{p|t,p'|t',r|s}^{(d_1,d_2),gr},\\
        &\mathscr{A}_{p|t,p'|t',r|s,\mathfrak{A}}^{gr}:=\bigoplus_{d_1,d_2\in\mathbb{N}}\mathscr{A}_{p|t,p'|t',r|s,\mathfrak{A}}^{(d_1,d_2),gr}\subseteq\mathscr{A}_{p|t,p'|t',r|s}^{gr}.\\
        &\mathscr{PD}^{k|l,gr,(d_1,d_2)}_{r|s,\mathfrak{A}}:=\bigoplus\limits_{\substack{A,B\in \Xi_{k|l,r|s},\\|A|=d_1,|B|=d_2}}\mathfrak{A}(t^{(A)}\otimes\partial^{(B)})\subseteq  \mathscr{PD}^{k|l,gr,(d_1,d_2)}_{r|s}.\\
        &\mathscr{PD}^{k|l,gr}_{r|s,\mathfrak{A}}:=\bigoplus_{d_1,d_2\in\mathbb{N}}\mathscr{PD}^{k|l,gr,(d_1,d_2)}_{r|s,\mathfrak{A}}.
    \end{aligned}
\end{equation}
Obviously, we have $\mathbb{K}\otimes_{\mathfrak{A}}\mathscr{A}_{p|t,p'|t',r|s,\mathfrak{A}}^{(d_1,d_2),gr}\cong\mathscr{A}_{p|t,p'|t',r|s}^{(d_1,d_2),gr}$ and $\mathbb{K}\otimes_{\mathfrak{A}}\mathscr{PD}^{k|l,gr,(d_1,d_2)}_{r|s,\mathfrak{A}}\cong\mathscr{PD}^{k|l,gr,(d_1,d_2)}_{r|s}$ as superalgebras.
\begin{lemma}\label{lemma: Aspans of PBW basis is Asuperalgebras}
    The $\mathfrak{A}$-modules $\mathscr{A}_{p|t,p'|t',r|s,\mathfrak{A}}^{gr}$, $\mathscr{PD}^{k|l,gr}_{r|s,\mathfrak{A}}$ are $\mathfrak{A}$-superalgebras.
\end{lemma}
\begin{proof}
   For any $ A,A'\in \tilde{\Xi}_{p|t,r|s}^{k|l},\ B,B'\in \tilde{\Xi}_{p'|t',r|s}^{k|l}$, by using relations (\ref{equation:relations for P}), (\ref{equation:relations for D}), (\ref{equation: the first relation for PD})-(\ref{equation: the fourth relation for PD}), we can write $(t^{(A)}\otimes\partial^{(B)})( t^{(A')}\otimes\partial^{(B')})=\sum a_{C,D}t^{(C)}\otimes\partial^{(D)}$, $a_{C,D}\in \mathfrak{A}$. Hence,  $\mathscr{A}_{p|t,p'|t',r|s,\mathfrak{A}}^{gr}$ is an $\mathfrak{A}$-superalgebras. The proof for $\mathscr{PD}^{k|l,gr}_{r|s,\mathfrak{A}}$ is similar.
\end{proof}

The specialization of the integral form $\mathscr{A}_{p|t,p'|t',r|s,\mathfrak{A}}^{(d_1,d_2),gr}$ and $\mathscr{A}_{p|t,p'|t',r|s,\mathfrak{A}}^{gr}$ at $q=1$ are
\begin{equation*}
\begin{aligned}
&\mathscr{A}_{p|t,p'|t',r|s}^{(d_1,d_2),gr}|_{q=1}:=\mathbb{C}_1\otimes_{\mathfrak{A}}\mathscr{A}_{p|t,p'|t',r|s,\mathfrak{A}}^{(d_1,d_2),gr} \cong \mathscr{A}_{p|t,p'|t',r|s,\mathfrak{A}}^{(d_1,d_2),gr}/\mathfrak{m}\mathscr{A}_{p|t,p'|t',r|s,\mathfrak{A}}^{(d_1,d_2),gr},\\
&\mathscr{A}_{p|t,p'|t',r|s}^{gr}|_{q=1}:=\mathbb{C}_1\otimes_{\mathfrak{A}}\mathscr{A}_{p|t,p'|t',r|s,\mathfrak{A}}^{,gr} \cong \bigoplus_{d_1,d_2\in\mathbb{N}}\mathscr{A}_{p|t,p'|t',r|s}^{(d_1,d_2),gr}|_{q=1},
\end{aligned}
\end{equation*}
where $\mathbb{C}_1=\mathbb{C}$ with $q$ acts as $1$ on $\mathbb{C}$. Then, we have
\begin{equation}\label{equation: direct sum decomposition of A at q=1}
    \begin{aligned}
        \mathscr{A}_{p|t,p'|t',r|s}^{gr}|_{q=1}\cong \mathscr{A}_{p|t,p'|t',r|s,\mathfrak{A}}^{gr}/(q-1)\mathscr{A}_{p|t,p'|t',r|s,\mathfrak{A}}^{gr}\cong \bigoplus\limits_{\substack{ A\in \tilde{\Xi}_{p|t,r|s}^{k|l},\ B\in \tilde{\Xi}_{p'|t',r|s}^{k|l}}}\mathbb{C}(\overline{t^{(A)}\otimes\partial^{(B)}}),
    \end{aligned}
\end{equation}
which implies that $\overline{t^{(A)}\otimes\partial^{(B)}}$,$A\in \tilde{\Xi}_{p|t,r|s}^{k|l},\ B\in \tilde{\Xi}_{p'|t',r|s}^{k|l}$ form a $\mathbb{C}$-basis of $\mathscr{A}_{p|t,p'|t',r|s}^{gr}|_{q=1}$. Let $\mathcal{T}$ be the superalgebra generated by $t_{i,j}$,$\partial_{i',j}$,$i\in \tilde{\mathbf{I}}_{p|t,k|l}$,$i'\in\tilde{\mathbf{I}}_{p'|t',k|l}$,$j\in \mathbf{I}_{r|s}$ modulo the relations:
\begin{equation}\label{equation: relation of A at q=1}
    \begin{aligned}
        &t_{a,b}^{2}=0,\quad\text{if }[a]+[b]=\bar{1},\quad\partial_{c',d}^{2}=0 \quad\text{if }[c']+[d]=\bar{1},\\
        &t_{a,b}t_{c,d}
        =
        (-1)^{([a]+[b])([c]+[d])}t_{c,d}t_{a,b},
        \\
        &\partial_{a',b}\partial_{c',d}
        =
        (-1)^{([a']+[b])([c']+[d])}\partial_{c',d}\partial_{a',b},
        \\
        &\partial_{c',d}t_{a,b}
        =
        (-1)^{([a]+[b])([c']+[d])}t_{a,b}\partial_{c',d},
    \end{aligned}   
\end{equation}
where $a,c\in \tilde{\mathbf{I}}_{p|t,k|l}$,$a',c'\in\tilde{\mathbf{I}}_{p'|t',k|l}$,$b,d\in\mathbf{I}_{r|s}$.
Then there exists an isomorphism of superalgebras:
\begin{equation*}
    \begin{aligned}
         \vartheta_1:\mathcal{T}&\xrightarrow[]{\sim} \mathscr{A}_{p|t,p'|t',r|s}^{gr}|_{q=1},\quad
        t_{i,j}\mapsto\overline{t_{i,j}\otimes1} ,\quad
        \partial_{i',j}\mapsto \overline{1\otimes\partial_{i',j}}.
    \end{aligned}
\end{equation*}

Similarly, the specialization of the integral form $\mathscr{PD}^{k|l,gr}_{r|s,\mathfrak{A}}$ at $q=1$ is
\begin{equation*}
    \begin{aligned}
\mathscr{PD}^{k|l,gr}_{r|s}|_{q=1}:=\mathbb{C}_1\otimes_{\mathfrak{A}}\mathscr{PD}^{k|l,gr}_{r|s,\mathfrak{A}}
    \end{aligned}
\end{equation*}
where $\mathbb{C}_1=\mathbb{C}$ with $q$ acts as $1$ on $\mathbb{C}$. Then $\mathscr{PD}^{k|l,gr}_{r|s}|_{q=1}$ is isomorphic to the superalgebra generated by $t_{a,b},\partial_{a,b}$, $a\in\mathbf{I}_{k|l}$, $b\in \mathbf{I}_{r|s}$ modulo the relations:
\begin{equation*}
    \begin{aligned}
        \partial_{a,b} t_{c,d}
        - (-1)^{([a]+[b])([c]+[d])} t_{c,d} \partial_{a,b}
        = (-1)^{[a]+[b]+[a][b]} \delta_{a,c}\delta_{b,d}.
    \end{aligned}
\end{equation*}
Set
$
D_{a,b}=(-1)^{[a]+[b]+[a][b]}\partial_{a,b},
$ $ a\in I_{k|l}$,$\ b\in I_{r|s}.
$, the defining relations become
\begin{equation*}
    \begin{aligned}
    D_{a,b} t_{c,d}
            - (-1)^{([a]+[b])([c]+[d])} t_{c,d}D_{a,b}
            =\delta_{a,c}\delta_{b,d}.
    \end{aligned}
\end{equation*}
where $a,c\in I_{k|l}$ and $b,d\in I_{r|s}$. Thus $\mathscr{PD}^{k|l,gr}_{r|s}|_{q=1}$  is isomorphic to Weyl-Clifford algebra.

Let $\tilde{k}:=\mathrm{max}\{p'+t'-l,0\}$, $\tilde{l}:=\mathrm{min}\{p'+t',l\}$, $\tilde{r}:=\mathrm{max}\{p+t-l,0\}$, $\tilde{s}:=\mathrm{min}\{p+t,l\}$.Then $\tilde{k}$ and $\tilde{l}$ are the numbers of even and odd
indices in $\mathbf{I}_{p'|t',k|l}$, respectively, while
$\tilde{r}$ and $\tilde{s}$ are the numbers of even and odd
indices in $\mathbf{I}_{p|t,k|l}$, respectively. In pariticular, $\tilde{k}+\tilde{l}=p'+t'$, $\tilde{r}+\tilde{s}=p+t$. Let 
\begin{equation*}
    \begin{aligned}
        \mathcal{P}^{\tilde{k}|\tilde{l}}_{\tilde{r}|\tilde{s}}|_{q=1}\cong S(
\mathbb C^{\tilde{k}|\tilde{l}}\otimes\mathbb C^{r|s}
\oplus
\mathbb C^{\tilde{r}|\tilde{s}}\otimes(\mathbb C^{r|s})^*
)
    \end{aligned}
\end{equation*}
 be the superalgebra defined in \cite[Section 6.1]{Zh20}. Let $\mathbf{U}(\mathfrak{gl}(r|s))$ be the universal enveloping
algebra of the general linear Lie superalgebra $\mathfrak{gl}(r|s)$.
It is a Hopf superalgebra with comultiplication, counit, and antipode
given by
\begin{equation*}
    \begin{aligned}
        \Delta(x)=1\otimes x+x\otimes 1,\quad \epsilon(x)=0,\quad
        S(x)=-x,\quad \forall x\in\mathfrak{gl}(r|s).
    \end{aligned}
\end{equation*}
We equip $\mathcal{P}^{\tilde{k}|\tilde{l}}_{\tilde{r}|\tilde{s}}|_{q=1}$ with the $\mathbf{U}(\mathfrak{gl}(r|s))$-module structure given in \cite[Section 6.1]{Zh20}. $\mathcal{P}^{\tilde{k}|\tilde{l}}_{\tilde{r}|\tilde{s}}|_{q=1}$ is generated by $T_{i,j},\overline{T}_{i',j}$,$i\in \mathbf{I}_{\tilde{k}|\tilde{l}}$, $i'\in \mathbf{I}_{\tilde{r}|\tilde{s}}$, $j\in \mathbf{I}_{r|s} $ modulo the following relations:
\begin{equation}\label{equation: relation of Pzy at q=1}
    \begin{aligned}
        &\overline{T}_{a',b}^{2}=0,\quad\text{if }[a']+[b]=\bar{1},\quad T_{c,d}^{2}=0 \quad\text{if }[c]+[d]=\bar{1},\\
        &\overline{T}_{a',b}\overline{T}_{c',d}
        =
        (-1)^{([a']+[b])([c']+[d])}\overline{T}_{c',d}\overline{T}_{a',b},
        \\
        &T_{a,b}T_{c,d}
        =
        (-1)^{([a]+[b])([c]+[d])}T_{c,d}T_{a,b},
        \\
        &T_{c,d}\overline{T}_{a',b}
        =
        (-1)^{([a']+[b])([c]+[d])}\overline{T}_{a',b}T_{c,d},
    \end{aligned}   
\end{equation}
where $a,c\in \mathbf{I}_{\tilde{k}|\tilde{l}}$,$a',c'\in\mathbf{I}_{\tilde{r}|\tilde{s}}$,$b,d\in\mathbf{I}_{r|s}$.

Thus, by the relations
\eqref{equation: relation of A at q=1} and
\eqref{equation: relation of Pzy at q=1},
we can verify that there is a superalgebra isomorphism given by
\begin{equation*}
    \begin{aligned}
\vartheta_2:\mathcal{P}^{\tilde{k}|\tilde{l}}_{\tilde{r}|\tilde{s}}|_{q=1}\xrightarrow[]{\sim}\mathcal{T},\quad
        \overline{T}_{i',j}\mapsto t_{k+l-p-t+i',j},\quad
        T_{i,j}\mapsto \partial_{k+l-p'-t'+i,j}.
    \end{aligned}
\end{equation*}
The parameters $\tilde{k},\tilde{l},\tilde{r},\tilde{s}$
are chosen so that these assignments preserve the
$\mathbb Z_2$-grading. We obtain a superalgebra isomorphism $\vartheta:=\vartheta_1\circ\vartheta_2:\mathcal{P}^{\tilde{k}|\tilde{l}}_{\tilde{r}|\tilde{s}}|_{q=1}\xrightarrow[]{\sim}\mathscr{A}_{p|t,p'|t',r|s}^{gr}|_{q=1}$. Thus, we induce a $\mathbf{U}(\mathfrak{gl}(r|s))$-module structure on
$\mathscr{A}_{p|t,p'|t',r|s}^{gr}|_{q=1}$ via $\vartheta$.
Explicitly,for $x\in \mathbf{U}(\mathfrak{gl}(r|s))$,
$f\in \mathscr{A}_{p|t,p'|t',r|s}^{gr}|_{q=1}$, $x\cdot f=\vartheta(x\cdot \vartheta^{-1}(f))$. With this action, $\vartheta$ is a
$\mathbf{U}(\mathfrak{gl}(r|s))$-module isomorphism.
Let
\begin{equation*}
    \begin{aligned}
(\mathscr{A}_{p|t,p'|t',r|s}^{gr}|_{q=1})^{\mathbf{U}(\mathfrak{gl}(r|s))}:=\{f\in \mathscr{A}_{p|t,p'|t',r|s}^{gr}|_{q=1}\mid x\cdot f=\epsilon(x)f \ \mathrm{for\ all}\ x\in \mathbf{U}(\mathfrak{gl}(r|s))\}
    \end{aligned}
\end{equation*}
denote its $\mathbf{U}(\mathfrak{gl}(r|s))$-invariant subalgebra.

\begin{lemma}\label{lemma: generators of the invariant subalg of A at q=1}
   The $\mathbf{U}(\mathfrak{gl}(r|s))$-invariant subalgebra $(\mathscr{A}_{p|t,p'|t',r|s}^{gr}|_{q=1})^{\mathbf{U}(\mathfrak{gl}(r|s))}$ is generated by $\overline{{L}_{i,j}^{gr}}$, $i\in \tilde{\mathbf{I}}_{p|t,k|l}$, $j\in \tilde{\mathbf{I}}_{p'|t',k|l}$. In particular, 
   \begin{equation*}
       \begin{aligned}
           (\mathscr{A}_{p|t,p'|t',r|s}^{(d,d),gr}|_{q=1})^{\mathbf{U}(\mathfrak{gl}(r|s))}=\mathrm{span_{\mathbb{C}}}\{\overline{L_{i_1,j_1}^{gr}\cdots L_{i_d,j_d}^{gr}},i_k\in  \tilde{\mathbf{I}}_{p|t,k|l},j_k\in \tilde{\mathbf{I}}_{p'|t',k|l} \}
       \end{aligned}
   \end{equation*}

\end{lemma}
\begin{proof}
    Let $\mathcal{X}^{\tilde{k}|\tilde{l}}_{\tilde{r}|\tilde{s}}|_{q=1}$ be the  $\mathbf{U}(\mathfrak{gl}(r|s))$-invariant subalgebra of $\mathcal{P}^{\tilde{k}|\tilde{l}}_{\tilde{r}|\tilde{s}}|_{q=1}$. It follows from \cite[Theorem 6.2]{Zh20} that
    $\mathcal{X}^{\tilde{k}|\tilde{l}}_{\tilde{r}|\tilde{s}}|_{q=1}$
    is generated by the following elements:
    \begin{equation*}
        \begin{aligned}
            X_{a,b}=\sum_{c\in\mathbf{I}_{r|s} }(-1)^{[a]([b]+[c])}T_{a,c}\overline{T}_{b,c}, \quad a\in \mathbf{I}_{\tilde{k}|\tilde{l}},b\in \mathbf{I}_{\tilde{r}|\tilde{s}}.
        \end{aligned}
    \end{equation*}
    Since $\vartheta$ is a $\mathbf{U}(\mathfrak{gl}(r|s))$-module isomorphism,
    we have $\vartheta( \mathcal{X}^{\tilde{k}|\tilde{l}}_{\tilde{r}|\tilde{s}}|_{q=1})=(\mathscr{A}_{p|t,p'|t',r|s}^{gr}|_{q=1})^{\mathbf{U}(\mathfrak{gl}(r|s))}$.
    A direct calculation gives
    \begin{equation*}
        \begin{aligned}
            \vartheta(X_{a,b})&=\sum_{c\in \mathbf{I}_{r|s}}(-1)^{[a]([b]+[c])}\overline{(1\otimes \partial_{k+l-p'-t'+a,c})}\cdot\overline{t_{k+l-p-t+b}\otimes1}\\
            &=\sum_{c\in \mathbf{I}_{r|s}}(-1)^{[j]([i]+[c])+([i]+[c])([j]+[c])}\overline{(t_{i,c}\otimes \partial_{j,c})}\\
            &=\sum_{c\in \mathbf{I}_{r|s}}(-1)^{[c]([i]+1)}\overline{(t_{i,c}\otimes \partial_{j,c})}=\overline{L_{i,j}},
        \end{aligned}
    \end{equation*}
    where $i=k+l-p-t+b\in \tilde{\mathbf{I}}_{p|t,k|l}$, $j=k+l-p'-t'+a\in  \tilde{\mathbf{I}}_{p'|t',k|l}$. Since $\vartheta$ is a superalgebra isomorphism, $(\mathscr{A}_{p|t,p'|t',r|s}^{gr}|_{q=1})^{\mathbf{U}(\mathfrak{gl}(r|s))}$ is generated by $\bar{L}_{i,j}^{gr}$,$i\in \tilde{\mathbf{I}}_{p|t,k|l}$, $j\in \tilde{\mathbf{I}}_{p'|t',k|l}$. The relations $\overline{L_{i,j}}\in\mathscr{A}_{p|t,p'|t',r|s}^{(1,1),gr}|_{q=1} $ and  $\mathscr{A}_{p|t,p'|t',r|s}^{(d_1,d_1),gr}|_{q=1}\cdot \mathscr{A}_{p|t,p'|t',r|s}^{(d_2,d_2),gr}|_{q=1}\subseteq\mathscr{A}_{p|t,p'|t',r|s}^{(d_1+d_2,d_1+d_2),gr}|_{q=1}$ implies that $(\mathscr{A}_{p|t,p'|t',r|s}^{(d,d),gr}|_{q=1})^{\mathbf{U}(\mathfrak{gl}(r|s))}=\mathrm{span_{\mathbb{C}}}\{\overline{L_{i_1,j_1}^{gr}\cdots L_{i_d,j_d}^{gr}},i_k\in  \tilde{\mathbf{I}}_{p|t,k|l},j_k\in \tilde{\mathbf{I}}_{p'|t',k|l} \}$.
\end{proof}

From the $\mathbf{U}(\mathfrak{gl}(r|s))$-module superalgebra structure on
$\mathcal{P}^{\tilde{k}|\tilde{l}}_{\tilde{r}|\tilde{s}}|_{q=1}$
given in \cite{Zh20}, we obtain the following explicit action formulas
on the generators:
\begin{equation*}
    \begin{aligned}
        e_{a,a+1}\cdot(\overline{ t_{i,j}\otimes 1})
&=-\delta_{j,a}\ \overline {t_{i,a+1}\otimes 1},\\
e_{a+1,a}\cdot(\overline {t_{i,j}\otimes 1})
&=(-1)^{[a]+[a+1]+1}\delta_{j,a+1}\
\overline {t_{i,a}\otimes 1},\\
e_{b,b}\cdot(\overline{ t_{i,j}\otimes 1})
&=-\delta_{b,j}\ \overline{ t_{i,j}\otimes 1},\\
e_{a,a+1}\cdot(\overline{1\otimes\partial_{i,j}})
&=\delta_{j,a+1}\ \overline{1\otimes\partial_{i,a}},\\
e_{a+1,a}\cdot(\overline{1\otimes\partial_{i,j}})
&=\delta_{j,a}\ \overline{1\otimes\partial_{i,a+1}},\\
e_{b,b}\cdot(\overline{1\otimes\partial_{i,j}})
&=\delta_{b,j}\ \overline{1\otimes\partial_{i,j}}.
    \end{aligned}
\end{equation*}
where $a\in \mathbf{I}_{r|s}',b\in \mathbf{I}_{r|s}$. It follows from (\ref{equation: formulas of the action on D}) and
(\ref{equation: formulas of the action on P}) that $\mathscr{A}_{p|t,p'|t',r|s,\mathfrak{A}}^{gr} $ is stable under the action of $1\otimes E_a,1\otimes F_a,1\otimes\frac{K_b-1}{q_b-1}\in1\otimes\mathbf{U}_q(\mathfrak{gl}(r|s)) $.
Moreover , since $\mathscr{A}_{p|t,p'|t',r|s}^{gr}$
is a $1\otimes\mathbf{U}_q(\mathfrak{gl}(r|s))$-module superalgebra
and $\mathscr{A}_{p|t,p'|t',r|s}^{gr}|_{q=1}$
is a $\mathbf{U}(\mathfrak{gl}(r|s))$-module superalgebra,
 for every
$f\in\mathscr{A}_{p|t,p'|t',r|s,\mathfrak{A}}^{gr} $, we have
\begin{equation}\label{equation: compatibility of the module structures under specialization}
    \begin{aligned}
        e_{a,a+1}\cdot\overline f
        =\overline{E_a\cdot f},\quad
        e_{a+1,a}\cdot\overline f
        =\overline{F_a\cdot f},\quad
        e_{b,b}\cdot\overline f
        =\overline{\frac{K_b-1}{q_b-1}\cdot f},
    \end{aligned}
\end{equation}
where $a\in \mathbf{I}_{r|s}'$, $b\in \mathbf{I}_{r|s}$.
In the above formulas, we abbreviate $1\otimes x$ as $x$
whenever $x\in \mathbf{U}_q(\mathfrak{gl}(r|s))$.

\begin{lemma}\label{lemma: compare the dim of invariant subspaces} 
    For any $d\geq 0$, we have $\mathrm{dim}_{\mathbb{K}}(\mathscr{A}_{p|t,p'|t',r|s}^{(d,d),gr})_{(\epsilon_R)}\leq \mathrm{dim}_{\mathbb{C}}(\mathscr{A}_{p|t,p'|t',r|s}^{(d,d),gr}|_{q=1})^{\mathbf{U}(\mathfrak{gl}(r|s))}$.
\end{lemma}
\begin{proof}
    If $r+s=0$, there are no generators $t_{i,j}$ or
    $\partial_{i',j}$, so
    $\mathscr A^{\mathrm{gr}}_{p|t,p'|t',r|s}=\mathbb K$.
    Both dimensions in the statement are therefore equal to
    $1$ for $d=0$ and to $0$ for $d>0$.
    Hence we may assume that $r+s>0$.
    \begin{equation*}
        \begin{aligned}
            \mathcal{Q}_d:\mathscr{A}_{p|t,p'|t',r|s}^{(d,d),gr}&\xrightarrow[]{}(\mathscr{A}_{p|t,p'|t',r|s}^{(d,d),gr})^{\oplus(3(r+s)-2)},\\
            f&\mapsto((E_a\cdot f)_{a=1}^{r+s-1},(F_a\cdot f)_{a=1}^{r+s-1},(\frac{K_b-1}{q_b-1}\cdot f)_{b=1}^{r+s}).
        \end{aligned}
    \end{equation*}
    Then we have $(\mathscr{A}_{p|t,p'|t',r|s}^{(d,d),gr})_{(\epsilon_R)}=\mathrm{Ker}(\mathcal{Q}_d)$. By definition, the 
$\{t^{(A)}\otimes\partial^{(B)},A\in \tilde{\Xi}_{p|t,r|s}^{k|l},\ B\in \tilde{\Xi}_{p'|t',r|s}^{k|l},|A|=|B|=d\}$
forms a $\mathbb{K}$-basis of
$\mathscr{A}_{p|t,p'|t',r|s}^{(d,d),gr}$.
This basis is finite, and we denote its cardinality by $\mathcal{N}_d$.
We fix an ordering of the basis and write it as
$h_1,\cdots,h_{\mathcal{N}_d}$. Then
    \begin{equation*}
        \begin{aligned}
            h_1^{(1)},\cdots,h_{\mathcal{N}_d}^{(1)},\cdots,h_1^{3(r+s)-2},\cdots,h_{\mathcal{N}_d}^{3(r+s)-2}
        \end{aligned}
    \end{equation*}
    form a $\mathbb{K}$-basis of
    $(\mathscr{A}_{p|t,p'|t',r|s}^{(d,d),gr})^{\oplus(3(r+s)-2)}$,
    where $h_j^{(i)}$ is the element of
    $(\mathscr{A}_{p|t,p'|t',r|s}^{(d,d),gr})^{\oplus(3(r+s)-2)}$
    whose the $i$th component is $h_j$ and the other components are $0$.

    Let $M(q)$ be the matrix of $\mathcal{Q}_d$ with respect to these two bases. It follows from (\ref{equation: formulas of the action on D}),(\ref{equation: formulas of the action on P}) that $M(q)\in \mathrm{Mat}_{(3(r+s)-2)\mathcal{N}_d\times\mathcal{N}_d}(\mathfrak{A})$.
    By Rank–nullity theorem, we have
    \begin{equation}\label{equation: relation of dimensions}
        \begin{aligned}
            \mathrm{dim}_{\mathbb{K}}(\mathscr{A}_{p|t,p'|t',r|s}^{(d,d),gr})_{(\epsilon_R)}=\mathrm{dim}_{\mathbb{K}}{\mathrm{Ker}(\mathcal{Q}_d)}=\mathrm{dim}_{\mathbb{K}}(\mathscr{A}_{p|t,p'|t',r|s}^{(d,d),gr})-\mathrm{rank}_{\mathbb{K}}(M(q)).
        \end{aligned}
    \end{equation}
    Restricting $\mathcal{Q}_d$ to
    $\mathscr{A}_{p|t,p'|t',r|s,\mathfrak{A}}^{(d,d),gr}$, we obtain
        \begin{equation*}
        \begin{aligned}
        \mathcal{Q}_{d,\mathfrak{A}}:\mathscr{A}_{p|t,p'|t',r|s,\mathfrak{A}}^{(d,d),gr}&\xrightarrow[]{}(\mathscr{A}_{p|t,p'|t',r|s,\mathfrak{A}}^{(d,d),gr})^{\oplus(3(r+s)-2)},\\
            f&\mapsto((E_a\cdot f)_{a=1}^{r+s-1},(F_a\cdot f)_{a=1}^{r+s-1},(\frac{K_b-1}{q_b-1}\cdot f)_{b=1}^{r+s}).
        \end{aligned}
    \end{equation*}
By the definition of $\mathscr{A}_{p|t,p'|t',r|s,\mathfrak{A}}^{(d,d),gr}$,
$h_1,\cdots,h_{\mathcal{N}_d}$ form an $\mathfrak{A}$-basis of
$\mathscr{A}_{p|t,p'|t',r|s,\mathfrak{A}}^{(d,d),gr}$. Hence,
    \begin{equation*}
        \begin{aligned}
            h_1^{(1)},\cdots,h_{\mathcal{N}_d}^{(1)},\cdots,h_1^{3(r+s)-2},\cdots,h_{\mathcal{N}_d}^{3(r+s)-2}
        \end{aligned}
    \end{equation*}
    form an $\mathfrak{A}$-basis of $\mathscr{A}_{p|t,p'|t',r|s,\mathfrak{A}}^{(d,d),gr})^{\oplus(3(r+s)-2)}$.
   Then, the matrix of $\mathcal{Q}_{d,\mathfrak{A}}$ with respect to these two bases is also $M(q)$.
    We induce a linear map via $\mathcal{Q}_d$ given by
    \begin{equation*}
        \begin{aligned}
        \tilde{\mathcal{Q}}_d:\mathscr{A}_{p|t,p'|t',r|s}^{(d,d),gr}|_{q=1}&\xrightarrow[]{}(\mathscr{A}_{p|t,p'|t',r|s,\mathfrak{A}}^{(d,d),gr})^{\oplus(3(r+s)-2)}/\mathfrak{m}(\mathscr{A}_{p|t,p'|t',r|s,\mathfrak{A}}^{(d,d),gr})^{\oplus(3(r+s)-2)},\\
            f&\mapsto \overline{\mathcal{Q}_d(f)}.
        \end{aligned}
    \end{equation*}
    By the direct sum decompositions
    (\ref{equation: direct sum decomposition of the integral form of A}) and
    (\ref{equation: direct sum decomposition of A at q=1}),
    we obtain a linear isomorphism
    \begin{equation*}
        \begin{aligned}
            (\mathscr{A}_{p|t,p'|t',r|s,\mathfrak{A}}^{(d,d),gr})^{\oplus(3(r+s)-2)}/\mathfrak{m}(\mathscr{A}_{p|t,p'|t',r|s,\mathfrak{A}}^{(d,d),gr})^{\oplus(3(r+s)-2)}&\cong(\mathscr{A}_{p|t,p'|t',r|s}^{(d,d),gr}|_{q=1})^{\oplus(3(r+s)-2)},\\
            \overline{(x_1,\cdots,x_{3(r+s)-2})}&\mapsto(\overline{x_1},\cdots,\overline{x_{3(r+s)-2}}).
        \end{aligned}
    \end{equation*}
 Via this isomorphism, we identify the elements of these two vector spaces.
Thus, by
    (\ref{equation: compatibility of the module structures under specialization}), we can rewrite the map as
    \begin{equation*}
        \begin{aligned}
        \tilde{\mathcal{Q}}_d:\mathscr{A}_{p|t,p'|t',r|s}^{(d,d),gr}|_{q=1}&\xrightarrow[]{}(\mathscr{A}_{p|t,p'|t',r|s}^{(d,d),gr}|_{q=1})^{\oplus(3(r+s)-2)},\\
            f&\mapsto((e_{a,a+1}\cdot\overline{f})_{a=1}^{r+s-1},(e_{a+1,a}\cdot\overline{f})_{a=1}^{r+s-1},(e_{b,b}\cdot \overline{ f})_{b=1}^{r+s}).
        \end{aligned}
    \end{equation*}
    Thus, $(\mathscr{A}_{p|t,p'|t',r|s}^{(d,d),gr}|_{q=1})^{\mathbf{U}(\mathfrak{gl}(r|s))}=\mathrm{Ker}(\tilde{\mathcal{Q}}_d)$.Since $\mathscr{A}_{p|t,p'|t',r|s}^{(d,d),gr}|_{q=1}$
    has the direct sum decomposition
    (\ref{equation: direct sum decomposition of A at q=1}),
    $\{\overline{h_1},\cdots, \overline{h_{\mathcal{N}_d}}\}$
    forms a $\mathbb{C}$-basis of
    $\mathscr{A}_{p|t,p'|t',r|s}^{(d,d),gr}|_{q=1}$.
    Hence,
    \begin{equation*}
        \begin{aligned}
            \overline{h_1^{(1)}},\cdots,\overline{h_{\mathcal{N}_d}^{(1)}},\cdots,\overline{h_1^{3(r+s)-2}},\cdots,\overline{h_{\mathcal{N}_d}^{3(r+s)-2}}
        \end{aligned}
    \end{equation*}
    forms a $\mathbb{C}$-basis of
    $(\mathscr{A}_{p|t,p'|t',r|s}^{(d,d),gr}|_{q=1})^{\oplus(3(r+s)-2)}$.
    By the definition of $\tilde{\mathcal{Q}}_d$, we have
    \begin{equation*}
        \begin{aligned}
            \tilde{\mathcal{Q}}_d(\overline{h_v})=\overline{\mathcal{Q}_d(h_v)}=\overline{\sum\limits_{\substack{1\leq i\leq 3(r+s)-2\\1\leq j\leq \mathcal{N}_d}}m_{(i-1)\mathcal{N}_d+j,v}(q)h_{j}^{(i)}}=\sum\limits_{\substack{1\leq i\leq 3(r+s)-2\\1\leq j\leq \mathcal{N}_d}}m_{(i-1)\mathcal{N}_d+j,v}(1)\overline{h_{j}^{(i)}}.
        \end{aligned}
    \end{equation*}
Therefore, the matrix of $\tilde{\mathcal{Q}}_d$ with respect to these two bases is $M(1)$. Hence
\begin{equation*}
        \begin{aligned}
            \mathrm{dim}_{\mathbb{C}}(\mathscr{A}_{p|t,p'|t',r|s}^{(d,d),gr}|_{q=1})^{\mathbf{U}(\mathfrak{gl}(r|s))}=\mathrm{dim}_{\mathbb{C}}{\mathrm{Ker}(\tilde{\mathcal{Q}}_d)}&=\mathrm{dim}_{\mathbb{C}}(\mathscr{A}_{p|t,p'|t',r|s}^{(d,d),gr}|_{q=1})-\mathrm{rank_{\mathbb{C}}(M(1))}\\
            &=\mathrm{dim}_{\mathbb{K}}(\mathscr{A}_{p|t,p'|t',r|s}^{(d,d),gr})-\mathrm{rank_{\mathbb{C}}(M(1))}.
        \end{aligned}
    \end{equation*}
Note that $\mathbb{K}=\mathbb{C}(q)$ and
$\mathrm{rank_{\mathbb{C}}(M(1))}\leq \mathrm{rank}_{\mathbb{C}(q)}(M(q))$,
we complete the proof of the lemma by comparing the above equality with
(\ref{equation: relation of dimensions}).
\end{proof}
\begin{proposition}\label{proposition: the generaotors of invariant alg of A}
    The invariant subspace 
    \begin{equation*}
        \begin{aligned}
            (\mathscr{A}_{p|t,p'|t',r|s}^{(d,d),gr})_{(\epsilon_R)}=\mathrm{span_{\mathbb{K}}}\{L_{i_1,j_1}^{gr}\cdots L_{i_d,j_d}^{gr},i_k\in  \tilde{\mathbf{I}}_{p|t,k|l},j_k\in \tilde{\mathbf{I}}_{p'|t',k|l} \}
        \end{aligned}
    \end{equation*}
     and the invariant subalgebra $(\mathscr{A}_{p|t,p'|t',r|s}^{gr})_{(\epsilon_R)}$ is generated by $L_{i,j}^{gr}$, $i\in\tilde{\mathbf{I}}_{p|t,k|l} $, $j\in \tilde{\mathbf{I}}_{p'|t',k|l}$.
\end{proposition}
\begin{proof}
    Proposition \ref{proposition: Lijgr Lij are invariants} implies that 
    \begin{equation*}
        \begin{aligned}
            \mathrm{span_{\mathbb{K}}}\{L_{i_1,j_1}^{gr}\cdots L_{i_d,j_d}^{gr},i_k\in  \tilde{\mathbf{I}}_{p|t,k|l},j_k\in \tilde{\mathbf{I}}_{p'|t',k|l} \}\subseteq (\mathscr{A}_{p|t,p'|t',r|s}^{(d,d),gr})_{(\epsilon_R)}.
        \end{aligned}
    \end{equation*}
 It follows from (\ref{lemma: generators of the invariant subalg of A at q=1}),(\ref{lemma: compare the dim of invariant subspaces}) that 
    \begin{equation*}
        \begin{aligned}
            \dim_{\mathbb{C}} (\mathscr{A}_{p|t,p'|t',r|s}^{(d,d),gr}|_{q=1})^{\mathbf{U}(\mathfrak{gl}(r|s))}
            &= \dim_{\mathbb{C}} \mathrm{span_{\mathbb{C}}}(\{\overline{L_{i_1,j_1}^{gr}\cdots L_{i_d,j_d}^{gr}},i_k\in  \tilde{\mathbf{I}}_{p|t,k|l},j_k\in \tilde{\mathbf{I}}_{p'|t',k|l} \})\\
            &\leq\mathrm{dim}_{\mathbb{K}}( \mathrm{span_{\mathbb{K}}}\{L_{i_1,j_1}^{gr}\cdots L_{i_d,j_d}^{gr},i_k\in  \tilde{\mathbf{I}}_{p|t,k|l},j_k\in \tilde{\mathbf{I}}_{p'|t',k|l} \})\\
            &\leq \mathrm{dim}_{\mathbb{K}}(\mathscr{A}_{p|t,p'|t',r|s}^{(d,d),gr})_{(\epsilon_R)}\leq \dim_{\mathbb{C}} (\mathscr{A}_{p|t,p'|t',r|s}^{(d,d),gr}|_{q=1})^{\mathbf{U}(\mathfrak{gl}(r|s))}.
        \end{aligned}
    \end{equation*}
    Therefore, all the above inequalities are equalities. In particular, we have 
    \begin{equation*}
        \begin{aligned}
            \mathrm{dim}_{\mathbb{K}}(\mathscr{A}_{p|t,p'|t',r|s}^{(d,d),gr})_{(\epsilon_R)}=\mathrm{dim}_{\mathbb{K}}( \mathrm{span_{\mathbb{K}}}\{L_{i_1,j_1}^{gr}\cdots L_{i_d,j_d}^{gr},i_k\in  \tilde{\mathbf{I}}_{p|t,k|l},j_k\in \tilde{\mathbf{I}}_{p'|t',k|l} \}).
        \end{aligned}
    \end{equation*}
    This proves the first assertion of the proposition. We set $K=\prod_{b=1}^{r+s}K_b^{(-1)^{[b]}}$. Then, it follows from (\ref{equation: formulas of the action on D}),(\ref{equation: formulas of the action on P}) that
    \begin{equation*}
        \begin{aligned}
            K\cdot t_{i,j}=q^{-1} t_{i,j},\quad K\cdot \partial_{i',j}=q\partial_{i',j},
        \end{aligned}
    \end{equation*}
    where $i\in \tilde{\mathbf{I}}_{p|t,k|l}$, $i'\in\tilde{\mathbf{I}}_{p'|t',k|l}$, $j\in \mathbf{\mathbf{I}}_{r|s}$. Thus, for any $f\in \mathscr{A}_{p|t,p'|t',r|s}^{(d_1,d_2),gr}$, we have 
  $K\cdot f=q^{d_2-d_1}f$.
    Moreover, if $f\in (\mathscr{A}_{p|t,p'|t',r|s}^{(d_1,d_2),gr})_{(\epsilon_R)}$, then $q^{d_2-d_1}f=K\cdot f=f$.
    It implies that
    \begin{equation*}
        \begin{aligned}
            (\mathscr{A}_{p|t,p'|t',r|s}^{(d_1,d_2),gr})_{(\epsilon_R)}=0 \quad \mathrm{for} \ d_1\neq d_2. 
        \end{aligned}
    \end{equation*} 
    Hence, we have
    \begin{equation*}
        \begin{aligned}
            (\mathscr{A}_{p|t,p'|t',r|s}^{gr})_{(\epsilon_R)}&=\bigoplus_{d\geq 0}(\mathscr{A}_{p|t,p'|t',r|s}^{(d,d),gr})_{(\epsilon_R)}
            \\
            &=\bigoplus_{d\geq 0}\mathrm{span_{\mathbb{K}}}\{L_{i_1,j_1}^{gr}\cdots L_{i_d,j_d}^{gr},i_k\in  \tilde{\mathbf{I}}_{p|t,k|l},j_k\in \tilde{\mathbf{I}}_{p'|t',k|l} \}.
        \end{aligned}
    \end{equation*}
\end{proof}

\subsection{main results}

Recall the $\mathbf{U}_q(\mathfrak{gl}(k|l))\otimes\mathbf{U}_q(\mathfrak{gl}(r|s))$-module isomorphism $P:\mathscr{PD}^{k|l,gr}_{r|s}\xrightarrow[]{}\mathscr{PD}^{k|l}_{r|s},~f\otimes g\mapsto fg $ defined in (\ref{equation: the module isomorphism between PDgr and PD}), we restrict it  to $\mathscr{A}_{p|t,p'|t',r|s}^{gr}$, then we obtain a $1\otimes\mathbf{U}_q(\mathfrak{gl}(r|s))$-module isomorphism
\begin{equation*}
    \begin{aligned}
        P:\mathscr{A}_{p|t,p'|t',r|s}^{gr}\xrightarrow[]{\sim}\mathscr{A}_{p|t,p'|t',r|s}.
    \end{aligned}
\end{equation*}

\begin{lemma}\label{lemma: P(D)P(D')-P(DD') belongs to what}
    For $f\otimes g\in \mathscr{PD}^{k|l,gr,(d,d)}_{r|s}$ and $f'\otimes g'\in \mathscr{PD}^{k|l,gr,(d',d')}_{r|s}$, we have
    \begin{equation*}
        \begin{aligned}
            P(f\otimes g)P(f'\otimes g')-P((f\otimes g)(f'\otimes g'))\in \bigoplus_{i=1}^{d+d'}\mathscr{PD}^{k|l,(d+d'-i,d+d'-i)}_{r|s}.
        \end{aligned}
    \end{equation*}
    In particular, for $f\otimes g\in \mathscr{A}_{p|t,p'|t',r|s}^{(d,d),gr}$ and $f'\otimes g'\in  \mathscr{A}_{p|t,p'|t',r|s}^{(d',d'),gr}$, we have
    \begin{equation*}
        \begin{aligned}
            P(f\otimes g)P(f'\otimes g')-P((f\otimes g)(f'\otimes g'))\in \bigoplus_{i=1}^{d+d'} \mathscr{A}_{p|t,p'|t',r|s}^{(d+d'-i,d+d'-i)}.
        \end{aligned}
    \end{equation*}
\end{lemma}
\begin{proof}
     Proposition \ref{proposition: the relations for PD} implies that the  relations of $\mathscr{PD}^{k|l,gr}_{r|s}$ differ from those of $\mathscr{PD}^{k|l}_{r|s}$ only by removing the constant term $(-1)^{[a]+[b]+[a][b]}$ on the right-hand side of (\ref{equation: the fourth relation for PD}). Using this fact and the same argument as in the proof of \cite[Proposition 3.13.3]{LSS24}, the result follows.
\end{proof}

\begin{lemma}\cite[Proposition 9.3.2]{LSS24}\label{lemma: transfer generators from a graded algebra to a filtered algebra}
Let $A=\bigoplus_{i\ge0} A^{(i)}$ be a graded algebra and let $B$ be a filtered algebra. For $i\geq 0$, set
$A^i := \bigoplus_{j=0}^{i} A^{(j)}$, then it induces a filtration on $A$. 
Suppose that
$F:A\xrightarrow[]{} B$
is a filtration-preserving linear map such that the induced map
$\mathrm{gr}(F):\mathrm{gr}(A)\to \mathrm{gr}(B)$
is an isomorphism of algebras. If $a_1,\dots,a_r\in A^{(1)}$ generate $A$, then 
$F(a_1),\dots,F(a_r)$ generate $B$.
\end{lemma}
The above lemma  allows us to transfer generators from a graded algebra to a filtered algebra via an 
isomorphism between their associated graded algebras. We have already known the generators of $(\mathscr{A}_{p|t,p'|t',r|s}^{gr})_{(\epsilon_R)}$ by Proposition \ref{proposition: the generaotors of invariant alg of A}, we will use this lemma to obtain the generators of $(\mathscr{A}_{p|t,p'|t',r|s})_{(\epsilon_R)}$.

\begin{theorem}\label{theorem: AeR is generated by  Lij}
    The superalgebra $(\mathscr{A}_{p|t,p'|t',r|s})_{(\epsilon_R)}$ is generated by $L_{i,j}$, $i\in\tilde{\mathbf{I}}_{p|t,k|l} $, $j\in \tilde{\mathbf{I}}_{p'|t',k|l}$.
\end{theorem}
\begin{proof}
    Set  $A:=(\mathscr{A}_{p|t,p'|t',r|s}^{gr})_{(\epsilon_R)}$, $B:=(\mathscr{A}_{p|t,p'|t',r|s})_{(\epsilon_R)}$.  we first define a grading on $A$ and a filtration of $B$. We set 
\begin{equation*}
    \begin{aligned}
        A^{(d)}:=\mathscr{A}_{p|t,p'|t',r|s}^{gr,(d,d)}\cap A=(\mathscr{A}_{p|t,p'|t',r|s}^{gr,(d,d)})_{(\epsilon_R)}\quad \mathrm{for} ~d\geq0.
    \end{aligned}
\end{equation*}
Then, by proof of Proposition \ref{proposition: the generaotors of invariant alg of A}, we have 
 $A=\bigoplus_{d\geq 0}A^{(d)}$. 
Set
\begin{equation*}
    \begin{aligned}
        A^i:=\bigoplus_{d=0}^{i}A^{(d)},\quad B^{i}:=P(A^i).
    \end{aligned}
\end{equation*}
Since $P$ is a  $\mathbf{U}_q(\mathfrak{gl}(k|l))\otimes\mathbf{U}_q(\mathfrak{gl}(r|s))$-module isomorphism, we have $P(A)=B$ and $B^{0}\subseteq B^{1}\subseteq B^2\subseteq\cdots$ is a filtration of $B$ as a vector space. Moreover, the map $\mathrm{gr}(P):\mathrm{gr}(A)\to \mathrm{gr}(B)$ induced by $P$
is an isomorphism of graded vector space. Proposition \ref{proposition: the generaotors of invariant alg of A} implies that $L_{i,j}^{gr}\in A^{(1)}$ for $i\in\tilde{\mathbf{I}}_{p|t,k|l} $, $j\in \tilde{\mathbf{I}}_{p'|t',k|l}$ generate $A$. Then, by Lemma \ref{lemma: transfer generators from a graded algebra to a filtered algebra} and the relations $P(L_{i,j}^{gr})=L_{i,j}$ for $i\in\tilde{\mathbf{I}}_{p|t,k|l} $, $j\in \tilde{\mathbf{I}}_{p'|t',k|l}$, it suffices to show that $\mathrm{gr}(P)$ is an homomorphism of algebras, which is equivalent to $P(D)P(D')-P(DD')\in B^{d+d'-1}$ for $D\in A^{(d)}$, $D'\in A^{(d')}$. We set $\tilde{B}^{(d)}:=P( \mathscr{A}_{p|t,p'|t',r|s}^{gr,(d,d)})$ and $\tilde{B}^{i}:=\bigoplus_{d=0}^i\tilde{B}^{(d)}$. By Lemma \ref{lemma: P(D)P(D')-P(DD') belongs to what} and the relation $P(A)=B$, we have 
\begin{equation*}
    \begin{aligned}
        P(D)P(D')-P(DD')\in \tilde{B}^{d+d'-1}\cap B\quad \mathrm{for}\  D\in A^{(d)},D'\in A^{(d')}.
    \end{aligned}
\end{equation*}
Note that $P$ is a bijection, so we have 
\begin{equation*}
    \begin{aligned}
        \tilde{B}^{i}\cap B=(\bigoplus_{d=0}^{i}P( \mathscr{A}_{p|t,p'|t',r|s}^{gr,(d,d)}))\cap P(A)=P((\bigoplus_{d=0}^{i} \mathscr{A}_{p|t,p'|t',r|s}^{gr,(d,d)}) \cap A)=P(A^{i})=B^i .
    \end{aligned}
\end{equation*}
Thus, we have $P(D)P(D')-P(DD')\in B^{d+d'-1}$ and then $\mathrm{gr}(P)$ is an homomorphism of algebras.
\end{proof}

\begin{corollary}\label{corollary: PDeR is generated by Lij}
     The superalgebra $(\mathscr{PD}^{k|l}_{r|s})_{(\epsilon_R)}$ is generated by $L_{i,j}$ for $i,j\in \mathbf{I}_{k|l}$.
\end{corollary}
\begin{proof}
The assertion follows immediately by setting $p=p'=k$ and $t=t'=l$ in Theorem \ref{theorem: AeR is generated by  Lij}.
\end{proof}

By (\ref{equation: the first relation for PD})-(\ref{equation: the fourth relation for PD}) in Proposition \ref{proposition: the relations for PD}, we can check that there exists an isomorphism of superalgebras given by 
\begin{equation*}
    \begin{aligned}
        \theta^{k|l}_{r|s}:\mathscr{PD}^{k|l}_{r|s}&\xrightarrow[]{}\mathscr{PD}^{r|s}_{k|l},\\
        t_{i,j}\mapsto(-1)^{[i][j]}&t_{j,i},\quad \partial_{i,j}\mapsto(-1)^{[i][j]}\partial_{j,i}.
    \end{aligned}
\end{equation*}
By viewing $\mathscr{P}^{k|l}_{r|s}$ as a subsuperalgebra  of $\mathscr{PD}^{k|l}_{r|s}$, $\theta^{k|l}_{r|s}$ restricts  to an isomorphism between $\mathscr{P}^{k|l}_{r|s}$ and $\mathscr{P}^{r|s}_{k|l}$, which induces an isomorphism
\begin{equation*}
    \begin{aligned}
        \tilde{\theta}^{k|l}_{r|s}:\mathrm{End}_{\mathbb{K}}(\mathscr{P}^{k|l}_{r|s})&\xrightarrow[]{}\mathrm{End}_{\mathbb{K}}(\mathscr{P}^{r|s}_{k|l}),\\
        F&\xrightarrow[]{}(\theta^{r|s}_{k|l})^{-1}\circ F\circ \theta^{r|s}_{k|l}.
    \end{aligned}
\end{equation*}
Note that $(\theta^{r|s}_{k|l})^{-1}= \theta^{k|l}_{r|s}$, then we have $\tilde{\theta}^{k|l}_{r|s}(D)(f)=(\theta^{r|s}_{k|l})^{-1}(D\cdot \theta^{r|s}_{k|l}(f))=\theta^{k|l}_{r|s}( D)(f)$ for  $f\in \mathscr{P}^{r|s}_{k|l}$, $D\in \mathscr{PD}^{k|l}_{r|s}$. Thus ,$\tilde{\theta}^{k|l}_{r|s}|_{\mathscr{PD}^{k|l}_{r|s}}=\theta^{k|l}_{r|s}$.

\begin{lemma}\label{lemma: how (1) impies (2)}
    The isomorphism $\theta^{k|l}_{r|s}$ has the following properties.
    
   \item[(1)]
    $\theta^{k|l}_{r|s}(L_{i,j})=R_{i,j}$ for $i,j\in \mathbf{I}_{k|l}$.
    
    \item[(2)]
    $\theta^{k|l}_{r|s}((x\otimes y)\cdot D)=(-1)^{[x][y]}(y\otimes x)\theta^{k|l}_{r|s}(D)$ for homogeneous $x\otimes y \in \mathbf{U}_q(\mathfrak{gl}(k|l))\otimes\mathbf{U}_q(\mathfrak{gl}(r|s))$, $D\in \mathscr{PD}^{k|l}_{r|s}$.

    \item[(3)] $\theta^{k|l}_{r|s}((\mathscr{PD}^{k|l}_{r|s})_{(\epsilon_R)})=(\mathscr{PD}^{r|s}_{k|l})_{(\epsilon_L)}$.

    \item[(4)] $\theta^{k|l}_{r|s}(\mathscr{L_{\bullet}})=\mathscr{R_\bullet}$, $\theta^{k|l}_{r|s}(\mathscr{R_{\bullet}})=\mathscr{L_\bullet}$.

    \item[(5)] $\theta^{k|l}_{r|s}((\mathscr{PD}^{k|l}_{r|s})^{\mathscr{R}_\bullet})=(\mathscr{PD}^{r|s}_{k|l})^{\mathscr{L}_\bullet}$.
\end{lemma}
\begin{proof}
    (1) follows from the definition of $\theta^{k|l}_{r|s}$ and direct calculations. Since $\mathscr{PD}^{k|l}_{r|s}$ is a  $ \mathbf{U}_q(\mathfrak{gl}(k|l))\otimes\mathbf{U}_q(\mathfrak{gl}(r|s))$-module superalgebra, it suffices to check for generators $D=t_{i,j},\partial_{i,j}$. For generators $x\otimes y$ of the form $x\otimes 1$ or $1\otimes y$ where $x,y=E_a,F_a,K_a$, (2) holds by
    the (\ref{equation: formulas of the action on D}) and (\ref{equation: formulas of the action on P}). Thus (2) holds for any $x\otimes y\in \mathbf{U}_q(\mathfrak{gl}(k|l))\otimes\mathbf{U}_q(\mathfrak{gl}(r|s))$.
    (3) follows from (2) and the bijection of $\theta^{k|l}_{r|s}$.

    It suffices to show that $\theta^{k|l}_{r|s}(\mathscr{L_{\bullet}})=\mathscr{R_\bullet}$ in (4), since $(\theta^{k|l}_{r|s})^{-1}=\theta^{r|s}_{k|l}$. It follows from (2) that for any $x\in \mathbf{U}_q(\mathfrak{gl}(k|l))$, $f\in \mathscr{P}^{r|s}_{k|l}$, we have
    \begin{equation*}
        \begin{aligned}
            \tilde{\theta}^{k|l}_{r|s}(\phi_U(x\otimes 1))(f)
            =(\theta^{r|s}_{k|l})^{-1}( (x\otimes 1)\cdot\theta^{r|s}_{k|l}(f))
            =(\theta^{r|s}_{k|l})^{-1}(\theta^{r|s}_{k|l}((1\otimes x)\cdot f))=\phi_U(1\otimes x)(f).
        \end{aligned}
    \end{equation*}
     Thus, $\tilde{\theta}^{k|l}_{r|s}(\phi_U(\mathbf{U}_q(\mathfrak{gl}(k|l))\otimes 1))=\phi_U(1\otimes \mathbf{U}_q(\mathfrak{gl}(k|l)))$. Then for any $f\in \mathscr{L}_{\bullet}=\phi_U(\mathbf{U}_q(\mathfrak{gl}(k|l))\otimes 1)\cap\mathscr{PD}^{k|l}_{r|s}$, we have
     \begin{equation*}
         \begin{aligned}
             \theta^{k|l}_{r|s}(f)=\tilde{\theta}^{k|l}_{r|s}(f)\in \phi_U(1\otimes\mathbf{U}_q(\mathfrak{gl}(k|l)))\cap\mathscr{PD}^{r|s}_{k|l}=\mathscr{R}_{\bullet}.
         \end{aligned}
     \end{equation*}
     Thus, $\theta^{k|l}_{r|s}(\mathscr{L_{\bullet}})\subseteq\mathscr{R_\bullet}$. Fix any $g\in \mathscr{R}_{\bullet}$,  since $\theta^{k|l}_{r|s}$ is surjective, $\theta^{k|l}_{r|s}(f)=g$ for some $f\in \mathscr{PD}^{k|l}_{r|s}$. Moreover, there exists $f'\in \phi_U(\mathbf{U}_q(\mathfrak{gl}(k|l))\otimes 1)$ such that $\tilde{\theta}^{k|l}_{r|s}(f')=g$. Since $\tilde{\theta}^{k|l}_{r|s}$ is bijective and $\tilde{\theta}^{k|l}_{r|s}|_{\mathscr{PD}^{k|l}_{r|s}}=\theta^{k|l}_{r|s}$, we have $f=f'\in  \phi_U(\mathbf{U}_q(\mathfrak{gl}(k|l))\otimes 1)\cap  \mathscr{PD}^{k|l}_{r|s}=\mathscr{L}_{\bullet}$. Thus, $\mathscr{R_\bullet}\subseteq\theta^{k|l}_{r|s}(\mathscr{L_{\bullet}})$.
     (5) follows from (4) immediately.
\end{proof}

\begin{theorem}\label{theorem: The double centralizer property}
     Recall $\mathscr{L_\bullet}$, $\mathscr{R}_\bullet$ in (\ref{definition: L. and R.}) and $L_{i,j}$, $R_{i,j}$ in (\ref{equation: definition of L_i,j}). Let $\langle L_{i,j}\rangle$ be the subalgebra of $\mathscr{PD}^{k|l}_{r|s}$ generated by $L_{i,j}$ for $i,j\in \mathbf{I}_{k|l}$ and $\langle R_{i,j}\rangle$ be the subalgebra of $\mathscr{PD}^{k|l}_{r|s}$ generated by $R_{i,j}$ for $i,j\in \mathbf{I}_{r|s}$ ,   we have
    
    \item[(1)]$(\mathscr{PD}^{k|l}_{r|s})^{\mathscr{R}_\bullet}=(\mathscr{PD}^{k|l}_{r|s})_{(\epsilon_R)}=\mathscr{L_\bullet}=\langle L_{i,j}\rangle$.

    \item[(2)]$(\mathscr{PD}^{k|l}_{r|s})^{\mathscr{L}_\bullet}=(\mathscr{PD}^{k|l}_{r|s})_{(\epsilon_L)}=\mathscr{R_\bullet}=\langle R_{i,j}\rangle.$
\end{theorem}

\begin{proof}
      Firstly, we show that $\mathscr{L}_{\bullet}\subseteq (\mathscr{PD}^{k|l}_{r|s})_{(\epsilon_R)}$. For any  $D\in \mathscr{L}_{\bullet}$, there exists $x\in\mathbf{U}_q(\mathfrak{gl}(k|l)) $ such that $D=\phi_U(x\otimes 1)$ by  the definition  of $\mathscr{L}_{\bullet}$. (\ref{equation: the relation of adjoint action and phi}) implies that for any $y\in \mathbf{U}_q(\mathfrak{gl}(r|s))$,  we have 
    \begin{equation*}
        \begin{aligned}
            (1\otimes y)\cdot D=\phi_{U}(x\otimes \mathrm{ad}_{y}(1))=\epsilon(y)\phi_U(x\otimes 1)=\epsilon(y)D.
        \end{aligned}
    \end{equation*}
    Thus, $D\in (\mathscr{PD}^{k|l}_{r|s})_{(\epsilon_R)}$. By Lemma \ref{lemma: Lij belongs to L.} , we have $\langle L_{i,j}\rangle\subseteq \mathscr{L}_{\bullet}$. Then Corollary \ref{corollary: PDeR is generated by Lij}  implies that $ (\mathscr{PD}^{k|l}_{r|s})_{(\epsilon_R)}= \mathscr{L}_{\bullet}=\langle L_{i,j}\rangle$ for $i,j\in \mathbf{I}_{k|l}$. Thus, (1) follows from Proposition \ref{proposition: invariant subspaces are equal} and (2) follows from (1) by Lemma \ref{lemma: how (1) impies (2)}.  
\end{proof}

\end{document}